\pdfoutput=1
\documentclass[3p,nopreprintline]{elsarticle}
\usepackage[T1]{fontenc}
\usepackage{amsmath,amssymb,amsthm,mathtools,bm}
\usepackage{booktabs,array,tabularx,longtable,multirow,makecell}
\usepackage{graphicx,subcaption,adjustbox}
\usepackage{siunitx,microtype}
\usepackage[normalem]{ulem}
\usepackage{xcolor}
\usepackage{aliascnt}
\usepackage[colorlinks,linkcolor=magenta!80!black,citecolor=blue!80!black,urlcolor=blue!80!black]{hyperref}
\usepackage[capitalise,noabbrev]{cleveref}
\usepackage{placeins}
\usepackage{algorithm}
\usepackage[noend]{algpseudocode}
\newcommand{\bR}{\mathbb R}

\newcommand{\HH}{\mathbb H}
\newcommand{\LL}{\mathbb L}
\newcommand{\PP}{\mathfrak P}
\newcommand{\pp}{\mathfrak p}
\newcommand{\U}{\mathbb U}
\newcommand{\X}{\mathbb X}
\newcommand{\Y}{\mathbb Y}
\newcommand{\cB}{\mathcal B}
\newcommand{\cK}{\mathcal K}
\newcommand{\cR}{\mathcal R}
\newcommand{\cL}{\mathcal L}
\newcommand{\cN}{\mathcal N}
\newcommand{\cT}{\mathcal T}
\newcommand{\cS}{\mathcal S}

\newcommand{\cA}{\mathcal A}
\newcommand{\bss}{\boldsymbol{s}}
\newcommand{\bst}{\boldsymbol{t}}
\newcommand{\bsg}{\boldsymbol{g}}
\newcommand{\bsv}{\boldsymbol{v}}
\newcommand{\uc}{u^\circ}
\newcommand{\sigc}{\sigma^\circ}
\newcommand{\norm}[1]{\left\|#1\right\|}

\newcommand{\rdiv}{\operatorname{div}}
\newcommand{\Hdiv}{H(\operatorname{div};\Omega)}
\newcommand{\HdivN}{H_{0,\Gamma_N}(\operatorname{div};\Omega)}
\newcommand{\clin}{c_{\mathrm{lin}}}
\newcommand{\Clin}{C_{\mathrm{lin}}}
\newcommand{\Crel}{C_{\mathrm{rel}}}
\newcommand{\Ceff}{C_{\mathrm{eff}}}
\newcommand{\MK}{M_{\cK}}

\newcommand{\pmin}{\mathfrak p_{\min}}
\newcommand{\pmax}{\mathfrak p_{\max}}
\newcommand{\csec}{c_{\mathrm{sec}}}

\newcommand{\bmnum}[1]{\bfseries #1}%
\newcommand{\uwnum}[1]{{\uwave{\num{#1}}}}%

\newtheorem{theorem}{Theorem}[section]
\newcommand{\sharedthm}[3]{\newaliascnt{#1}{theorem}\newtheorem{#1}[#1]{#2}%
  \aliascntresetthe{#1}\crefname{#1}{#2}{#3}\Crefname{#1}{#2}{#3}}
\sharedthm{lemma}{Lemma}{Lemmas}
\sharedthm{proposition}{Proposition}{Propositions}
\sharedthm{corollary}{Corollary}{Corollaries}
\newtheorem{assumption}{Assumption}
\theoremstyle{remark}
\sharedthm{remark}{Remark}{Remarks}
\crefname{assumption}{Assumption}{Assumptions}
\Crefname{assumption}{Assumption}{Assumptions}
\crefname{algorithm}{Algorithm}{Algorithms}
\renewcommand{\arraystretch}{1.12}

\hypersetup{pdftitle={Variationally correct physics-informed operator learning for nonlinear partial differential equations},pdfauthor={Yuan Qiu; Wolfgang Dahmen; Peng Chen}}
\begin{document}
\begin{frontmatter}
\title{Variationally correct physics-informed operator learning for nonlinear partial differential equations}
\author[1]{Yuan Qiu}\ead{yuan.qiu@gatech.edu}
\author[2]{Wolfgang Dahmen}\ead{dahmen@math.sc.edu}
\author[1]{Peng Chen\corref{cor1}}\ead{pchen402@gatech.edu}
\affiliation[1]{organization={School of Computational Science and Engineering},addressline={Georgia Institute of Technology},city={Atlanta},postcode={GA 30332},country={USA}}
\affiliation[2]{organization={Department of Mathematics},addressline={University of South Carolina},city={Columbia},postcode={SC 29208},country={USA}}
\cortext[cor1]{Corresponding author.}
\begin{abstract}
Physics-informed operator learning trains neural operators on the residual of a parametric partial differential equation (PDE) instead of labeled solution data, but its reliability depends on a variationally correct residual: one that is equivalent to the error in a norm dictated by the equation, the PDE-compliant norm. For nonlinear equations written as first-order systems, we prove that the first-order system least-squares residual is equivalent to the error in this norm near a solution branch if the linearized residual is stable and the residual derivative is Lipschitz continuous, possibly in a norm stronger than the error norm. We verify these conditions for monotone reaction, solution- and gradient-dependent diffusion, and compressible Neo-Hookean hyperelasticity. The stability mechanisms vary across these classes, as do the regions where the equivalence holds, from arbitrary bounded sets to a neighborhood of the branch. In a conforming finite element space, we show that the loss excess of a state near the finite element solution is equivalent to its squared distance from that solution. We use this loss to train a reduced basis neural operator with values in this space, with solution snapshots needed only for the basis and the loss and its symbolically derived gradient assembled by the finite element library. At inference, this loss estimates the error. Numerically, the square root of the loss excess tracks the error within factors comparable to those of the linear case, and the operator achieves higher accuracy in the PDE-compliant norm than grid-based baselines evaluated in the same finite element space.
\end{abstract}
\end{frontmatter}
\section{Introduction}\label{sec:introduction}
Uncertainty quantification, inverse problems, and design optimization require the solution of a parametric partial differential equation (PDE) for many parameter values. For nonlinear PDEs, each solve typically requires an iterative procedure, making solution evaluation across many parameter values computationally expensive. Neural operators such as the deep operator network (DeepONet)~\cite{lu2021learning} and the Fourier neural operator (FNO)~\cite{li2020fourier} provide fast surrogates for these solves by learning the parameter-to-solution map. Data-driven neural operators learn from solution snapshots, shifting the computational cost to an offline training stage that may itself require many PDE solves. Their computational advantage over classical solvers depends on the problem structure and target accuracy, with neural operators being most competitive at low-to-moderate accuracy requirements, even when data-generation and training costs are fully amortized~\cite{huangStuart2026costaccuracy}. Alternatively, physics-informed operator learning trains on the PDE residual and avoids the need for solution snapshots, since the residual of a prediction can be evaluated at any parameter value. It extends physics-informed neural networks~\cite{raissi2019pinn} and their variational formulations~\cite{kharazmi2021hp} from a single PDE instance to the parameter-to-solution map~\cite{wang2021learning,li2024physics}. The usefulness of this approach hinges on two questions:
\begin{quote}
\itshape
Does a small residual imply a small error, and in which norm?\\
How can the residual be evaluated and differentiated accurately and efficiently enough to serve as a training loss?
\end{quote}

For linear problems, the answer to the first question rests on a stable variational formulation, in which errors in the trial norm are equivalent to residuals in the dual test norm. For several PDE models, a first-order formulation is stable with a test space of $L^2$ type, which is self-dual. If $\HH$ is the state space, $s^*(\pp)$ the solution for the parameter $\pp$, and $\cB_\pp s-S_\pp$ the residual in an $L^2$-type product space $\LL^2$, then
\begin{equation}\label{eq:linear_vc}
 c_0\norm{s-s^*(\pp)}_\HH^2
 \le \norm{\cB_\pp s-S_\pp}_{\LL^2}^2
 \le C_0\norm{s-s^*(\pp)}_\HH^2 ,
\end{equation}
so that the first-order system least-squares (FOSLS) residual loss is equivalent to the error in a norm dictated by the PDE~\cite{caiLazarovManteuffelMcCormick1994,caiManteuffelMcCormick1997,caiStarke2003,bochevGunzburger2009}. Losses with this property are called variationally correct~\cite{BDO,CDG2025}. In~\cite{QCD} we combined a FOSLS loss, variational lifts of mixed boundary data, and a reduced basis neural operator (RBNO) that predicts reduced basis coefficients in a conforming finite element space for the full first-order state. Fluxes and stresses are direct outputs of the model rather than quantities recovered by postprocessing the predicted primal field. Each prediction can be evaluated with the same residual functional as the finite element solution, and solution snapshots are needed only to construct the reduced basis. In contrast to data-driven reduced models~\cite{bhattacharya2021model,dal2020data,fresca2021comprehensive}, the network is trained by minimizing the residual loss, which provides an estimate of the prediction error.

The choice of residual decides whether a physics-informed loss controls the error. A pointwise strong-form residual can fail to distinguish the intended PDE solution from solutions of a different problem, allowing a neural network to achieve a small training loss while approximating the wrong solution. In the first benchmark of \cref{sec:numerics}, we consider $-\rdiv(\pp\nabla u)+\kappa u^3=f$, with $f=f_2+\rdiv f_1$, where the conductivity $\pp$ and the field $f_1$ are piecewise constant and jump across the boundaries of sixteen square inclusions. Away from these interfaces, their derivatives vanish, so pointwise evaluation of the second-order residual, as commonly used in physics-informed neural networks, yields only $-\pp\Delta u+\kappa u^3-f_2$. This residual does not capture the interface contribution of $f_1$ or enforce the required transmission condition. It therefore vanishes not only at the intended solution, but also at the solution of the problem with $f_1$ omitted, and more generally at any sufficiently regular function satisfying the equation within each subdomain and the external boundary conditions, without necessarily satisfying the transmission condition. The first-order formulation addresses this ambiguity by requiring the total flux $\sigma=\pp\nabla u+f_1$ to belong to $\Hdiv$, which ensures continuity of its normal component across the interfaces, as discussed in \cref{rem:interfaces}. Moreover, \cref{prop:semilinear} shows that the resulting FOSLS residual controls the $H^1$ error of the primal component $u$.

In this paper, we extend the variationally correct physics-informed operator learning framework developed in~\cite{QCD} from linear to nonlinear PDEs. This extension requires addressing two main issues. First, a residual--error equivalence such as~\eqref{eq:linear_vc} generally holds only locally. Nonlinear equations may admit multiple solutions, their linearizations may become singular along a solution branch, and constitutive laws may restrict the admissible states. Second, establishing this equivalence requires assumptions that depend on the form of the nonlinearity. For reaction terms, which enter the lower-order part of the operator, monotonicity provides the needed control. For solution- and gradient-dependent diffusion, the nonlinearity enters the principal part, and Lipschitz continuity of the linearization generally requires a stronger norm than the one used to measure the error. In finite-strain elasticity, the stress law is nonpolynomial and defined only for orientation-preserving deformations, while the linearized operator may lose stability as the load varies. To accommodate these different structures, we express the residual as the sum of a linear first-order operator and a nonlinear term acting on the primal variable. We then identify two conditions on the linearization that ensure local residual--error equivalence and verify these conditions for each of the three problem classes.

Local stability analysis has its roots in the approximation theory of nonsingular solution branches~\cite{brezziRappazRaviart1980}, and residual-based error estimates are well established for nonlinear reduced basis methods~\cite{veroyPrudhommePatera2003,quarteroni2015reduced,haasdonk2017tutorial}. In our setting, the loss measures the residual of the full first-order state and serves as the training objective for a parameter-to-solution map. The residual--error equivalence must therefore hold uniformly over the parameter set and remain valid under reduced basis approximation and numerical integration. Related questions for individual PDE instances have been studied through neural least-squares methods~\cite{bersetche2023deep,opschoor2024fosls}, nonlinear least-squares finite element methods~\cite{bertrand2025lsfem}, and a posteriori error bounds for physics-informed neural networks~\cite{ernst2025posteriori}. For operator learning, recent work has connected residual-based adaptivity to variational objectives~\cite{toscano2025adaptivity} and trained reduced basis neural operators using physics-informed residuals~\cite{zheng2025rebano}.

The second question concerns how to evaluate and differentiate the residual efficiently during training. In the linear case, the loss of a reduced state is a quadratic form with matrices that can be assembled once~\cite{QCD}. For nonlinear problems, this structure is generally unavailable, so we retain finite element assembly in the training loop. For each network prediction, the finite element code evaluates the loss and assembles its gradient with respect to the finite element coefficients (or degrees of freedom). The corresponding derivative form is obtained by symbolic differentiation and compiled once by the finite element library. This gradient is then propagated to the network parameters through the chain rule~\cite{bouziani2023physics,mitusch2019dolfin}, without requiring a forward or adjoint solve. Galerkin residuals have previously been used as data-free training losses for operator networks~\cite{lee2025feonet}, with an error analysis for linear elliptic problems involving the condition number of the finite element matrix~\cite{hong2024error}. In our formulation, the loss is the squared $L^2$ norm of the first-order residual evaluated at a conforming state. Under the local stability assumptions, it is equivalent to the squared error in the PDE-compliant norm, with constants independent of the mesh.

\subsection{Contributions}
The main contributions are as follows:
\begin{enumerate}
\item \emph{Variational correctness for nonlinear first-order systems.}
In \cref{thm:local_vc}, we establish local equivalence between the FOSLS residual and the solution error under stability and Lipschitz continuity assumptions on the linearization. The error is measured in the PDE-compliant norm of the first-order system, while the neighborhood of the solution branch may be defined in a stronger norm to accommodate gradient-dependent nonlinearities. We further account for errors in residual evaluation through an additive term in \cref{prop:consistency}. Within this neighborhood, the residual therefore provides an assessment of prediction accuracy without requiring a reference solution.

\item \emph{Verification for three classes of nonlinearity.}
We verify these assumptions for monotone reaction, nonlinear diffusion with solution-dependent coefficients of Kirchhoff type or gradient-dependent coefficients of $p$-Laplace type, and compressible Neo-Hookean hyperelasticity, with the results summarized in \cref{tab:example_summary}. For monotone reaction and Kirchhoff-type diffusion, residual--error equivalence holds on every bounded set of states. For the reaction problem, we further obtain an error bound computable from the problem data and the prediction. For the hyperelasticity problem, the equivalence holds locally near the solution branch and relies on the stability of the tangent operator, which we assess numerically.

\item \emph{Physics-informed training through finite element residuals.}
In \cref{lem:loss_excess}, we show that, for conforming finite element states near the finite element solution, the loss excess relative to that solution is equivalent to the squared distance from it, allowing the equivalence constants to be assessed numerically. We use this loss to train RBNO through finite element assembly and differentiation, without nonlinear solves during training or solution data beyond the snapshots used to construct the reduced basis. The same procedure applies to other predictors whose outputs lie in the conforming space. In our benchmarks, one evaluation of the loss and its gradient costs between $1/500$ and $1/30$ of a finite element solve.

\item \emph{Numerical validation across the three problem classes.}
We study semilinear diffusion--reaction, quasilinear diffusion, and compressible hyperelasticity benchmarks using the same data and discretizations as the corresponding linear benchmarks in~\cite{QCD}. The solution error and the square root of the loss excess remain comparable, with ratios similar to those observed in the linear case. Discrete stability tests also support the stability assumptions at the computed states. When evaluated in the same finite element space, RBNO achieves smaller PDE-compliant norm errors than grid-based baselines, including PCA-Net, FNO, and FC-PINO. These comparisons further show that small grid-based residuals can coexist with substantial PDE-compliant norm errors, even when the predictions are accurate in $L^2$.
\end{enumerate}

\subsection{Organization}
The remainder of the paper is organized as follows. We introduce the residual formulation with boundary lifting in \cref{sec:setup}. In \cref{sec:local_vc}, we establish local residual--error equivalence and account for errors in residual evaluation. We then verify the assumptions for three classes of nonlinear PDEs in \cref{sec:examples}, with the main results summarized in \cref{tab:example_summary}. In \cref{sec:rbno}, we develop the conforming finite element formulation, derive loss-excess estimates, and describe the discrete stability assessment and RBNO training procedure. We present numerical experiments in \cref{sec:numerics} to validate the error estimates, examine stability at the computed states, and compare RBNO with grid-based neural operators. Finally, we conclude in \cref{sec:conclusions}. Readers primarily interested in the method and numerical results may focus on \cref{sec:setup}, the statement of \cref{thm:local_vc}, \cref{tab:example_summary}, and \cref{sec:rbno,sec:numerics}.

\section{Parametric nonlinear FOSLS formulation}\label{sec:setup}
We formulate each model as a first-order system for the flux or stress $\sigma$ and the primal variable $u$. We remove the boundary data through lifts and measure the resulting residual in an $L^2$ product space.

Let $\Omega\subset\bR^d$, $d\in\{2,3\}$, be a bounded Lipschitz domain. We partition its boundary into a Dirichlet part $\Gamma_D$ of positive measure and a Neumann part $\Gamma_N=\partial\Omega\setminus\Gamma_D$. We use $(\cdot,\cdot)$ and $\norm{\cdot}$ for the $L^2(\Omega)$ inner product and norm, including their vector- and matrix-valued counterparts. All other norms carry a subscript. For scalar diffusion problems, we define
\[
 \U=H^1_{0,\Gamma_D}(\Omega),\qquad
 \Sigma=\HdivN
 :=\{\tau\in\Hdiv:\tau\cdot n=0\text{ on }\Gamma_N\},
\]
where the boundary conditions are understood in the trace sense. We equip the state space $\HH=\Sigma\times\U$ with the PDE-compliant norm
\begin{equation}\label{eq:state_norm}
 \norm{[\tau,v]}_\HH^2
 =\norm{\tau}_{H(\rdiv)}^2+\norm{v}_{H^1}^2
 =\norm{\tau}^2+\norm{\rdiv\tau}^2
  +\norm{v}^2+\norm{\nabla v}^2.
\end{equation}
For the elasticity problem, we take $\U=H^1_{0,\Gamma_D}(\Omega;\bR^d)$ and let $\Sigma$ consist of matrix fields with rows in $\Hdiv$ and vanishing traction on $\Gamma_N$. The norm is defined as in~\eqref{eq:state_norm}, with divergence taken row-wise. We impose no symmetry on the stress, since the first Piola--Kirchhoff stress is generally nonsymmetric.

We measure residuals in the self-dual product space
\[
 \LL^2=
 \begin{cases}
 L^2(\Omega)^d\times L^2(\Omega), & \text{for scalar diffusion},\\
 L^2(\Omega)^{d\times d}\times L^2(\Omega)^d, & \text{for elasticity}.
 \end{cases}
\]
Throughout, $C_P$ denotes the Poincar\'e constant satisfying
$\norm{v}\le C_P\norm{\nabla v}$ for $v\in\U$.
We summarize the notation in \cref{tab:notation}.

\subsection{Boundary data, sources, and lifts}
We remove inhomogeneous boundary data using the variational lifts of~\cite{QCD}. The Dirichlet lift $w\in H^1(\Omega)$ satisfies $w=u_D$ on $\Gamma_D$, and the flux lift $z\in\Hdiv$ satisfies $z\cdot n=g$ on $\Gamma_N$. For scalar problems, we take $w$ harmonic with $\nabla w\cdot n=0$ on $\Gamma_N$, and $z=\nabla\xi$, where $\xi\in\U$ solves $(\nabla\xi,\nabla v)=\langle g,v\rangle_{\Gamma_N}$ for all $v\in\U$; in particular $\rdiv z=0$. We use the same lifts for every candidate state and for the solution, so that the difference of two physical states equals the difference of their homogeneous parts.

We allow volume sources in $\U'$, provided they are given in the flux-free form
\begin{equation}\label{eq:flux_free}
 f(v)=(f_2,v)-(f_1,\nabla v),\qquad v\in\U,\qquad f_2\in L^2(\Omega),\ f_1\in L^2(\Omega)^d,
\end{equation}
used in~\cite{QCD}. The first-order unknown is then the shifted flux, the constitutive flux plus $f_1$. For the diffusion operator $-\rdiv(\pp\nabla u)$ this gives $\sigma=\pp\nabla u+f_1$, and the equilibrium equation is $-\rdiv\sigma=f_2$, plus reaction terms where present. If $f_1$ vanishes in a neighborhood of $\partial\Omega$, as in all our examples, $\sigma\cdot n$ coincides with the physical normal flux on the boundary. With $u=\uc+w$ and $\sigma=\sigc+z$, the trainable state is
\begin{equation}\label{eq:lifts}
 s=[\sigc,\uc]\in\HH .
\end{equation}
The lifts remove all boundary terms from the loss, but a nonlinear constitutive or reaction term is still evaluated at the physical state $\uc+w$.

\begin{remark}[Interfaces]\label{rem:interfaces}
If $\pp$ or $f_1$ jumps across an interface $\Gamma$, the weak form requires the normal component of $\pp\nabla u+f_1$ to be continuous across $\Gamma$. In the first-order system this transmission condition is part of the requirement $\sigma\in\Hdiv$. The residual of the second-order equation, evaluated pointwise off $\Gamma$, does not contain it, since the derivatives of a piecewise constant $\pp$ and $f_1$ vanish there.
\end{remark}

\subsection{Abstract residual}
Let $\mu$ be a probability measure on the parameter set $\PP$. We write the lifted residual as
\begin{equation}\label{eq:abstract_residual}
 \cR(s;\pp)=\cB_\pp s+\cK_\pp(\uc+w)-S_\pp,
 \qquad \cL(s;\pp):=\norm{\cR(s;\pp)}_{\LL^2}^2 ,
\end{equation}
where $\cB_\pp:\HH\to\LL^2$ is bounded and linear, $S_\pp$ collects the lifted data, and $\cK_\pp$ is the nonlinear part, acting on the physical primal variable and taking values in $\LL^2$. For the diffusion--reaction equation $-\rdiv(\pp\nabla u)+\kappa u^3=f$ of \cref{sec:sdr}, for instance,
\[
 \cB_\pp[\sigc,\uc]=(\sigc-\pp\nabla\uc,\,-\rdiv\sigc),\qquad
 \cK_\pp(u)=(0,\,\kappa u^3),\qquad
 S_\pp=(f_1+\pp\nabla w-z,\,f_2).
\]
This form covers a reaction term, in which $\cK_\pp$ acts on $u$, a diffusion coefficient that depends on $u$ or on $\nabla u$, and a stress law that depends on $\nabla u$; \cref{sec:examples} treats models of each kind. We do not consider constitutive laws that depend nonlinearly on the flux or stress. The map $\cK_\pp$ is defined on an open subset of a Banach space $\Y$ that embeds continuously into the primal space. We use $\Y=\U$ for the semilinear model and $\Y=W^{1,\infty}(\Omega)\cap\U$ for the gradient-dependent models, with $\norm{v}_{W^{1,\infty}}=\norm{v}_{L^\infty}+\norm{\nabla v}_{L^\infty}$ and the Frobenius norm for matrix-valued gradients. When $\norm{\cdot}_\Y$ is stronger than $\norm{\cdot}_\U$, the residual is defined only for states with primal component in $\Y$.

A selected solution branch $s^*(\pp)=[\sigma^{\circ,*}(\pp),u^{\circ,*}(\pp)]$ satisfies
\begin{equation}\label{eq:fiber_nonlinear}
 \cR(s^*(\pp);\pp)=0,\qquad \pp\in\PP .
\end{equation}
We assume that $s^*(\pp)$ is strongly measurable and belongs to $\X=L^2_\mu(\PP;\HH)$. Existence, selection, and regularity of the branch are properties of the individual model, which we discuss in \cref{sec:examples}.

\paragraph{Parameter sets}
The estimates below hold uniformly in $\pp$ provided that the coefficients admit uniform upper and positive lower bounds on $\PP$. For the two benchmarks with lognormal fields, we follow~\cite{QCD} and restrict $\PP$ to fields satisfying $m\le\hat m$ almost everywhere in $\Omega$. This restriction provides the coefficient upper bounds listed in \cref{tab:benchmarks}.

\begin{table}[tb]
\centering\small
\begin{tabularx}{\linewidth}{@{}lX@{}}\toprule
Notation & Meaning\\\midrule
$\HH=\Sigma\times\U$, $\LL^2$ & State space with the PDE-compliant norm~\eqref{eq:state_norm}; residual space\\
$s=[\sigc,\uc]$, $s^*(\pp)$ & Lifted state, with flux or stress $\sigc$ and primal variable $\uc$; selected solution branch\\
$w$, $z$; $f_1$, $f_2$ & Dirichlet and flux lifts; flux-free splitting~\eqref{eq:flux_free} of the source\\
$\cR$, $\cL$; $\cB_\pp$, $\cK_\pp$, $S_\pp$ & Residual and loss~\eqref{eq:abstract_residual}; linear part, nonlinear part, and data of the residual\\
$\cA_{\bar u,\pp}$ & Derivative of the residual at a state with primal component $\bar u$\\
$\Y$, $\cN_\pp$, $R_0$ & Regularity space; closed $\Y$-ball of radius $R_0$ about the branch\\
$\clin,\Clin$; $\MK$ & Squared bounds of the linearized residual; Lipschitz constant of $\cK_\pp'$\\
$R$, $\gamma$; $\Crel,\Ceff$ & Equivalence radius and margin; squared equivalence constants\\
$J_{h,q}$ & Loss as it is computed for a finite element state, with a quadrature rule of degree $q$\\
$C_P$, $C_4$, $C_S$ & Constants of $\U$ in $L^2$, $L^4$, and $L^6$ with respect to $\norm{\nabla\cdot}$\\
$\pmin,\pmax$; $\kappa$ & Bounds of the coefficient field $\pp$; strength of the cubic reaction\\
$c_B$, $C_B$ & Squared constants of the bounded-set estimate for the semilinear model\\
$M_a$; $M_F$, $j_0$, $K_J$, $L_J$ & Regularity bound for quasilinear diffusion; deformation bounds for elasticity\\
$\HH_h$, $\HH_r$; $\psi_j$, $\phi_k$ & Conforming finite element space and reduced space; their bases\\
$\eta_h$, $E_h$ & Residual of the finite element solution; error relative to it\\
$\bss$, $\bst$; $\bsg$; $X_h$ & Coefficient vectors of $s,t\in\HH_h$; assembled derivative~\eqref{eq:assembled_gradient} of the loss; Gram matrix of $\HH$\\
$s_r(\pp;\theta)$, $\bss_r(\pp;\theta)$ & RBNO and its coefficient vector, the output of the network, as in~\cite{QCD}\\
$\Phi_r$, $\Pi_r$, $P_r$ & POD basis functions, their coefficient vectors, and the $\HH$-orthogonal projector onto $\HH_r$\\
$\eta$, $\chi$ & Residual of a state; floor-corrected ratio of error and residual~\eqref{eq:numerical_metrics}\\
\bottomrule\end{tabularx}
\caption{Notation. Spaces are set in blackboard type; operators on the first-order state, the loss, and sets of states in calligraphic type; the parameter in Fraktur; coefficient vectors in bold.}\label{tab:notation}
\end{table}

\section{Local residual--error equivalence}\label{sec:local_vc}
In this section, we establish the main result of the paper: under stability and Lipschitz continuity assumptions on the linearization, the residual is equivalent to the solution error in a neighborhood of the solution branch. We state these assumptions in \cref{subsec:assumptions} and prove the equivalence in \cref{thm:local_vc}. We verify the assumptions for three classes of nonlinearity in \cref{sec:examples}. In \cref{sec:lifted}, we quantify how approximate residual evaluation and the use of a numerical reference solution affect the residual--error bounds.

\subsection{Stability and regularity of the linearization}\label{subsec:assumptions}
For each parameter, we fix the closed neighborhood of the selected branch
\begin{equation}\label{eq:tube_def}
 \cN_\pp=\{\bar u\in\Y:\norm{\bar u-u^{\circ,*}(\pp)}_\Y\le R_0\},
\end{equation}
with a radius $R_0>0$ independent of $\pp$, and assume that the domain of $\bar u\mapsto \cK_\pp(\bar u+w)$ contains $\cN_\pp$. The set $\cN_\pp$ is convex; no compactness is required.

\begin{assumption}[Linearized stability]\label{ass:local_stability}
The map $\bar u\mapsto \cK_\pp(\bar u+w)$ is continuously Fr\'echet differentiable on $\cN_\pp$ with respect to the topology of $\Y$, and its derivative at each $\bar u\in\cN_\pp$ extends to a bounded operator $\cK_\pp'(\bar u):\U\to\LL^2$. Defining
\begin{equation}
\label{Aubar}
 \cA_{\bar u,\pp}[\tau,v]:=\cB_\pp[\tau,v]+\cK_\pp'(\bar u)v ,
\end{equation}
there are constants $0<\clin\le\Clin<\infty$, independent of $\bar u$ and $\pp$, such that
\begin{equation}\label{eq:uniform_linearized_stability}
 \clin\norm{t}_\HH^2\le\norm{\cA_{\bar u,\pp}t}_{\LL^2}^2\le\Clin\norm{t}_\HH^2,
 \qquad t\in\HH,\quad \bar u\in\cN_\pp .
\end{equation}
\end{assumption}

\begin{assumption}[Lipschitz derivative]\label{ass:lipschitz_derivative}
There is a constant $\MK\ge0$, independent of $\pp$, such that
\begin{equation}\label{eq:uniform_lip_derivative}
 \norm{\bigl(\cK_\pp'(\bar u_1)-\cK_\pp'(\bar u_2)\bigr)v}_{\LL^2}
 \le \MK\norm{\bar u_1-\bar u_2}_\Y\norm{v}_\U ,
 \qquad \bar u_1,\bar u_2\in\cN_\pp,\ v\in\U .
\end{equation}
\end{assumption}

\Cref{ass:local_stability} states that the linearized problem is well posed in $\HH$, uniformly near the branch and in the parameter. \Cref{ass:lipschitz_derivative} bounds the change of the linearization with the linearization point. The constants $\clin$ and $\MK$ enter the lower bound of the theorem below through the margin $\sqrt{\clin}-\tfrac12\MK R$, where $R\le R_0$ bounds the distance of the primal component of a state from the branch in the norm of $\Y$.

\paragraph{Example}
For the diffusion--reaction model, with $\cK_\pp(u)=(0,\kappa u^3)$ as in \cref{sec:setup}, the linearized residual
\[
 \cA_{\bar u,\pp}[\tau,v]=\bigl(\tau-\pp\nabla v,\;-\rdiv\tau+3\kappa(\bar u+w)^2v\bigr)
\]
is the residual operator of a linear diffusion--reaction problem with the nonnegative reaction coefficient $3\kappa(\bar u+w)^2$. The energy argument for such problems gives~\eqref{eq:uniform_linearized_stability} (\cref{sec:sdr}), and H\"older's inequality with the embedding $H^1(\Omega)\subset L^6(\Omega)$ gives~\eqref{eq:uniform_lip_derivative} with $\Y=\U$. When the flux or stress depends nonlinearly on $\nabla u$, as for the $p$-Laplacian and for elasticity in \cref{sec:examples}, the difference $(\cK_\pp'(\bar u_1)-\cK_\pp'(\bar u_2))v$ contains $\nabla(\bar u_1-\bar u_2)$ multiplied by $\nabla v$. Since $\nabla v$ is only in $L^2$, a bound in terms of $\norm{v}_\U$ then requires $\nabla(\bar u_1-\bar u_2)\in L^\infty$, and we use $\Y=W^{1,\infty}(\Omega)\cap\U$ for these models and for the quasilinear diffusion model of \cref{sec:quasilinear}.

In~\eqref{eq:uniform_linearized_stability} the direction $t$ ranges over the whole state space, and only the linearization point is restricted to $\cN_\pp$. This distinction matters for quasilinear diffusion and elasticity, where $\cN_\pp$ is a $W^{1,\infty}$ neighborhood and an $H^1$ bound on the error does not place a state in it. The lower bound in the theorem below uses stability only at the branch, and stability at the branch implies uniform stability on $\cN_\pp$ for small $R_0$: if $\beta_*\norm{t}_\HH\le\norm{\cA_{u^{\circ,*},\pp}t}_{\LL^2}\le\Lambda_*\norm{t}_\HH$ for $t\in\HH$, then~\eqref{eq:uniform_lip_derivative} gives
\[
 (\beta_*-\MK R_0)\norm{t}_\HH\le\norm{\cA_{\bar u,\pp}t}_{\LL^2}\le(\Lambda_*+\MK R_0)\norm{t}_\HH ,
 \qquad \bar u\in\cN_\pp ,
\]
which is a uniform bound when $\MK R_0<\beta_*$.

\subsection{The residual--error estimate}\label{subsec:equivalence}
\begin{theorem}[Local residual--error equivalence]\label{thm:local_vc}
Let \cref{ass:local_stability,ass:lipschitz_derivative} hold, and let $0<R\le R_0$ satisfy
\begin{equation}\label{eq:R_smallness}
 \gamma:=\sqrt{\clin}-\tfrac12\MK R>0 .
\end{equation}
Let $s=[\tau,v]\in\HH$ be a state whose flux component $\tau\in\Sigma$ is arbitrary and whose primal component satisfies $v\in\Y$ and $\norm{v-u^{\circ,*}(\pp)}_\Y\le R$. With $e=s-s^*(\pp)$ and $e_u=v-u^{\circ,*}(\pp)$,
\begin{equation}\label{eq:local_vc_pointwise}
 \bigl(\sqrt{\clin}-\tfrac12\MK\norm{e_u}_\Y\bigr)\norm{e}_\HH
 \le\norm{\cR(s;\pp)}_{\LL^2}
 \le\sqrt{\Clin}\,\norm{e}_\HH .
\end{equation}
In particular, with the reliability and efficiency constants $\Crel=\gamma^2$ and $\Ceff=\Clin$,
\begin{equation}\label{eq:local_vc_unsquared}
 \Crel\norm{s-s^*(\pp)}_\HH^2\le\cL(s;\pp)\le\Ceff\norm{s-s^*(\pp)}_\HH^2 .
\end{equation}
\end{theorem}
\begin{proof}
Since $\cN_\pp$ is convex, the segment $u^{\circ,*}+te_u$, $0\le t\le1$, lies in $\cN_\pp$. By~\eqref{eq:fiber_nonlinear} and the fundamental theorem of calculus,
\begin{equation}\label{eq:residual_split}
 \cR(s;\pp)=\cR(s;\pp)-\cR(s^*(\pp);\pp)=\cA_{u^{\circ,*},\pp}e+\omega,\qquad
 \omega=\int_0^1\bigl(\cK_\pp'(u^{\circ,*}+te_u)-\cK_\pp'(u^{\circ,*})\bigr)e_u\,dt ,
\end{equation}
and~\eqref{eq:uniform_lip_derivative} yields
\begin{equation}\label{eq:taylor_remainder}
 \norm{\omega}_{\LL^2}\le\tfrac12\MK\norm{e_u}_\Y\norm{e_u}_\U\le\tfrac12\MK\norm{e_u}_\Y\norm{e}_\HH .
\end{equation}
With the lower bound in~\eqref{eq:uniform_linearized_stability} at the branch, the reverse triangle inequality gives the lower bound in~\eqref{eq:local_vc_pointwise},
\[
 \norm{\cR(s;\pp)}_{\LL^2}\ge\norm{\cA_{u^{\circ,*},\pp}e}_{\LL^2}-\norm{\omega}_{\LL^2}\ge\bigl(\sqrt{\clin}-\tfrac12\MK\norm{e_u}_\Y\bigr)\norm{e}_\HH .
\]
For the upper bound we write $\cR(s;\pp)=\int_0^1\cA_{u^{\circ,*}+te_u,\pp}\,e\,dt$ and use the upper bound in~\eqref{eq:uniform_linearized_stability} on the segment. Squaring and using $\norm{e_u}_\Y\le R$ gives~\eqref{eq:local_vc_unsquared}.
\end{proof}

The condition~\eqref{eq:R_smallness} limits the radius $R$ by the ratio $2\sqrt{\clin}/\MK$ of stability to nonlinearity. For a linear problem, $\MK=0$, and the estimate is global and reduces to~\eqref{eq:linear_vc}. If~\eqref{eq:uniform_linearized_stability} holds only at the branch, the same proof gives the lower bound and the upper factor $\sqrt{\Clin}+\tfrac12\MK R$. The theorem yields local uniqueness of the zero in the admissible set
\begin{equation}\label{eq:admissible_set}
 \cS_\pp(R)=\{[\tau,v]\in\Sigma\times\Y:\ \norm{v-u^{\circ,*}(\pp)}_\Y\le R\}
\end{equation}
and the a posteriori bound
\begin{equation}\label{eq:certificate}
 \norm{s-s^*(\pp)}_\HH\le\gamma^{-1}\cL(s;\pp)^{1/2},\qquad s\in\cS_\pp(R).
\end{equation}
The bound becomes a certificate once the margin $\gamma$ and the membership $s\in\cS_\pp(R)$ are known independently of the residual. For the semilinear model of \cref{sec:sdr} both are available, and \cref{cor:semilinear_certificate} gives a fully computable version of~\eqref{eq:certificate}. A margin estimated from sampled eigenvalues or from fitted residual--error ratios may serve as a practical orientation but does not have a rigorous status.

\subsection{Evaluation errors}\label{sec:lifted}
We evaluate the residual approximately, through numerical integration and approximations of the problem data and boundary lifts. These approximations affect the computed residual, while modifying the constitutive law changes the underlying PDE. We denote by $J_{h,q}(s;\pp)$ the computed loss of a finite element state $s$, evaluated using a quadrature rule of degree $q$.

\begin{proposition}[Evaluation and reference errors]\label{prop:consistency}
Let $s\in\cS_\pp(R)$ and $E=\norm{s-s^*(\pp)}_\HH$. If $\eta^-\le\cL(s;\pp)^{1/2}\le\eta^+$ for computable numbers $\eta^\pm$, then
\begin{equation}\label{eq:residual_defect}
 \eta^-/\sqrt{\Ceff}\le E\le\eta^+/\sqrt{\Crel}.
\end{equation}
Admissible choices are $\eta^\pm=\norm{\widetilde\cR}_{\LL^2}\pm\varepsilon_R$ when an approximate residual satisfies $\norm{\widetilde\cR-\cR(s;\pp)}_{\LL^2}\le\varepsilon_R$, and $\eta^\pm=\sqrt{(J_{h,q}\pm\varepsilon_q)_+}$, with $(x)_+=\max\{x,0\}$, when the computed loss $J_{h,q}=J_{h,q}(s;\pp)$ satisfies $|J_{h,q}-\cL(s;\pp)|\le\varepsilon_q$. If a reference $s_h$ satisfies $\norm{s_h-s^*(\pp)}_\HH\le\varepsilon_h$, then $E_h=\norm{s-s_h}_\HH$ satisfies $|E-E_h|\le\varepsilon_h$.
\end{proposition}
\begin{proof}
By~\eqref{eq:local_vc_unsquared}, $\sqrt{\Crel}\,E\le\cL(s;\pp)^{1/2}\le\sqrt{\Ceff}\,E$, which gives~\eqref{eq:residual_defect}. For the first choice of $\eta^\pm$, the bounds on $\cL(s;\pp)^{1/2}=\norm{\cR(s;\pp)}_{\LL^2}$ follow from the triangle inequality, and for the second from $|J_{h,q}-\cL(s;\pp)|\le\varepsilon_q$. The last statement is the reverse triangle inequality $|E-E_h|\le\norm{s_h-s^*(\pp)}_\HH$.
\end{proof}

\section{Three classes of nonlinearity}\label{sec:examples}
For a given model, \cref{ass:local_stability,ass:lipschitz_derivative} raise two questions: in which norm is the derivative of the nonlinear part Lipschitz, and why is the linearization stable? We answer them for three classes that differ in where the nonlinearity enters: the lower-order term in semilinear diffusion--reaction, the principal part through a coefficient that depends on the solution or its gradient in quasilinear diffusion, and the stress law in compressible hyperelasticity. We collect the answers in \cref{tab:example_summary}. In its first two columns, where stability is proved, the estimate holds on every set of states that is bounded in the norm of the neighborhood, without a smallness condition on the error; in the third, where stability is assumed, it holds in a neighborhood of the branch. Only for the reaction is the error bound computable without information on the solution; for the other two classes the constants involve a bound on the solution in $W^{1,\infty}$ or on the inverse of the tangent at the branch.

\begin{table}[!htb]
\centering\small
\begin{tabularx}{\linewidth}{@{}>{\raggedright\arraybackslash}p{.21\linewidth}>{\raggedright\arraybackslash}X>{\raggedright\arraybackslash}X>{\raggedright\arraybackslash}X@{}}\toprule
 & Semilinear diffusion--reaction (\cref{sec:sdr}) & Quasilinear diffusion (\cref{sec:quasilinear}) & Compressible hyperelasticity (\cref{sec:hyperelastic})\\\midrule
Nonlinear term & reaction $\kappa u^3$ & flux $\pp(1+u^2)\nabla u$ & stress $P(I+\nabla u;\pp)$\\\addlinespace
It enters & the lower-order term & the principal part & the principal part; defined only for $\det F>0$\\\addlinespace
Norm of the neighborhood $\cN_\pp$ (\cref{ass:lipschitz_derivative}) & $H^1$, the norm of the error & $W^{1,\infty}$ & $W^{1,\infty}$, with $\det F\ge j_0$\\\addlinespace
Stability of the linearization (\cref{ass:local_stability}) & proved: monotonicity (\cref{prop:semilinear}) & proved: factorization of the tangent (\cref{prop:quasilinear}) & assumed: coercive or invertible tangent at the branch (\cref{prop:hyperelastic})\\\addlinespace
Estimate holds on & every bounded set in $H^1$ & every bounded set in $W^{1,\infty}$ & a neighborhood of the branch\\\addlinespace
Computable error bound & yes (\cref{cor:semilinear_certificate}) & given a bound $M_*$ on the solution in $W^{1,\infty}$ & no; $a_*$ and $C_*$ are not known\\\addlinespace
Benchmark (\cref{tab:benchmarks}) & heat conduction with a cubic reaction & flow with a solution-dependent permeability & Neo-Hookean body under traction\\
\bottomrule
\end{tabularx}
\caption{The three classes side by side: how the nonlinear term enters, the two properties required by \cref{thm:local_vc}, and what they imply for the residual--error estimate. Each column corresponds to one subsection and to one benchmark of \cref{sec:numerics}.}\label{tab:example_summary}
\end{table}

\subsection{Semilinear diffusion--reaction}\label{sec:sdr}
Let $\pp\in L^\infty(\Omega)$ with $0<\pmin\le\pp\le\pmax$, and $\kappa\in L^\infty(\Omega)$ with $0\le\kappa\le\kappa_{\max}$. We consider
\begin{equation}\label{eq:diffusion_reaction_note}
 -\rdiv(\pp\nabla u)+\kappa u^3=f\ \ \text{in }\Omega,\qquad
 u=u_D\ \ \text{on }\Gamma_D,\qquad
 \pp\nabla u\cdot n=g\ \ \text{on }\Gamma_N,
\end{equation}
with $f$ given in the flux-free form~\eqref{eq:flux_free} and $f_1$ vanishing near $\partial\Omega$, as in \cref{sec:setup}. In the numerical study, we use a piecewise constant conductivity and a constant reaction coefficient. With the lifts of \cref{sec:setup}, the solution is the function $u\in w+\U$ with
\[
 (\pp\nabla u+f_1-z,\nabla v)+(\kappa u^3,v)=(f_2,v),\qquad v\in\U .
\]
It is the unique minimizer over $w+\U$ of the strictly convex and coercive functional $\tfrac12(\pp\nabla u,\nabla u)+\tfrac14(\kappa,u^4)+(f_1-z,\nabla u)-(f_2,u)$. Since $d\le3$, $H^1(\Omega)\subset L^6(\Omega)$ and $u^3\in L^2(\Omega)$, so the shifted flux $\sigma=\pp\nabla u+f_1$ belongs to $\Hdiv$. With $\sigc=\sigma-z$ and $\uc=u-w$ the lifted residual is
\begin{equation}\label{eq:diffusion_reaction_residual}
 \cR_{\mathrm{dr}}([\sigc,\uc];\pp)=
 \begin{pmatrix}\sigc+z-\pp\nabla(\uc+w)-f_1\\ -\rdiv\sigc+\kappa(\uc+w)^3-f_2\end{pmatrix},
 \qquad \cL_{\mathrm{dr}}=\norm{\cR_{\mathrm{dr}}}_{\LL^2}^2 .
\end{equation}
Let $u^*=u^{\circ,*}+w$ denote the physical solution. Testing the weak equation with $v=u^*-w=u^{\circ,*}\in\U$ gives the a priori bound
\begin{equation}\label{eq:dr_branch_bound}
 \norm{\nabla(u^*-w)}\le\frac{a+\sqrt{a^2+4\pmin b}}{2\pmin},\qquad
 a=\norm{\pp\nabla w}+\norm{f_1-z}+C_P\norm{f_2},\quad b=\tfrac{27}{256}(\kappa,w^4).
\end{equation}
Indeed, with $v=u^*-w$ the weak equation gives
\[
 \pmin\norm{\nabla v}^2\le a\norm{\nabla v}-\bigl(\kappa,(u^*)^4-(u^*)^3w\bigr)\le a\norm{\nabla v}+b,
\]
because pointwise $(u^*)^4-(u^*)^3w\ge t^4-\alpha t^3\ge-\tfrac{27}{256}\alpha^4$ with $t=|u^*|$ and $\alpha=|w|$; the function $t\mapsto t^4-\alpha t^3$ attains its minimum over $t\ge0$ at $t=\tfrac34\alpha$. The bound~\eqref{eq:dr_branch_bound} is the larger root of $\pmin x^2-ax-b$. The branch is therefore bounded in $H^1$ uniformly over any parameter set on which the data are bounded.

Below, $C_4$ and $C_S$ denote the constants in $\norm{v}_{L^4}\le C_4\norm{\nabla v}$ and $\norm{v}_{L^6}\le C_S\norm{\nabla v}$ for $v\in\U$. We compute them for the geometry of the numerical study in \cref{lem:embedding}. H\"older's inequality gives $C_4\le|\Omega|^{1/12}C_S$; \cref{lem:embedding} bounds $C_4$ directly and more sharply, $C_4<0.70$ against $C_S<1.03$ on the unit square.

The following proposition rests on the secant form of the cubic, $u^3-u^{*3}=(u^2+uu^*+u^{*2})(u-u^*)$. By this identity, the residual of any state is a linear diffusion--reaction operator, with the nonnegative reaction coefficient $\csec$, applied to the error of the state. The energy argument for linear problems therefore applies without linearization and without a restriction on the size of the error.

\begin{proposition}[Equivalence on bounded sets]\label{prop:semilinear}
Let $s=[\tau,v]\in\Sigma\times\U$ be any state with physical primal component $u=v+w$, and write $e=s-s^*(\pp)=[e_\sigma,e_u]$, $\eta=\norm{\cR_{\mathrm{dr}}(s;\pp)}_{\LL^2}$, $\csec=\kappa(u^2+uu^*+u^{*2})\ge0$, and $K=\sqrt{(1+C_P^2)/\pmin}$.
\begin{enumerate}
\item[(a)] $\norm{\pp^{1/2}\nabla e_u}\le K\eta$ and $\norm{\csec^{1/2}e_u}\le K\eta$. Consequently
\[
 \norm{e_u}_{H^1}\le\frac{1+C_P^2}{\pmin}\,\eta,\qquad
 \norm{e_\sigma}\le\bigl(1+\sqrt{\pmax}\,K\bigr)\eta .
\]
\item[(b)] $\norm{\rdiv e_\sigma}\le\eta+\norm{\csec e_u}$. If $\norm{\csec e_u}\le\zeta\eta$, then $c_B\norm{e}_\HH^2\le\cL_{\mathrm{dr}}(s;\pp)$ with
\begin{equation}\label{eq:dr_constants}
 c_B^{-1}=\frac{(1+C_P^2)^2}{\pmin^2}+\Bigl(1+\sqrt{\pmax}\,K\Bigr)^2+(1+\zeta)^2 .
\end{equation}
Admissible choices are $\zeta=K\norm{\csec}_{L^\infty}^{1/2}$ and $\zeta=C_SK\norm{\csec}_{L^3}/\sqrt{\pmin}$.
\item[(c)] $\cL_{\mathrm{dr}}(s;\pp)\le C_B\norm{e}_\HH^2$ with $C_B=2\max\{1,\pmax^2+C_S^2\norm{\csec}_{L^3}^2\}$.
\end{enumerate}
No smallness condition on $e$ is required, and the bounds in (a) do not depend on $\kappa$.
\end{proposition}
\begin{proof}
Since $\cR_{\mathrm{dr}}(s^*)=0$ and $\kappa(u^3-u^{*3})=\csec e_u$, the residual of $s$ is $\cR_{\mathrm{dr}}(s)=(r_1,r_2)$ with $r_1=e_\sigma-\pp\nabla e_u$ and $r_2=-\rdiv e_\sigma+\csec e_u$. Integration by parts, using $e_\sigma\cdot n=0$ on $\Gamma_N$ and $e_u=0$ on $\Gamma_D$, gives the energy identity
\[
 (\pp\nabla e_u,\nabla e_u)+(\csec e_u,e_u)=(r_2,e_u)-(r_1,\nabla e_u).
\]
By Poincar\'e's inequality the right-hand side is at most $(C_P\norm{r_2}+\norm{r_1})\norm{\nabla e_u}\le\sqrt{1+C_P^2}\,\eta\,\norm{\nabla e_u}\le K\eta\norm{\pp^{1/2}\nabla e_u}$. Both terms on the left are nonnegative, which gives the two bounds in (a). The bound on $e_u$ follows from $\norm{e_u}_{H^1}^2\le(1+C_P^2)\norm{\nabla e_u}^2\le(1+C_P^2)\pmin^{-1}\norm{\pp^{1/2}\nabla e_u}^2$, and the bound on $e_\sigma$ from $e_\sigma=r_1+\pp\nabla e_u$ and $\norm{\pp\nabla e_u}\le\sqrt{\pmax}\norm{\pp^{1/2}\nabla e_u}$. Part (b) follows from $\rdiv e_\sigma=\csec e_u-r_2$; the two choices of $\zeta$ come from $\norm{\csec e_u}\le\norm{\csec^{1/2}}_{L^\infty}\norm{\csec^{1/2}e_u}$ and from $\norm{\csec e_u}\le\norm{\csec}_{L^3}\norm{e_u}_{L^6}\le C_S\norm{\csec}_{L^3}\norm{\nabla e_u}$. Part (c) follows from $\norm{r_1}^2\le2\norm{e_\sigma}^2+2\pmax^2\norm{\nabla e_u}^2$ and $\norm{r_2}^2\le2\norm{\rdiv e_\sigma}^2+2C_S^2\norm{\csec}_{L^3}^2\norm{\nabla e_u}^2$.
\end{proof}

The flux bound in (a) depends on $\sqrt{\pmax/\pmin}$ rather than on $\pmax/\pmin$, because the energy identity controls $\pp^{1/2}\nabla e_u$, and the reaction enters only through the divergence of the flux error. For the conductivities of the numerical study, $\pmin=0.1$ and $\pmax=10$, and with $C_P=1/\pi$, this gives $c_B^{-1/2}\approx16$ for $\kappa=0$, against more than $100$ if $\norm{e_\sigma}$ is bounded through $\pmax\norm{\nabla e_u}$. Since $c_B$ decreases like $\pmin^2$, the benchmark of \cref{sec:numerics} keeps the conductivity bounded below.

For a finite element candidate, $u$ is bounded, and the factor $\zeta$ can be computed from the candidate and its residual alone.

\begin{corollary}[A posteriori error bound]\label{cor:semilinear_certificate}
In the setting of \cref{prop:semilinear}, let $u\in L^\infty(\Omega)$, let $\eta=\norm{\cR_{\mathrm{dr}}(s;\pp)}_{\LL^2}$ be the residual of $s$, and put $G=\sqrt{1+C_P^2}\,\eta/\pmin$. Then $\norm{\csec e_u}\le\zeta(\eta)\eta$ with
\begin{equation}\label{eq:zeta_computable}
 \zeta(\eta)=\kappa_{\max}\,\frac{\sqrt{1+C_P^2}}{\pmin}\Bigl(3C_P\norm{u}_{L^\infty}^2+3C_4^2\norm{u}_{L^\infty}G+C_S^3G^2\Bigr),
\end{equation}
and consequently
\begin{equation}\label{eq:semilinear_certificate}
 \norm{s-s^*(\pp)}_\HH\le\Delta(s):=c_B(\zeta(\eta))^{-1/2}\,\eta .
\end{equation}
Since $\Delta$ is increasing in $\eta$, the bound remains valid when $\eta$ is replaced by an upper bound $\eta^+$ as in \cref{prop:consistency}.
\end{corollary}
\begin{proof}
Writing $u^*=u-e_u$ gives $\csec e_u=\kappa(3u^2e_u-3ue_u^2+e_u^3)$, hence
$\norm{\csec e_u}\le\kappa_{\max}(3\norm{u}_{L^\infty}^2\norm{e_u}+3\norm{u}_{L^\infty}\norm{e_u}_{L^4}^2+\norm{e_u}_{L^6}^3)$. Each norm of $e_u$ is bounded by the corresponding constant times $\norm{\nabla e_u}\le G$, which follows from part (a) of \cref{prop:semilinear}.
\end{proof}

All quantities in~\eqref{eq:zeta_computable} are available once a prediction has been evaluated: the bounds of the sampled conductivity, $\kappa_{\max}$, the constants of \cref{lem:embedding}, the residual, and $\norm{u}_{L^\infty}$, which for piecewise quadratic fields is bounded by the largest Bernstein coefficient on each element. No a priori knowledge of the branch is needed. The bound combines worst-case constants, and its leading factor grows like $\kappa_{\max}\norm{u}_{L^\infty}^2/\pmin$; in the experiments, we therefore use the residual itself as the error estimator.

The local hypotheses also hold with $\Y=\U$: for physical states with $L^6$ norms at most $M_6$, H\"older's inequality gives $\MK\le6\kappa_{\max}C_S^2M_6$ with respect to $\norm{\nabla\cdot}$, and \cref{ass:local_stability} follows from \cref{prop:semilinear} applied to the linearized operator, whose reaction coefficient $3\kappa(\bar u+w)^2$ is nonnegative. \Cref{prop:semilinear} makes the local theory unnecessary for this model: the residual is equivalent to the error on bounded sets, and there is no branch to select.

\subsection{Quasilinear diffusion and its Kirchhoff structure}\label{sec:quasilinear}
We consider
\begin{equation}\label{eq:nonlinear_poisson}
 -\rdiv\bigl(\pp(1+u^2)\nabla u\bigr)=f\ \ \text{in }\Omega,
 \qquad u=u_D\ \ \text{on }\Gamma_D,\qquad
 \pp(1+u^2)\nabla u\cdot n=g\ \ \text{on }\Gamma_N,
\end{equation}
with $\pp\in L^\infty(\Omega)$, $0<\pmin\le\pp\le\pmax$, and $f\in L^2(\Omega)$. The flux is $\sigma=\pp(1+u^2)\nabla u$. With lifts $w\in W^{1,\infty}(\Omega)$ and $z\in\Hdiv$ as in \cref{sec:setup}, and $f_z=f+\rdiv z$, the lifted residual is
\begin{equation}\label{eq:diffusion_residual}
 \cR_{\mathrm{diff}}([\sigc,\uc];\pp)=
 \begin{pmatrix}
 \sigc+z-\pp\bigl(1+(\uc+w)^2\bigr)\nabla(\uc+w)\\
 \rdiv\sigc+f_z
 \end{pmatrix}.
\end{equation}
The constitutive term is defined for primal components in $\Y=W^{1,\infty}(\Omega)\cap\U$. We assume that the branch has this regularity, $\sup_\pp\norm{u^{\circ,*}(\pp)}_{W^{1,\infty}}\le M_*$, and that $\norm{w}_{W^{1,\infty}}\le M_w$. Such bounds follow from regularity theory for the Kirchhoff variable introduced below when the coefficient, the data, and the domain are smooth enough~\cite[Ch.~8--9]{gilbargTrudinger2001}; on polygonal domains with mixed boundary conditions, the corners where Dirichlet and Neumann parts meet are the main obstruction~\cite[Ch.~3]{grisvard1985}. On the neighborhood $\cN_\pp$ of radius $R_0$ the states $a=\bar u+w$ satisfy
\begin{equation}\label{eq:diffusion_a_bound}
 \norm{a}_{W^{1,\infty}}\le M_a:=M_*+M_w+R_0 .
\end{equation}
The tangent constitutive map has the factorized form
\begin{equation}\label{eq:diffusion_linearized}
 T_av=\pp\bigl((1+a^2)\nabla v+2av\nabla a\bigr)=\pp\nabla\bigl((1+a^2)v\bigr),
\end{equation}
in which the apparently lower-order term $2\pp av\nabla a$ combines with the principal part into the gradient of a product. This is the key to the stability estimate. The linearized residual is $\cA_{\bar u,\pp}[\tau,v]=(\tau-T_av,\rdiv\tau)$ with $a=\bar u+w$, and the following lemma reduces its stability to the invertibility of the second-order operator $v\mapsto-\rdiv(T_av)$ in weak form. We use the lemma again for the $p$-Laplacian and for hyperelasticity.

\begin{lemma}[FOSLS reduction]\label{lem:fosls_reduction}
Let $T:\U\to L^2(\Omega)^d$ be bounded with $\norm{T}\le L_T$, define $B:\U\to\U'$ by $\langle Bv,\varphi\rangle=(Tv,\nabla\varphi)$, and assume that $B$ is an isomorphism with $\norm{B^{-1}}_{\U'\to\U}\le C_I$. Then $\cA[\tau,v]=(\tau-Tv,\rdiv\tau)$ satisfies $c\norm{[\tau,v]}_\HH^2\le\norm{\cA[\tau,v]}_{\LL^2}^2\le C\norm{[\tau,v]}_\HH^2$ with
\begin{equation}\label{eq:reduction_constants}
 c^{-1}=C_I^2+(1+L_TC_I)^2+1,\qquad C=2\max\{1,L_T^2\}.
\end{equation}
The same holds row-wise for vector-valued problems.
\end{lemma}
\begin{proof}
Put $r_1=\tau-Tv$, $r_2=\rdiv\tau$, and $\eta^2=\norm{r_1}^2+\norm{r_2}^2$. Since $\tau\cdot n=0$ on $\Gamma_N$ and $\varphi=0$ on $\Gamma_D$, we have $\langle Bv,\varphi\rangle=-(r_2,\varphi)-(r_1,\nabla\varphi)$, so that $\norm{v}_\U\le C_I\eta$. Then $\norm{\tau}\le(1+L_TC_I)\eta$ and $\norm{\rdiv\tau}\le\eta$, which gives the lower constant. The upper constant follows from $\norm{\tau-Tv}^2\le2\norm{\tau}^2+2L_T^2\norm{v}_\U^2$.
\end{proof}

\begin{proposition}[Tangent stability]\label{prop:quasilinear}
Under~\eqref{eq:diffusion_a_bound}, \cref{ass:local_stability,ass:lipschitz_derivative} hold with
\begin{equation}\label{eq:diffusion_constants}
 \MK=6\pmax M_a,\qquad
 L_T=\pmax(1+3M_a^2),\qquad
 C_I=\Bigl(1+\tfrac{3\sqrt3}{8}M_a\Bigr)\frac{1+C_P^2}{\pmin},
\end{equation}
where $\clin$ and $\Clin$ are given by~\eqref{eq:reduction_constants}. No injectivity or smallness assumption on the tangent is needed.
\end{proposition}
\begin{proof}
For $a_1,a_2$ in the bounded set,
\[
 (T_{a_1}-T_{a_2})v=\pp(a_1+a_2)(a_1-a_2)\nabla v+2\pp(a_1-a_2)v\nabla a_1+2\pp a_2v\nabla(a_1-a_2),
\]
and each of the three terms is bounded in $L^2$ by $2\pmax M_a\norm{a_1-a_2}_{W^{1,\infty}}\norm{v}_{H^1}$, which gives $\MK$. The bound on $L_T$ follows from~\eqref{eq:diffusion_linearized} in the same way. Let $b=1+a^2$. Multiplication by $b$ is an isomorphism of $\U$ because $b\ge1$ and $b,b^{-1}\in W^{1,\infty}$, and $\norm{b^{-1}v}_{H^1}\le(1+\norm{\nabla b^{-1}}_{L^\infty})\norm{v}_{H^1}$ with $\norm{\nabla b^{-1}}_{L^\infty}=\norm{2a\nabla a/(1+a^2)^2}_{L^\infty}\le\tfrac{3\sqrt3}{8}M_a$, since $2|t|/(1+t^2)^2\le3\sqrt3/8$. Let $B_a$ be the weak operator of \cref{lem:fosls_reduction} for $T=T_a$, and let $B_0$ be the weak operator $v\mapsto(\pp\nabla v,\nabla\cdot)$ of the linear problem, which is $B_a$ for $a=0$. Then~\eqref{eq:diffusion_linearized} gives
\begin{equation}\label{eq:kirchhoff_factorization}
 B_av=B_0(bv),\qquad B_a^{-1}\ell=b^{-1}B_0^{-1}\ell,\qquad v\in\U,\ \ell\in\U',
\end{equation}
and the Lax--Milgram lemma gives $\norm{B_0^{-1}}_{\U'\to\U}\le(1+C_P^2)/\pmin$. \Cref{lem:fosls_reduction} completes the proof.
\end{proof}

The same structure is visible at the level of the PDE. The Kirchhoff transformation
\begin{equation}\label{eq:kirchhoff}
 \tilde u=Q(u):=u+\tfrac13u^3,\qquad \nabla\tilde u=(1+u^2)\nabla u,\qquad -\rdiv(\pp\nabla\tilde u)=f,
\end{equation}
turns~\eqref{eq:nonlinear_poisson} into a linear problem with Dirichlet data $Q(u_D)$ and unchanged Neumann data, and $u$ is recovered pointwise from the strictly increasing map $Q$. The flux is the same in both formulations, $\pp(1+u^2)\nabla u=\pp\nabla\tilde u$. For the benchmark of \cref{sec:numerics}, where $u_D=1-x_1$ on the vertical edges and $g=0$, the transformed problem is the Darcy problem of~\cite{QCD} with Dirichlet data $Q(u_D)=\tfrac43(1-x_1)$. For this particular law a linear FOSLS model for $\tilde u$ is therefore an alternative to the nonlinear model; we use the benchmark to test the nonlinear formulation on a problem whose stability is fully understood.

The transformation also gives a two-sided estimate in the original variables on bounded sets. If two physical states $u$ and $u^*$ are bounded by $M$ in $W^{1,\infty}$, then $Q(u)-Q(u^*)=b_{\mathrm{sec}}(u-u^*)$ with $b_{\mathrm{sec}}=1+\tfrac13(u^2+uu^*+u^{*2})\ge1$, $\norm{b_{\mathrm{sec}}}_{L^\infty}\le1+M^2$, and $\norm{\nabla b_{\mathrm{sec}}}_{L^\infty}\le2M^2$. Since the flux is $\pp\nabla Q(u)$ and $\cR_{\mathrm{diff}}(s^*(\pp);\pp)=0$, the residual of a state $s$ with error $e=[e_\sigma,e_u]=s-s^*(\pp)$ is
\[
 \cR_{\mathrm{diff}}(s;\pp)=\bigl(e_\sigma-\pp\nabla(b_{\mathrm{sec}}e_u),\;\rdiv e_\sigma\bigr).
\]
This is the operator of \cref{lem:fosls_reduction} with $Tv=\pp\nabla(b_{\mathrm{sec}}v)$, and the proof of \cref{prop:quasilinear} with $b_{\mathrm{sec}}$ in place of $b$ gives the equivalence without a smallness condition on $u-u^*$.

\paragraph{Gradient-dependent conductivity}
The argument through the tangent operator does not use the Kirchhoff transformation. For the regularized $p$-Laplacian
\begin{equation}\label{eq:nonlinear_plaplacian}
 -\rdiv\bigl(\pp\,\Psi_p(\nabla u)\bigr)=f,
 \qquad\Psi_p(F)=(1+|F|^2)^{(p-2)/2}F,\qquad p>1,
\end{equation}
the tangent weak operator is coercive on sets with bounded gradients, and the second derivative of $\Psi_p$ is bounded there. \Cref{lem:fosls_reduction} and a mean-value estimate then verify both assumptions, with explicit constants that are uniform for $p$ in a compact subinterval of $(1,\infty)$; see \ref{app:plaplacian}. Stability is thus proved as for the Kirchhoff law, but the estimate follows from \cref{thm:local_vc} and holds in a neighborhood of the branch; \cref{tab:example_summary} does not list this case.

\subsection{Compressible Neo-Hookean hyperelasticity}\label{sec:hyperelastic}
Let $u:\Omega\to\bR^d$ be the displacement, $F=I+\nabla u$ the deformation gradient, and
\begin{equation}\label{eq:neo_hookean_stress_note}
 P(F;\pp)=\mu(F-F^{-\top})+\lambda_L\log(\det F)F^{-\top}
\end{equation}
the first Piola--Kirchhoff stress. The Lam\'e parameters are $\mu=\pp/[2(1+\nu)]$ and $\lambda_L=\pp\nu/[(1+\nu)(1-2\nu)]$, where the parameter field $\pp$ is the Young's modulus and the Poisson ratio $\nu$ is fixed. Bounds $0<\pmin\le\pp\le\pmax$ and $0\le\nu\le\nu_{\max}<1/2$ give $0<\mu_{\min}\le\mu\le\mu_{\max}$ and $0\le\lambda_L\le\lambda_{L,\max}$. The equations are $\rdiv\sigma+f=0$ and $\sigma=P(F;\pp)$, with prescribed displacement on $\Gamma_D$ and traction on $\Gamma_N$. After lifting the boundary data, the residual is
\begin{equation}\label{eq:hyperelastic_residual}
 \cR_{\mathrm{hyp}}([\sigc,\uc];\pp)=
 \begin{pmatrix}\sigc+z-P(I+\nabla(\uc+w);\pp)\\ \rdiv\sigc+f_z\end{pmatrix},\qquad f_z=f+\rdiv z,
\end{equation}
and we work with $\Y=W^{1,\infty}(\Omega;\bR^d)\cap\U$. The tangent stress is
\begin{equation}\label{eq:neo_hookean_tangent_note}
 DP(F;\pp)[H]=\mu H+(\mu-\lambda_L\log\det F)F^{-\top}H^\top F^{-\top}+\lambda_L\operatorname{tr}(F^{-1}H)F^{-\top}.
\end{equation}
We verify the two assumptions in three steps. We first fix a neighborhood of the branch on which the deformation gradient stays bounded and its determinant stays positive. On this neighborhood, we bound the first two derivatives of the stress, which gives \cref{ass:lipschitz_derivative}. Finally, the stability of the tangent operator at the branch, which we assume, extends to the neighborhood and gives \cref{ass:local_stability} through \cref{lem:fosls_reduction}.

\paragraph{Admissible deformations}
We assume that $\norm{w}_{W^{1,\infty}}\le M_w$ and that the selected branch satisfies
\begin{equation}\label{eq:hyperelastic_branch_regularity}
 \sup_\pp\norm{u^{\circ,*}(\pp)}_{W^{1,\infty}}\le M_*,\qquad
 \det F_*(x;\pp)\ge j_*>0\ \ \text{a.e.}
\end{equation}
We choose $M_F\ge1$, $j_0>0$, and $R_0>0$ with
\begin{equation}\label{eq:hyperelastic_calibration}
 1+M_*+M_w+R_0\le M_F,\qquad R_0\le\frac{j_*}{2\sqrt d\,M_F^{d-1}},\qquad j_0\le\frac{j_*}{2}.
\end{equation}
Then every state in $\cN_\pp$ has $\norm{F}_{\mathrm{op}}\le M_F$ and $\det F\ge j_0$, and so does every point of the segment between two such states; we give the argument, based on $\norm{\operatorname{cof}F}_{\mathrm F}\le\sqrt d\,M_F^{d-1}$, in \ref{app:hyperelastic_derivatives}. A positive determinant expresses local orientation preservation; global injectivity of the deformation is not needed for the estimates. With
\begin{equation}\label{eq:hyperelastic_K_J}
 K_J=\frac{M_F^{d-1}}{j_0},\qquad L_J=\max\{|\log j_0|,\,d\log M_F\},
\end{equation}
we have $\norm{F^{-1}}_{\mathrm{op}}\le K_J$ and $|\log\det F|\le L_J$ on $\cN_\pp$, and the tangent and its derivative satisfy
\begin{align}
 \norm{DP(F;\pp)}&\le C^{(1)}:=\mu_{\max}(1+K_J^2)+\lambda_{L,\max}(d+L_J)K_J^2,
 \label{eq:hyperelastic_DP_op_bound}\\
 \norm{D^2P(F;\pp)}&\le C^{(2)}:=\bigl[2\mu_{\max}+\lambda_{L,\max}(2L_J+3\sqrt d)\bigr]K_J^3,
 \label{eq:hyperelastic_D2P_bound}
\end{align}
with Frobenius norms on arguments and values. Since $F$ is affine in $\bar u$ and segments stay admissible, the mean-value theorem gives $\norm{(DP(F_1;\pp)-DP(F_2;\pp))[\nabla v]}\le C^{(2)}\norm{\bar u_1-\bar u_2}_{W^{1,\infty}}\norm{v}_{H^1}$, so \cref{ass:lipschitz_derivative} holds with $\MK=C^{(2)}$.

\paragraph{Tangent stability}
Let $B_*:\U\to\U'$ be the tangent operator at the branch, $\langle B_*v,\varphi\rangle=\int_\Omega DP(F_*;\pp)[\nabla v]:\nabla\varphi\,dx$.

\begin{proposition}[Hyperelastic tangent stability]\label{prop:hyperelastic}
In addition to~\eqref{eq:hyperelastic_branch_regularity} and~\eqref{eq:hyperelastic_calibration}, assume one of the following.
\begin{enumerate}
\item[(a)] Coercivity at the branch: $\langle B_*v,v\rangle\ge a_*\norm{v}_{H^1}^2$ for $v\in\U$, with $a_*>0$ independent of $\pp$, and $C^{(2)}R_0\le a_*/2$.
\item[(b)] Invertibility at the branch: $\norm{B_*^{-1}}_{\U'\to\U}\le C_*$ uniformly in $\pp$, and $C_*C^{(2)}R_0<1$.
\end{enumerate}
Then \cref{ass:local_stability} holds, and \cref{lem:fosls_reduction} applies with $L_T=C^{(1)}$ and
\begin{equation}\label{eq:hyperelastic_inverse}
 C_I=2/a_*\ \text{in case (a)},\qquad C_I=\frac{C_*}{1-C_*C^{(2)}R_0}\ \text{in case (b)}.
\end{equation}
\end{proposition}
\begin{proof}
For every state in $\cN_\pp$ the tangent operator satisfies $\norm{B_{\bar u}-B_*}_{\U\to\U'}\le C^{(2)}R_0$. In case (a) this gives $\langle B_{\bar u}v,v\rangle\ge\tfrac12a_*\norm{v}_{H^1}^2$, and the Lax--Milgram lemma gives the first bound. In case (b) we write $B_{\bar u}=B_*(I+B_*^{-1}(B_{\bar u}-B_*))$ and use a Neumann series. \Cref{lem:fosls_reduction} with $T_{\bar u}v=DP(F;\pp)[\nabla v]$ completes the proof.
\end{proof}

Coercivity implies invertibility but not conversely: under increasing compressive load the tangent form can lose coercivity while the tangent remains invertible, so condition (b) is the relevant one for larger deformations. For small deformations (a) follows from Korn's inequality $\norm{\nabla v}^2\le C_{\mathrm{Korn}}\norm{\varepsilon(v)}^2$ on $\U$: at $F=I$ the tangent is the isotropic elasticity tensor $DP(I;\pp)[H]=\mu(H+H^\top)+\lambda_L\operatorname{tr}(H)I$, which is coercive with $a_0=2\mu_{\min}/(C_{\mathrm{Korn}}(1+C_P^2))$, and if the segment between $I$ and $F_*$ stays admissible, (a) holds with $a_*=a_0-C^{(2)}(M_*+M_w)$ whenever this number is positive. Larger deformations require an independent bound on the tangent inverse~\cite{marsdenHughes1983,ciarlet1988}, which the convergence of a nonlinear solver does not provide; the geometric and logarithmic terms in~\eqref{eq:neo_hookean_tangent_note} act on $\nabla v$ and belong to the principal symbol, so they cannot be treated as a compact perturbation of linear elasticity. Since $C^{(1)}$ and $C^{(2)}$ grow linearly in $\lambda_L$, the constants are not uniform in the incompressible limit $\nu\to1/2$.

With a valid $\clin$ and $\MK=C^{(2)}$, \cref{thm:local_vc} applies for every $R\le R_0$ with $\sqrt{\clin}>C^{(2)}R/2$. We have no analytic value of $a_*$ or $C_*$ for the benchmark, and we test the stability of the tangent at the computed states in \cref{subsec:exp_stability}.

\section{Finite element solution and reduced basis neural operator}\label{sec:rbno}
The estimates in \cref{sec:local_vc,sec:examples} are independent of the discretization and apply to admissible states in any conforming finite element space $\HH_h=\Sigma_h\times\U_h\subset\HH$. Examples include Raviart--Thomas or Brezzi--Douglas--Marini elements for the flux combined with continuous Lagrange elements for the primal variable. For elasticity, we use the flux elements row-wise to represent the generally nonsymmetric first Piola--Kirchhoff stress. We incorporate the homogeneous boundary conditions into these spaces.

In \cref{subsec:loss_excess}, we relate the loss excess above the finite element solution to the squared error relative to that solution, providing the basis for the numerical comparisons in \cref{sec:numerics}. We then introduce a discrete stability diagnostic in \cref{sec:stability_check} to examine \cref{ass:local_stability} at computed states. Finally, in \cref{subsec:pod,subsec:gradient}, we describe the RBNO construction of~\cite{QCD}, whose notation we follow, and its training through the finite element residual.

\subsection{Loss floor and loss excess}\label{subsec:loss_excess}
Since the exact solution is unavailable in computations, we assess the residual--error relation using the finite element solution $s_h\in\HH_h$ on the mesh used for residual evaluation as a reference. The error $E_h=\norm{s-s_h}_\HH$ differs from the exact error by at most $\varepsilon_h=\norm{s_h-s^*(\pp)}_\HH$, as in \cref{prop:consistency}. If $s_h$ satisfies the admissibility conditions of \cref{thm:local_vc}, then $\varepsilon_h\le\cL(s_h;\pp)^{1/2}/\sqrt{\Crel}$. Thus, $E_h$ closely approximates the exact error when it is large relative to this discretization-error bound.

The finite element solution generally has a nonzero residual. Its loss therefore remains even when a predictor reproduces $s_h$ exactly. For linear problems, Galerkin orthogonality gives the exact decomposition
\[
 \cL(s;\pp)
 =\cL(s_h;\pp)
 +\norm{\cB_\pp(s-s_h)}_{\LL^2}^2,
 \qquad s\in\HH_h,
\]
as shown in~\cite{QCD}. The first term is the minimum loss attainable in the finite element space, while the second measures the error relative to $s_h$. The following lemma establishes a local counterpart for nonlinear problems: under its assumptions, the loss excess $\cL(s;\pp)-\cL(s_h;\pp)$ is equivalent to $\norm{s-s_h}_\HH^2$ near $s_h$. This excess is the part of the loss that training can reduce within that neighborhood; reducing the finite element loss floor requires improving the approximation space.

\begin{lemma}[Loss excess]\label{lem:loss_excess}
Let \cref{ass:lipschitz_derivative} hold, with $\cK_\pp$ differentiable on $\cN_\pp$ as in \cref{ass:local_stability}, and let $s_h\in\HH_h$ be a stationary point of the loss on $\HH_h$,
\begin{equation}\label{eq:discrete_normal}
 \bigl(\cR(s_h;\pp),\cA_{s_h,\pp}t\bigr)_{\LL^2}=0,\qquad t\in\HH_h,
\end{equation}
where $\cA_{s_h,\pp}$ is the derivative of $\cR(\cdot;\pp)$ at $s_h$, that is, the operator~\eqref{Aubar} with $\bar u$ the primal component of $s_h$. Let $s\in\HH_h$, put $e_h=s-s_h$ with primal part $e_{h,u}$, and assume that the segment between the primal components of $s_h$ and $s$ lies in $\cN_\pp$. Then, with $\eta_h=\norm{\cR(s_h;\pp)}_{\LL^2}=\cL(s_h;\pp)^{1/2}$,
\begin{equation}\label{eq:loss_excess}
 \Bigl|\,\cL(s;\pp)-\cL(s_h;\pp)-\norm{\cA_{s_h,\pp}e_h}_{\LL^2}^2\Bigr|
 \le \MK\norm{e_{h,u}}_\Y\norm{e_h}_\HH\Bigl(\norm{\cA_{s_h,\pp}e_h}_{\LL^2}+\tfrac14\MK\norm{e_{h,u}}_\Y\norm{e_h}_\HH+\eta_h\Bigr).
\end{equation}
If, in addition, the stability bounds~\eqref{eq:uniform_linearized_stability} hold at the primal component of $s_h$, then
\begin{align*}
 \cL(s;\pp)-\cL(s_h;\pp)&\ge\Bigl(\sqrt{\clin}-\tfrac12\MK\norm{e_{h,u}}_\Y\Bigr)_+^2\norm{e_h}_\HH^2-\MK\eta_h\norm{e_{h,u}}_\Y\norm{e_h}_\HH,\\
 \cL(s;\pp)-\cL(s_h;\pp)&\le\Bigl(\sqrt{\Clin}+\tfrac12\MK\norm{e_{h,u}}_\Y\Bigr)^2\norm{e_h}_\HH^2+\MK\eta_h\norm{e_{h,u}}_\Y\norm{e_h}_\HH .
\end{align*}
\end{lemma}
\begin{proof}
As in~\eqref{eq:residual_split}, $\cR(s;\pp)=\cR(s_h;\pp)+\cA_{s_h,\pp}e_h+\omega$ with $\norm{\omega}_{\LL^2}\le\tfrac12\MK\norm{e_{h,u}}_\Y\norm{e_h}_\HH$. Expanding $\norm{\cR(s;\pp)}_{\LL^2}^2$ and using~\eqref{eq:discrete_normal} with $t=e_h$ gives
\[
 \cL(s;\pp)-\cL(s_h;\pp)=\norm{\cA_{s_h,\pp}e_h+\omega}_{\LL^2}^2+2\bigl(\cR(s_h;\pp),\omega\bigr)_{\LL^2}.
\]
Expanding the square and bounding the terms that contain $\omega$ by the Cauchy--Schwarz inequality gives~\eqref{eq:loss_excess}. For the two-sided bounds we use instead
\[
 \bigl|\norm{\cA_{s_h,\pp}e_h+\omega}_{\LL^2}-\norm{\cA_{s_h,\pp}e_h}_{\LL^2}\bigr|\le\tfrac12\MK\norm{e_{h,u}}_\Y\norm{e_h}_\HH,\qquad
 2\bigl|\bigl(\cR(s_h;\pp),\omega\bigr)_{\LL^2}\bigr|\le\MK\eta_h\norm{e_{h,u}}_\Y\norm{e_h}_\HH,
\]
together with~\eqref{eq:uniform_linearized_stability} at the primal component of $s_h$.
\end{proof}

For a linear residual,~\eqref{eq:discrete_normal} gives the normal equations of the least-squares problem on $\HH_h$. For nonlinear problems, we compute $s_h$ by Gauss--Newton iteration,
\[
 s_{k+1}=s_k+\alpha_k\delta_k,
 \qquad \alpha_k\in(0,1],
\]
where the correction $\delta_k\in\HH_h$ solves
\[
 (\cA_k\delta_k,\cA_kt)_{\LL^2}
 =-(\cR(s_k;\pp),\cA_kt)_{\LL^2},
 \qquad \forall\,t\in\HH_h,
 \qquad \cA_k=\cA_{s_k,\pp}.
\]
The fixed points of this iteration satisfy~\eqref{eq:discrete_normal}. Conversely, where the linearization is stable, the correction is unique, and every solution of~\eqref{eq:discrete_normal} is a fixed point. In the full space $\HH$, the linearized operators considered in \cref{sec:examples} are surjective onto $\LL^2$, so the Gauss--Newton correction coincides with the Newton correction for $\cR(s;\pp)=0$. On $\HH_h$, the correction instead solves a discrete least-squares problem and generally does not cancel the linearized residual exactly.

When we terminate the iteration at a prescribed tolerance, the computed state satisfies~\eqref{eq:discrete_normal} only approximately. We measure this stationarity defect by
\[
 \rho_h(t)
 =(\cR(s_h;\pp),\cA_{s_h,\pp}t)_{\LL^2}.
\]
The estimate~\eqref{eq:loss_excess} then acquires the additional term
\[
 2|\rho_h(e_h)|
 \le 2\norm{\rho_h}_{\HH_h'}\norm{e_h}_\HH
\]
on its right-hand side.

When $\Y=\U$, the perturbation terms in~\eqref{eq:loss_excess}, relative to $\norm{e_h}_\HH^2$, are of order $O(\norm{e_h}_\HH+\eta_h)$. For sufficiently small error and residual floor, the loss excess is therefore equivalent to $\norm{e_h}_\HH^2$, with constants close to those of the linearized problem. For gradient-dependent models, we use the inverse inequality
\[
 \norm{e_{h,u}}_\Y\le C_h\norm{e_h}_\HH,
\]
where $C_h$ grows under mesh refinement. The relative perturbation is then of order $O(C_h(\norm{e_h}_\HH+\eta_h))$, provided that $C_h\norm{e_h}_\HH\le1$. For the semilinear model, we obtain an exact version on bounded sets.

\begin{remark}[Loss excess for the semilinear model]\label{rem:semilinear_excess}
Let $s,s_h\in\HH_h$, with $s_h$ satisfying~\eqref{eq:discrete_normal}, and let $e_h=s-s_h=[e_{h,\sigma},e_{h,u}]$, with physical primal components $u$ and $u_h$. Since $u^3-u_h^3=(u^2+uu_h+u_h^2)e_{h,u}$, the residual satisfies exactly $\cR_{\mathrm{dr}}(s)=\cR_{\mathrm{dr}}(s_h)+\cA_{\mathrm{sec}}e_h$ with the linear operator
\[
 \cA_{\mathrm{sec}}e_h=\bigl(e_{h,\sigma}-\pp\nabla e_{h,u},\,-\rdiv e_{h,\sigma}+\kappa(u^2+uu_h+u_h^2)e_{h,u}\bigr)
\]
of the proof of \cref{prop:semilinear}, with $u_h$ in place of $u^*$. It differs from the derivative $\cA_{s_h,\pp}$ only in the reaction coefficient, by $\kappa(u^2+uu_h-2u_h^2)=\kappa e_{h,u}(u+2u_h)$. Hence~\eqref{eq:discrete_normal} with $t=e_h$ gives the exact relation
\[
 \cL_{\mathrm{dr}}(s)-\cL_{\mathrm{dr}}(s_h)=\norm{\cA_{\mathrm{sec}}e_h}_{\LL^2}^2+2\bigl(r_2(s_h),\kappa\,e_{h,u}^2(u+2u_h)\bigr),
\]
where $r_2(s_h)$ is the second component of $\cR_{\mathrm{dr}}(s_h)$. The last term is bounded by $\delta_h\norm{\nabla e_{h,u}}^2$ with $\delta_h=2\kappa_{\max}C_S^2\norm{u+2u_h}_{L^6}\eta_h$. Hence the loss excess is equivalent to $\norm{e_h}_\HH^2$ with the constants of \cref{prop:semilinear}, perturbed by $\delta_h$, on every bounded set of states.
\end{remark}

In \cref{subsec:numerical_metrics,subsec:exp_ratio}, we compare the squared error $\norm{s-s_h}_\HH^2$ with the loss excess through the ratio $\chi$ in~\eqref{eq:numerical_metrics}.

\subsection{Discrete stability constant}\label{sec:stability_check}
We examine the stability of the linearized residual on $\HH_h$ through a generalized eigenvalue problem. Let $\{\psi_j\}_{j=1}^{N_h^s}$ be a basis of $\HH_h$, and define
\[
 (X_h)_{ij}=(\psi_j,\psi_i)_\HH,\qquad
 (G_h(\bar u,\pp))_{ij}
 =(\cA_{\bar u,\pp}\psi_j,\cA_{\bar u,\pp}\psi_i)_{\LL^2}.
\]
With exact integration, the smallest generalized eigenvalue satisfies
\begin{equation}\label{eq:discrete_eigen_check}
 \lambda_{\min,h}(\bar u,\pp)
 =\min_{\bst\ne0}\frac{\bst^\top G_h\bst}{\bst^\top X_h\bst}
 =\inf_{0\ne t_h\in\HH_h}
 \frac{\norm{\cA_{\bar u,\pp}t_h}_{\LL^2}^2}{\norm{t_h}_\HH^2},
\end{equation}
where $\bst$ is the coefficient vector of $t_h$. Thus, $\lambda_{\min,h}$ is the coercivity constant of the Gauss--Newton bilinear form
\[
 (t,t')\mapsto
 (\cA_{\bar u,\pp}t,\cA_{\bar u,\pp}t')_{\LL^2}
\]
on $\HH_h$. Equivalently, its square root is the best constant in the lower bound
$\norm{\cA_{\bar u,\pp}t_h}_{\LL^2}\ge
\lambda_{\min,h}^{1/2}\norm{t_h}_\HH$.

Since $\HH_h\subset\HH$, we have
$\lambda_{\min,h}(\bar u,\pp)\ge\lambda_{\min}(\bar u,\pp)$,
where $\lambda_{\min}$ denotes the corresponding infimum over $\HH$. For a fixed linearization point and parameter, the discrete eigenvalue is nonincreasing under nested refinement and converges to $\lambda_{\min}$ if the approximation spaces are dense in $\HH$, as shown in \ref{app:diagnostic}. It therefore provides an upper bound on the continuous coercivity constant and cannot certify the lower stability bound in~\eqref{eq:uniform_linearized_stability}. To isolate the effect of mesh refinement, we must keep both the parameter and the linearization point fixed, for example by representing the same coarse-mesh state on each finer mesh. If we instead linearize at the finite element solution computed on each mesh, changes in $\lambda_{\min,h}$ reflect changes in both the space and the state.

Numerical integration introduces a further dependence on the quadrature rules used to assemble the matrices. If we project the residual onto a discrete space before evaluating its norm, the eigenvalue instead measures stability of the projected operator. Subject to these qualifications, the discrete eigenvalue provides a diagnostic of lower stability along solver iterates and at network predictions. A small value may indicate a nearly singular linearization, loss of information through residual projection, or underintegration. In \cref{subsec:exp_stability}, we use this diagnostic to examine \cref{ass:local_stability} at the computed states, particularly for compressible hyperelasticity, where we do not have an analytic lower bound on the branch tangent.

\subsection{Reduced basis neural operator}\label{subsec:pod}
From snapshots $s_h(\pp_i)$, $i=1,\dots,N_{\mathrm{POD}}$, we compute an $\HH$-orthonormal basis by proper orthogonal decomposition (POD), as in~\cite{QCD}. Let $\bss_h(\pp_i)\in\bR^{N_h^s}$ be the coefficient vectors of the snapshots, $X_h$ the Gram matrix of the $\HH$ inner product, and $\mathbf S=(\bss_h(\pp_1),\dots,\bss_h(\pp_{N_{\mathrm{POD}}}))$ the snapshot matrix. With the eigenpairs of the correlation matrix,
\[
 \mathbf C\bsv_k=\lambda_k\bsv_k,\qquad \mathbf C=\frac1{N_{\mathrm{POD}}}\mathbf S^\top X_h\mathbf S,\qquad \lambda_1\ge\lambda_2\ge\cdots,
\]
and orthonormal eigenvectors $\bsv_k$, the vectors $\pi_k=\mathbf S\bsv_k/\sqrt{N_{\mathrm{POD}}\lambda_k}$ are orthonormal in the $X_h$ inner product, and $\phi_k\in\HH_h$ is the function with coefficient vector $\pi_k$. We set $\Phi_r=(\phi_1,\dots,\phi_r)$, $\HH_r=\operatorname{span}\{\phi_k\}_{k=1}^r$, and $\Pi_r=(\pi_1,\dots,\pi_r)\in\bR^{N_h^s\times r}$, so that a function $\Phi_r\bss_r\in\HH_r$ has the finite element coefficients $\Pi_r\bss_r$. The $\HH$-orthogonal projector $P_r$ onto $\HH_r$ satisfies $P_rs_h=\Phi_r\Pi_r^\top X_h\bss_h$, and
\begin{equation}\label{eq:empirical_pod_identity}
 \frac1{N_{\mathrm{POD}}}\sum_{i=1}^{N_{\mathrm{POD}}}\norm{(I-P_r)s_h(\pp_i)}_\HH^2
 =\sum_{k>r}\lambda_k
 =\frac1{N_{\mathrm{POD}}}\sum_i\norm{s_h(\pp_i)}_\HH^2-\sum_{k=1}^r\lambda_k .
\end{equation}
The last expression gives the tail without computing the trailing modes. The identity concerns the snapshots. We measure the projection error at other parameters on the test set in \cref{sec:numerics} and refer to~\cite{cohenDeVore2015} for the approximation of parametric solution manifolds in general.

The reduced basis neural operator is
\begin{equation}\label{eq:nonlinear_hypothesis}
 s_r(\pp;\theta)=\Phi_r\,\bss_r(\pp;\theta)=\sum_{k=1}^r\bigl(\bss_r(\pp;\theta)\bigr)_k\,\phi_k,\qquad \bss_r(\pp;\theta)\in\bR^r,
\end{equation}
where the coefficient vector $\bss_r(\pp;\theta)$ is the output of a neural network with parameters $\theta$. Its distance from the finite element solution $s_h=s_h(\pp)$ splits exactly into a part due to the basis and a part due to the coefficients,
\begin{equation}\label{eq:state_pythagoras}
 \norm{s_r(\pp;\theta)-s_h}_\HH^2=\norm{(I-P_r)s_h}_\HH^2+\norm{s_r(\pp;\theta)-P_rs_h}_\HH^2 ,
\end{equation}
and the network can reduce only the second part.

For a strongly measurable candidate $s\in\X$, we define the expected and the empirical loss
\begin{equation}\label{eq:ideal_nonlinear_loss}
 \cL_\mu(s)=\int_\PP\cL(s(\pp);\pp)\,d\mu(\pp),\qquad
 \cL_N(s)=\frac1N\sum_{i=1}^N\cL(s(\pp_i);\pp_i).
\end{equation}
If $s(\pp)\in\cS_\pp(R)$ (see~\eqref{eq:admissible_set}) for almost every $\pp$, integrating~\eqref{eq:local_vc_unsquared} gives
\begin{equation}\label{eq:lifted_local_vc}
 \Crel\norm{s-s^*}_\X^2\le\cL_\mu(s)\le\Ceff\norm{s-s^*}_\X^2 ,
\end{equation}
and the same argument applies to the empirical loss when every sampled prediction is admissible. We train RBNO by minimizing
\begin{equation}\label{eq:nonlinear_training}
 \cL_N\bigl(s_r(\cdot;\theta)\bigr)=\frac1N\sum_{i=1}^N\norm{\cR(s_r(\pp_i;\theta);\pp_i)}_{\LL^2}^2
\end{equation}
over $\theta$, where $\pp_1,\dots,\pp_N$ are training parameters for which we compute no solutions. We need snapshots only to construct $\HH_r$ and to validate the trained model; training itself requires only residual evaluations.

For the reduced minimization problem, we fix $\rho>0$ and use the bounded set
\begin{equation}\label{eq:tube_T}
 \cT_\pp=\{s\in\cS_\pp(R):\norm{s-s^*(\pp)}_\HH\le\rho\}.
\end{equation}
If $\HH_r\subset\Sigma\times\Y$, as for the finite element spaces used here, whose primal components are continuous piecewise polynomials, the norm of $\Y$ is continuous on the finite-dimensional space $\HH_r$. The intersection $\HH_r\cap\cT_\pp$ is then closed and bounded in $\HH_r$, hence compact, and no compact embedding of $\Y$ is needed; this gives the existence of a reduced minimizer in the following proposition.

\begin{proposition}[RBNO error and reduced quasi-optimality]\label{prop:rbno_error}
Let the assumptions of \cref{thm:local_vc} hold. Every strongly measurable $s_r(\cdot;\theta)\in\X$ with $s_r(\pp;\theta)\in\cS_\pp(R)$ for almost every $\pp$ satisfies
\begin{equation}\label{eq:rbno_certified_error}
 \Crel\norm{s_r(\cdot;\theta)-s^*}_\X^2\le\cL_\mu\bigl(s_r(\cdot;\theta)\bigr)\le\Ceff\norm{s_r(\cdot;\theta)-s^*}_\X^2 .
\end{equation}
If, for a fixed parameter, $\HH_r\subset\Sigma\times\Y$ and $\cT_{\pp,r}:=\HH_r\cap\cT_\pp\ne\emptyset$, then a minimizer $s_r^\dagger(\pp)$ of $\cL(\cdot;\pp)$ over $\cT_{\pp,r}$ exists and
\begin{equation}\label{eq:rbno_quasioptimal}
 \norm{s_r^\dagger(\pp)-s^*(\pp)}_\HH\le\sqrt{\Ceff/\Crel}\,\inf_{t\in \cT_{\pp,r}}\norm{t-s^*(\pp)}_\HH .
\end{equation}
\end{proposition}
\begin{proof}
The first statement is~\eqref{eq:lifted_local_vc}. The set $\cT_{\pp,r}$ is compact, as noted before the proposition, and $\cL(\cdot;\pp)$ is continuous on it, so it attains its minimum there. For every $t\in \cT_{\pp,r}$ we have $\Crel\norm{s_r^\dagger-s^*}_\HH^2\le\cL(s_r^\dagger)\le\cL(t)\le\Ceff\norm{t-s^*}_\HH^2$.
\end{proof}

When $P_rs^*(\pp)\in\cT_\pp$, which for $\Y\ne\U$ is a condition in the norm of $\Y$, the infimum in~\eqref{eq:rbno_quasioptimal} is at most $\norm{(I-P_r)s^*(\pp)}_\HH$. In the linear case the minimizer over $\HH_r$ is unique~\cite{QCD}. Here $s_r^\dagger(\pp)$ need not be unique, and integrating it over the parameter set requires a measurable selection, which the estimate~\eqref{eq:rbno_certified_error} for the network does not need. An optimizer applied to the reduced loss, for instance one started from the projection of a finite element solution, returns a stationary point, which need not be $s_r^\dagger$.

Training minimizes the empirical loss. Let $\mathcal F_r$ be a class of admissible reduced predictors, and let $\varepsilon_{\mathrm{gen}}=\sup_{s\in\mathcal F_r}|\cL_\mu(s)-\cL_N(s)|$ be the generalization gap. If $\widehat s\in\mathcal F_r$ minimizes the empirical loss up to $\varepsilon_{\mathrm{opt}}$, then~\eqref{eq:rbno_certified_error} gives, for every $s\in\mathcal F_r$,
\begin{align*}
 \Crel\norm{\widehat s-s^*}_\X^2&\le\cL_\mu(\widehat s)\le\cL_N(\widehat s)+\varepsilon_{\mathrm{gen}}\le\cL_N(s)+\varepsilon_{\mathrm{opt}}+\varepsilon_{\mathrm{gen}}\\
 &\le\cL_\mu(s)+\varepsilon_{\mathrm{opt}}+2\varepsilon_{\mathrm{gen}}\le\Ceff\norm{s-s^*}_\X^2+\varepsilon_{\mathrm{opt}}+2\varepsilon_{\mathrm{gen}},
\end{align*}
and hence
\begin{equation}\label{eq:learning_bound}
 \norm{\widehat s-s^*}_\X^2\le\frac{\Ceff}{\Crel}\inf_{s\in\mathcal F_r}\norm{s-s^*}_\X^2+\frac{2\varepsilon_{\mathrm{gen}}+\varepsilon_{\mathrm{opt}}}{\Crel}.
\end{equation}
Bounds for $\varepsilon_{\mathrm{gen}}$ in the linear case are given in~\cite{QCD}; the optimization error is outside the scope of the PDE analysis.

\subsection{Training through the finite element residual}\label{subsec:gradient}
In the linear case, the loss of a reduced state is a quadratic form in the reduced coefficients, with matrices assembled once for each training parameter~\cite{QCD}. For nonlinear problems, this representation is generally unavailable. With a cubic reaction term, the loss is a polynomial of degree six in the reduced coefficients. A direct evaluation using a preassembled dense sixth-order tensor would require $O(r^6)$ operations per parameter, which is prohibitive for $r=128$. For the Neo-Hookean law, the loss is nonpolynomial. We therefore retain finite element assembly in the training loop to evaluate the loss and its derivative at each predicted state.

Let $\bss\in\bR^{N_h^s}$ denote the coefficient vector of $s\in\HH_h$ in the basis $\{\psi_j\}$ introduced in \cref{sec:stability_check}, and write $J(\bss;\pp)=J_{h,q}(s;\pp)$. For RBNO, these coefficients are given by $\bss=\Pi_r\bss_r(\pp;\theta)$. Using symbolic differentiation in the form language~\cite{alnaes2012ufl}, we obtain
\begin{equation}\label{eq:state_gradient}
 DJ_{h,q}(s;\pp)[t]
 =2\bigl(\cR(s;\pp),\cA_st\bigr)_{\LL^2},
 \qquad t\in\HH_h,
\end{equation}
where $\cA_s$ is the derivative of $\cR(\cdot;\pp)$ at $s$, and we evaluate the inner product using the quadrature rule of the loss. The form compiler generates the element code for this linear form once, before training~\cite{logg2012automated,Baratta_DOLFINx_the_next_2023}. At each predicted state, we assemble the form as a load vector to obtain
\begin{equation}\label{eq:assembled_gradient}
 \bsg(\bss;\pp)
 =\bigl(DJ_{h,q}(s;\pp)[\psi_j]\bigr)_{j=1}^{N_h^s}
 \in\bR^{N_h^s}.
\end{equation}
Since $s$ depends linearly on $\bss$, this vector is the Euclidean gradient of $J$ with respect to the finite element coefficients. Using the same quadrature rule for the loss and its derivative ensures that $\bsg$ differentiates the evaluated loss $J_{h,q}$, up to floating-point errors. We verify this gradient using directional finite differences.

\begin{figure}[tb]
\centering
\includegraphics[width=\linewidth]{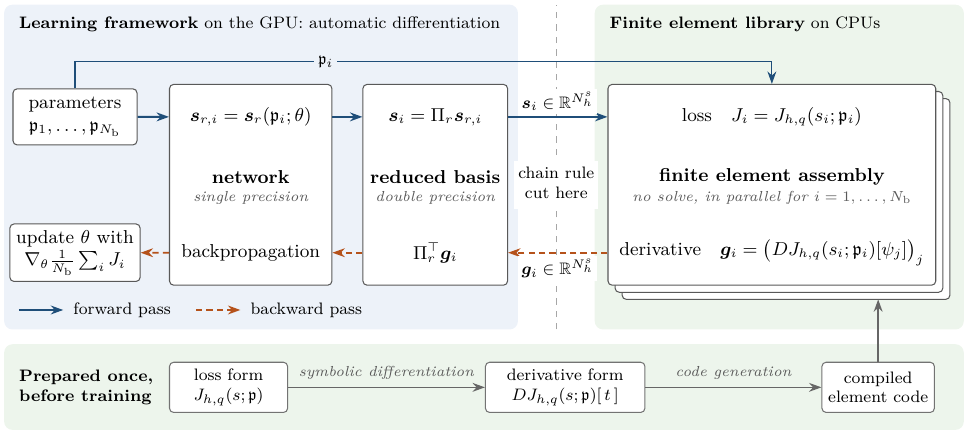}
\caption{Computational graph for the training step in \cref{alg:training}. Solid and dashed arrows indicate the forward and backward passes, respectively. The learning framework differentiates the network and the reduced basis map, while the finite element library assembles the loss and coefficient gradient for each parameter without solving a PDE. For each batch, the learning framework passes the state coefficients $\bss_i$ to the finite element evaluator, which returns the losses $J_i$ and gradients $\bsg_i$. Bottom: before training, we obtain the derivative form~\eqref{eq:state_gradient} by symbolic differentiation and compile it into element code.}\label{fig:training_pipeline}
\end{figure}

For one parameter, the chain rule gives
\begin{equation}\label{eq:chain_rule}
 \nabla_\theta J\bigl(\Pi_r\bss_r(\pp;\theta);\pp\bigr)
 =
 \bigl(D_\theta\bss_r(\pp;\theta)\bigr)^\top
 \Pi_r^\top
 \bsg\bigl(\Pi_r\bss_r(\pp;\theta);\pp\bigr).
\end{equation}
The gradient of the training objective~\eqref{eq:nonlinear_training} is the average of this expression over the training parameters. In each training step, we approximate this average using a batch of $N_{\mathrm b}$ parameters. The finite element code returns the scalar loss $J$ and the coefficient gradient $\bsg$, and the learning framework propagates $\Pi_r^\top\bsg$ through the network. Thus, finite element assembly acts as a differentiable operation in the computational graph, with its derivative supplied by the finite element library~\cite{bouziani2023physics,mitusch2019dolfin}. We summarize one training step in \cref{alg:training} and illustrate the division of work in \cref{fig:training_pipeline}.

\begin{algorithm}[tb]
\caption{One RBNO training step using the finite element residual}\label{alg:training}
\begin{algorithmic}[1]
\Require Batch $\pp_1,\dots,\pp_{N_{\mathrm b}}$ drawn from the $N$ training parameters; network $\bss_r(\cdot;\theta)$; POD matrix $\Pi_r$
\State $\bss_{r,i}\gets\bss_r(\pp_i;\theta)$ for $i=1,\dots,N_{\mathrm b}$ \Comment{GPU, single precision}
\State $\bss_i\gets\Pi_r\bss_{r,i}$ for $i=1,\dots,N_{\mathrm b}$ \Comment{reconstruction in double precision}
\For{$i=1,\dots,N_{\mathrm b}$, in parallel on CPU workers}
  \State Assemble $J_i=J(\bss_i;\pp_i)$ and $\bsg_i=\bsg(\bss_i;\pp_i)$ using~\eqref{eq:state_gradient} and~\eqref{eq:assembled_gradient}
\EndFor
\State Backpropagate $\tfrac{1}{N_{\mathrm b}}\Pi_r^\top\bsg_i$ through each network output to obtain $\nabla_\theta\bigl(\tfrac{1}{N_{\mathrm b}}\sum_iJ_i\bigr)$
\State Update $\theta$ using the optimizer
\end{algorithmic}
\end{algorithm}

For each parameter, a training step requires one assembly of the loss and one assembly of the linear form~\eqref{eq:state_gradient}. This cost is comparable to assembling the right-hand side of a Gauss--Newton step, without matrix assembly or factorization. We report the measured times in \cref{tab:cost}. Since these evaluations are independent across the batch, we distribute them among CPU workers while keeping the network on the GPU. Hyper-reduction~\cite{chaturantabut2010deim,farhat2014ecsw} could reduce the evaluation cost further, with the resulting approximation error accounted for through \cref{prop:consistency}.

No Gram matrix is needed in the chain rule. The vector $X_h^{-1}\bsg$ would represent the derivative through the $\HH$ inner product, whereas backpropagation requires the coefficient gradient $\bsg$ itself. Since $\HH_h$ incorporates the homogeneous boundary conditions, we need no boundary penalties or separately chosen collocation points, in contrast to physics-informed training with pointwise residuals. The same evaluator supplies both the training loss and the residual used to assess prediction accuracy under the stability and consistency assumptions. As we show in \cref{subsec:exp_baselines}, inaccurate grid-based evaluation can compromise this error assessment. The procedure applies to any predictor with outputs in a conforming finite element space, replacing $\Pi_r$ by the corresponding map from network outputs to finite element coefficients and differentiating that map through the chain rule.

\section{Numerical experiments}\label{sec:numerics}
We examine the theory on one benchmark for each problem class in \cref{sec:examples}. We describe the benchmarks, finite element discretizations, and error measures in \cref{subsec:problem_configs,subsec:numerical_metrics}, and report the accuracy and computational cost of RBNO in \cref{subsec:exp_accuracy}. In \cref{subsec:exp_ratio}, we compare the error with the loss excess to assess the estimates of \cref{lem:loss_excess}. In \cref{subsec:exp_stability}, we examine \cref{ass:local_stability} through discrete stability constants. We then compare RBNO with PCA-Net, FNO, and FC-PINO in a common finite element space in \cref{subsec:exp_baselines}, showing that small grid-based residuals can coexist with large PDE-compliant norm errors. We provide additional results on reduced basis approximation, model architectures and training settings, training histories, RBNO predictions, and pointwise residuals in \ref{app:numerical_results}.

We perform all computations on a Linux server with two AMD EPYC 7543 CPUs and an NVIDIA RTX A6000 GPU. We use DOLFINx~0.10.0~\cite{Baratta_DOLFINx_the_next_2023} with Python~3.11 for finite element discretization, nonlinear solves, and loss and gradient assembly. We sample Gaussian random fields using an adaptation of hIPPYlib~\cite{villa2021hippylib} to DOLFINx, which realizes the covariance operator through an elliptic problem with Robin boundary conditions. We implement the networks in PyTorch~2.9.1~\cite{paszke2019pytorch} and use neuraloperator~2.0.0~\cite{kossaifi2025librarylearningneuraloperators} for the FNO and FC-PINO baselines. The code is available at \url{https://github.com/qy849/vcno-rbno-nonlinear}.

\subsection{Benchmarks}\label{subsec:problem_configs}
\begin{figure}[!htb]
\centering
\includegraphics[width=.24\linewidth]{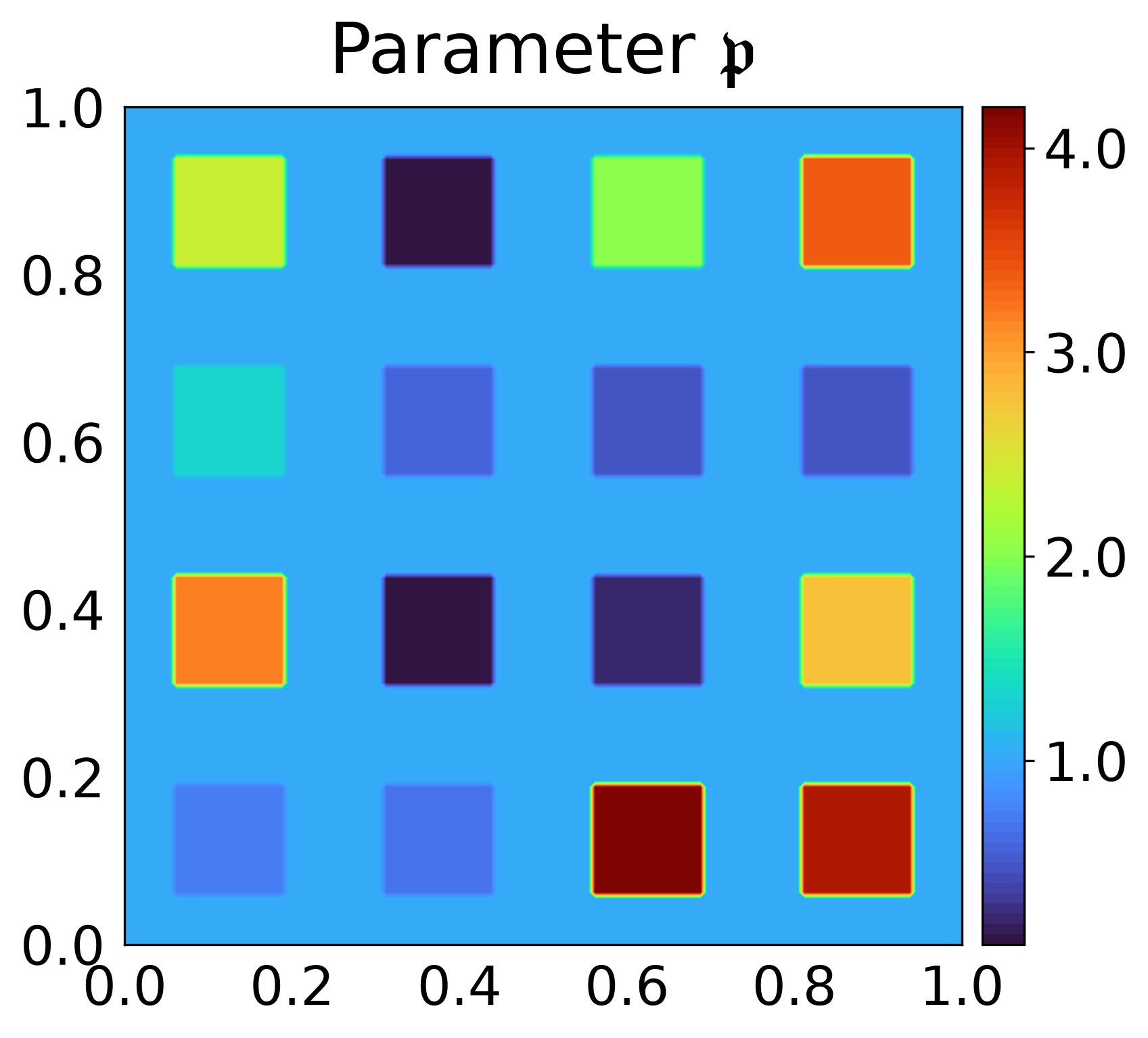}
\includegraphics[width=.24\linewidth]{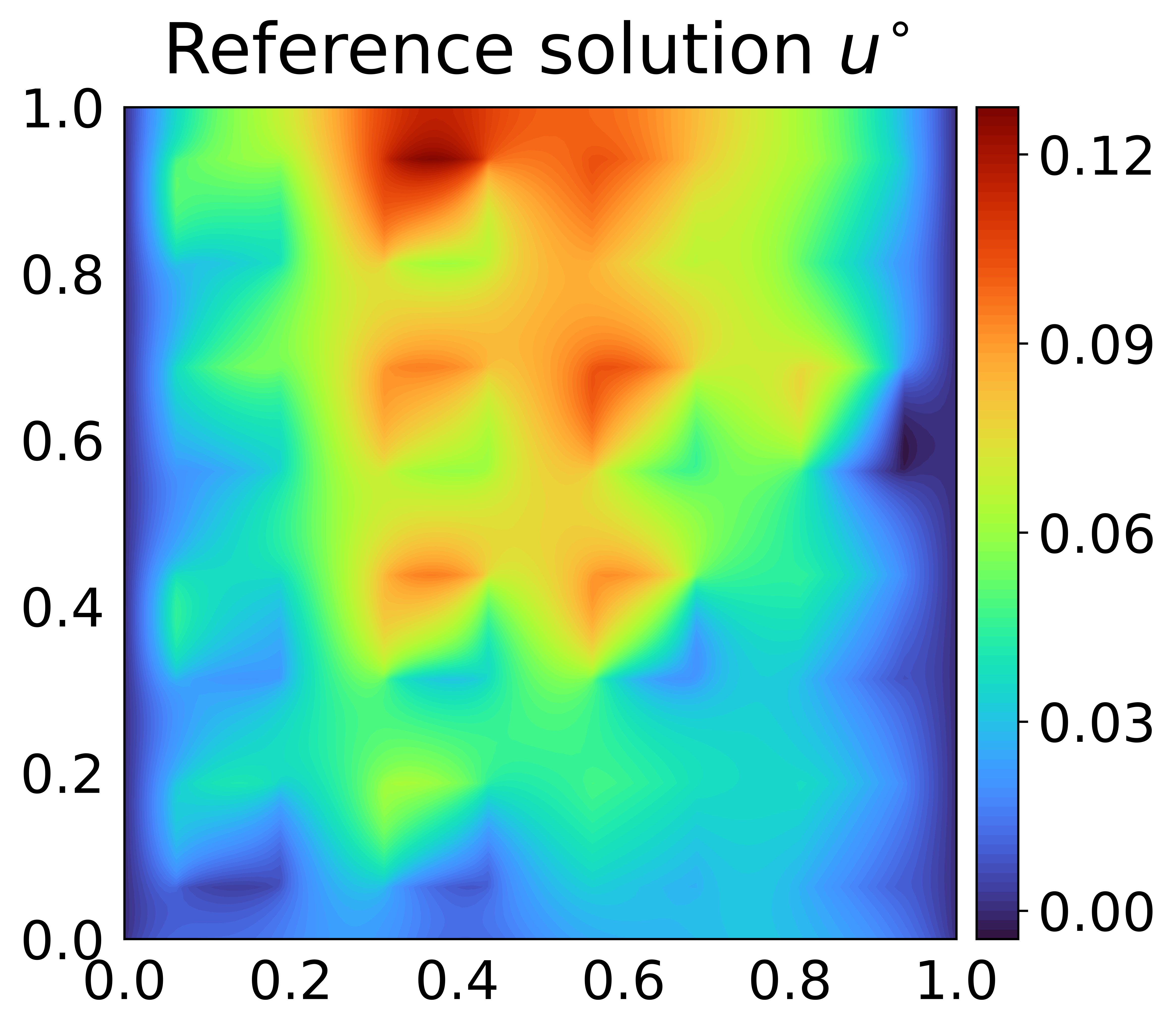}
\includegraphics[width=.24\linewidth]{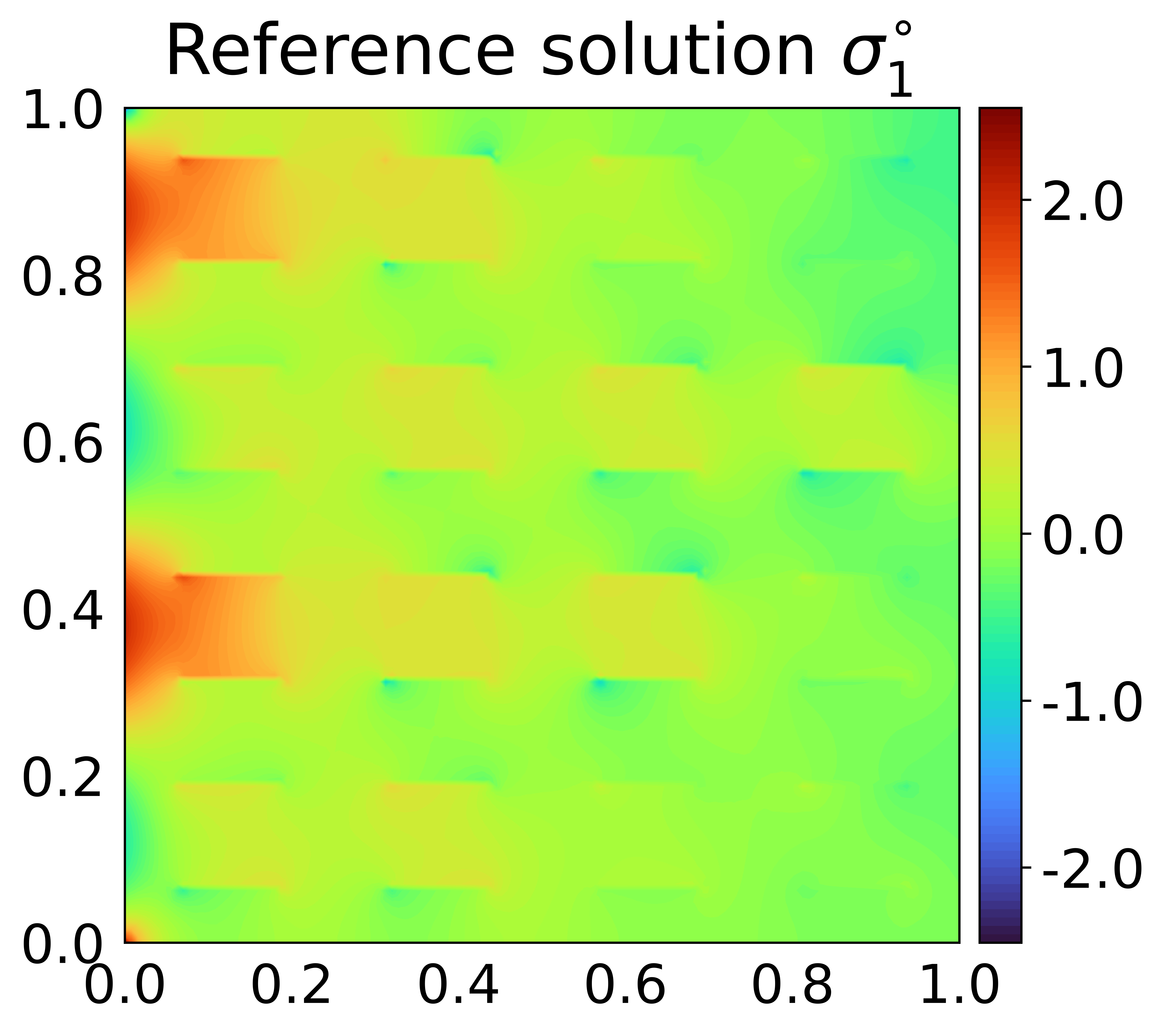}
\includegraphics[width=.24\linewidth]{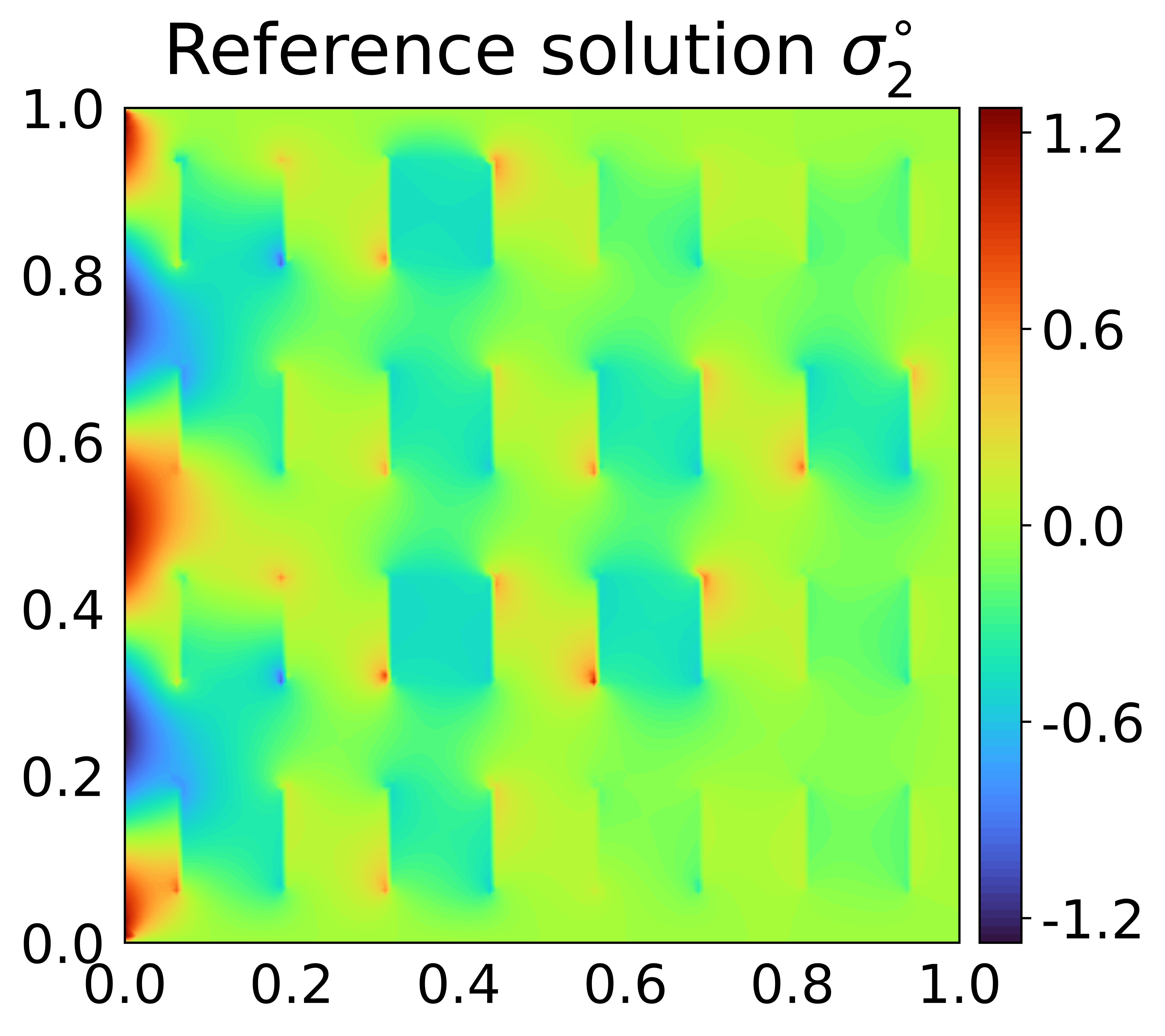}

\vspace{1em}
\includegraphics[width=.24\linewidth]{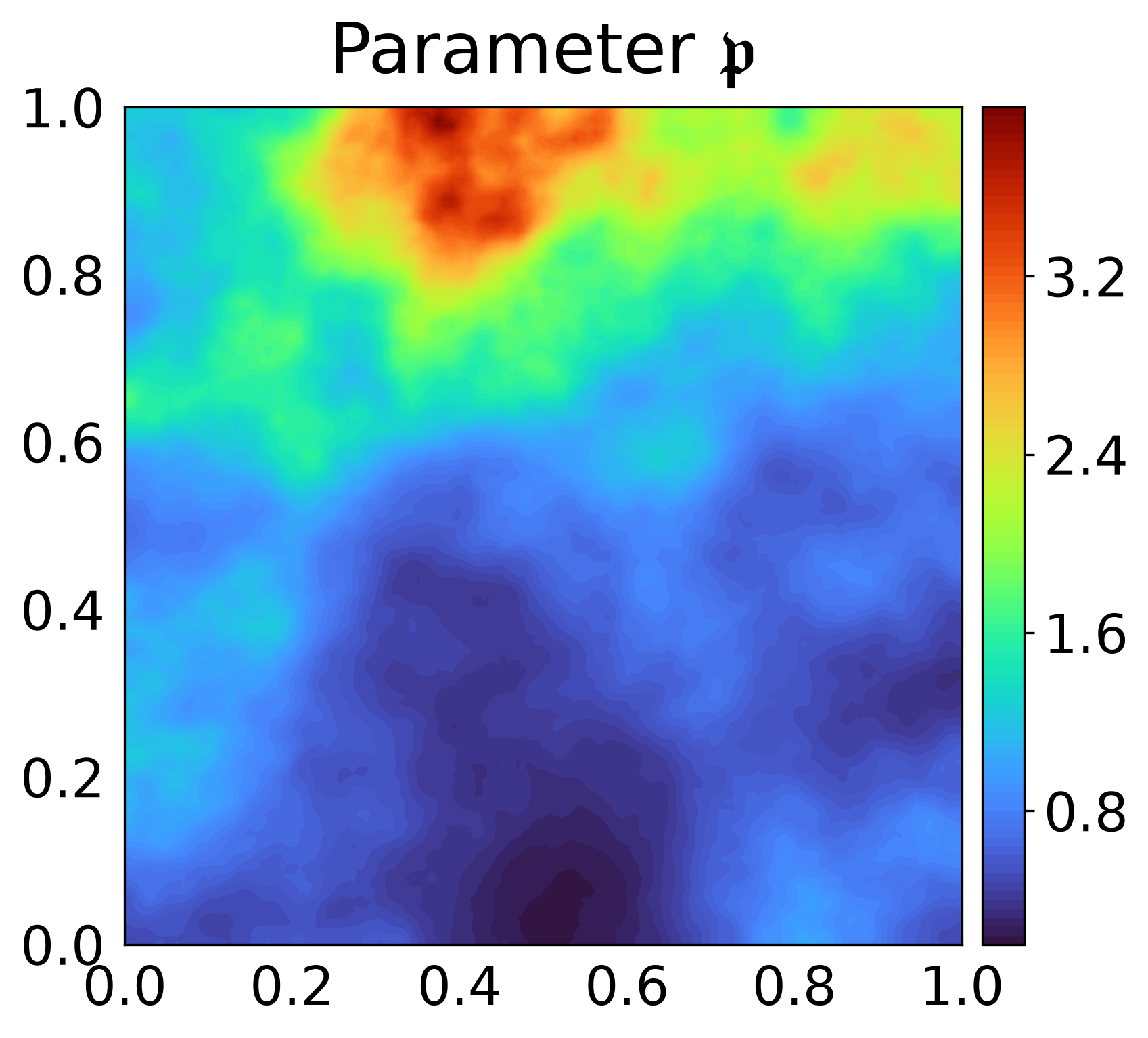}
\includegraphics[width=.24\linewidth]{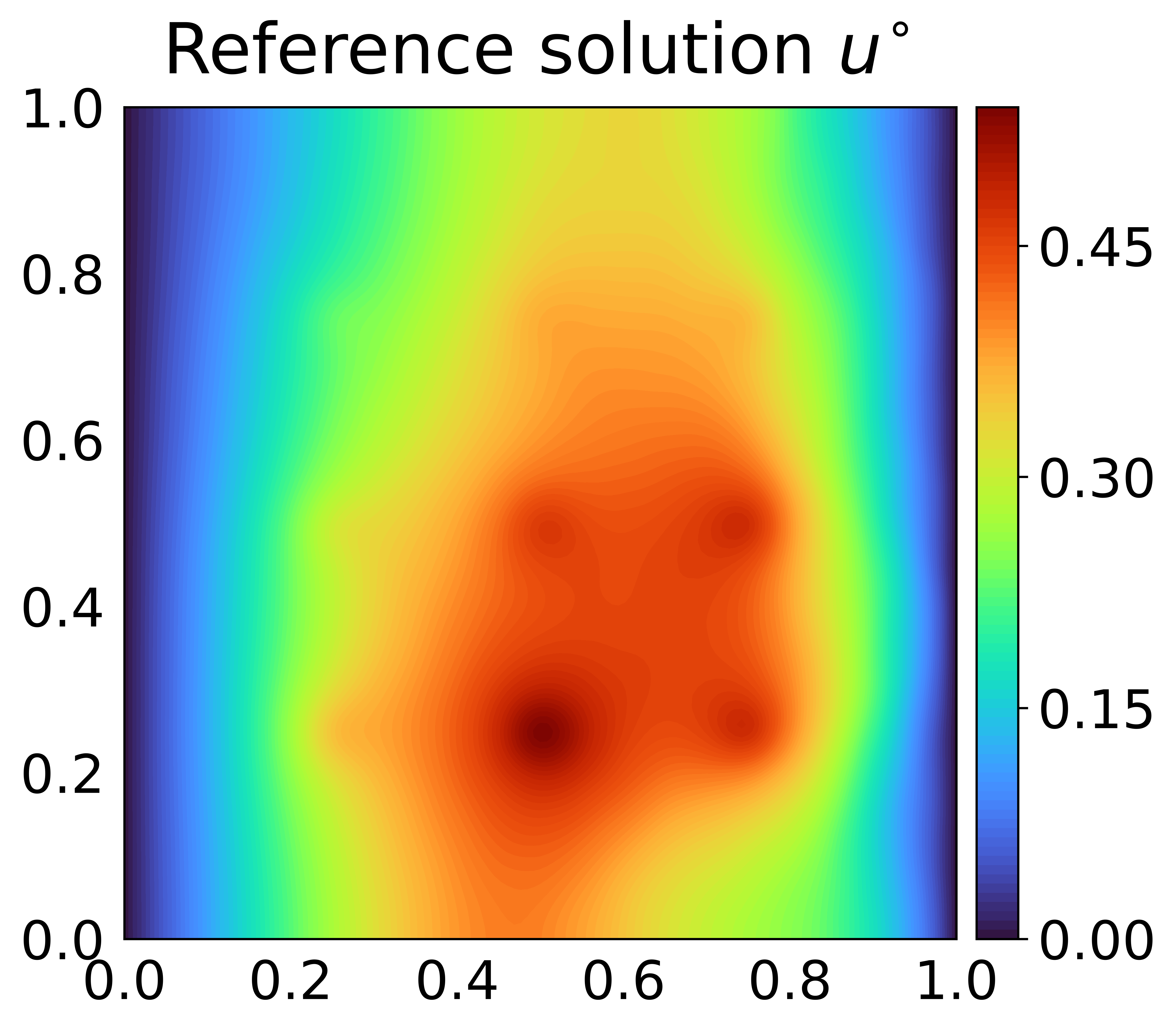}
\includegraphics[width=.24\linewidth]{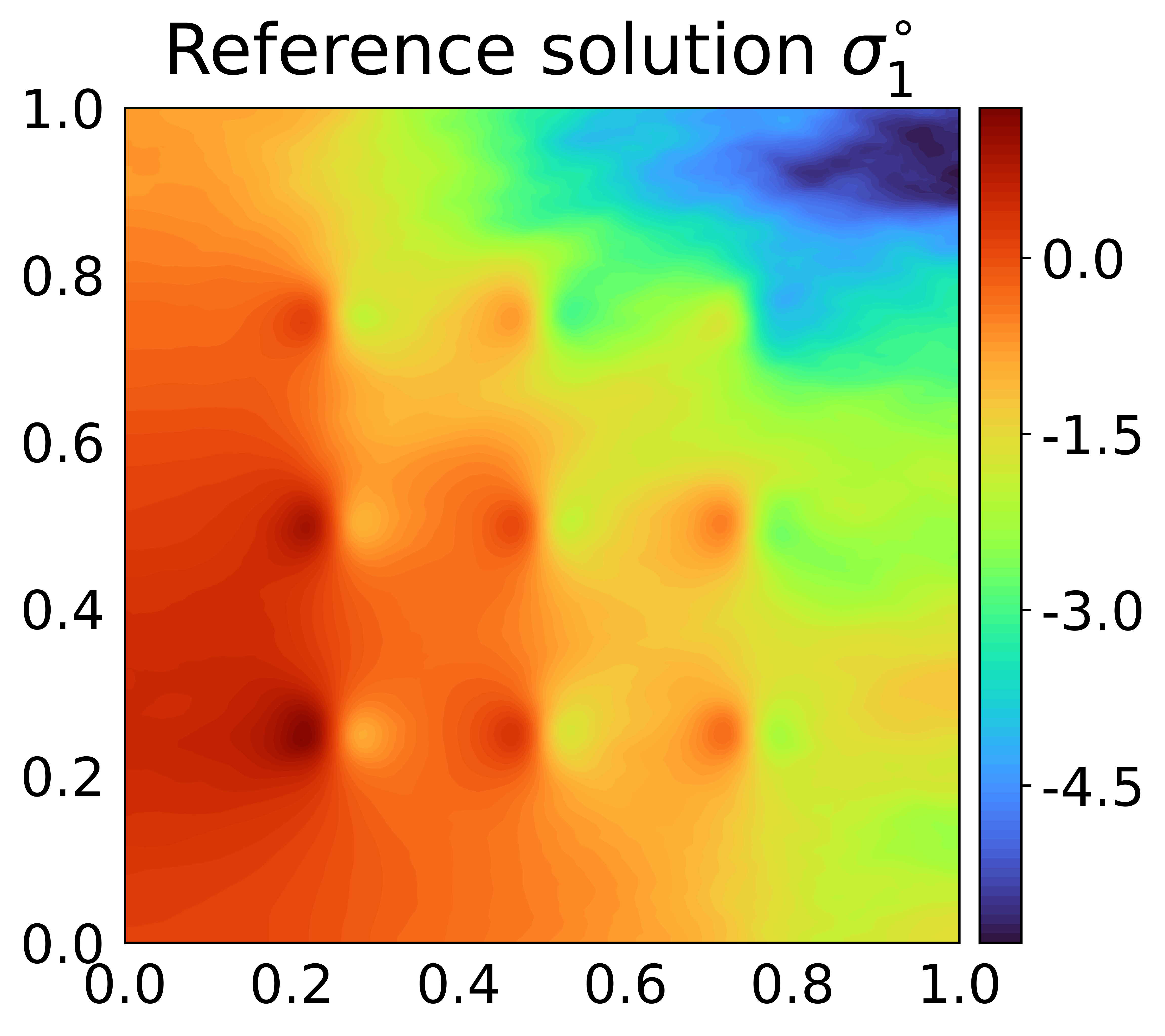}
\includegraphics[width=.24\linewidth]{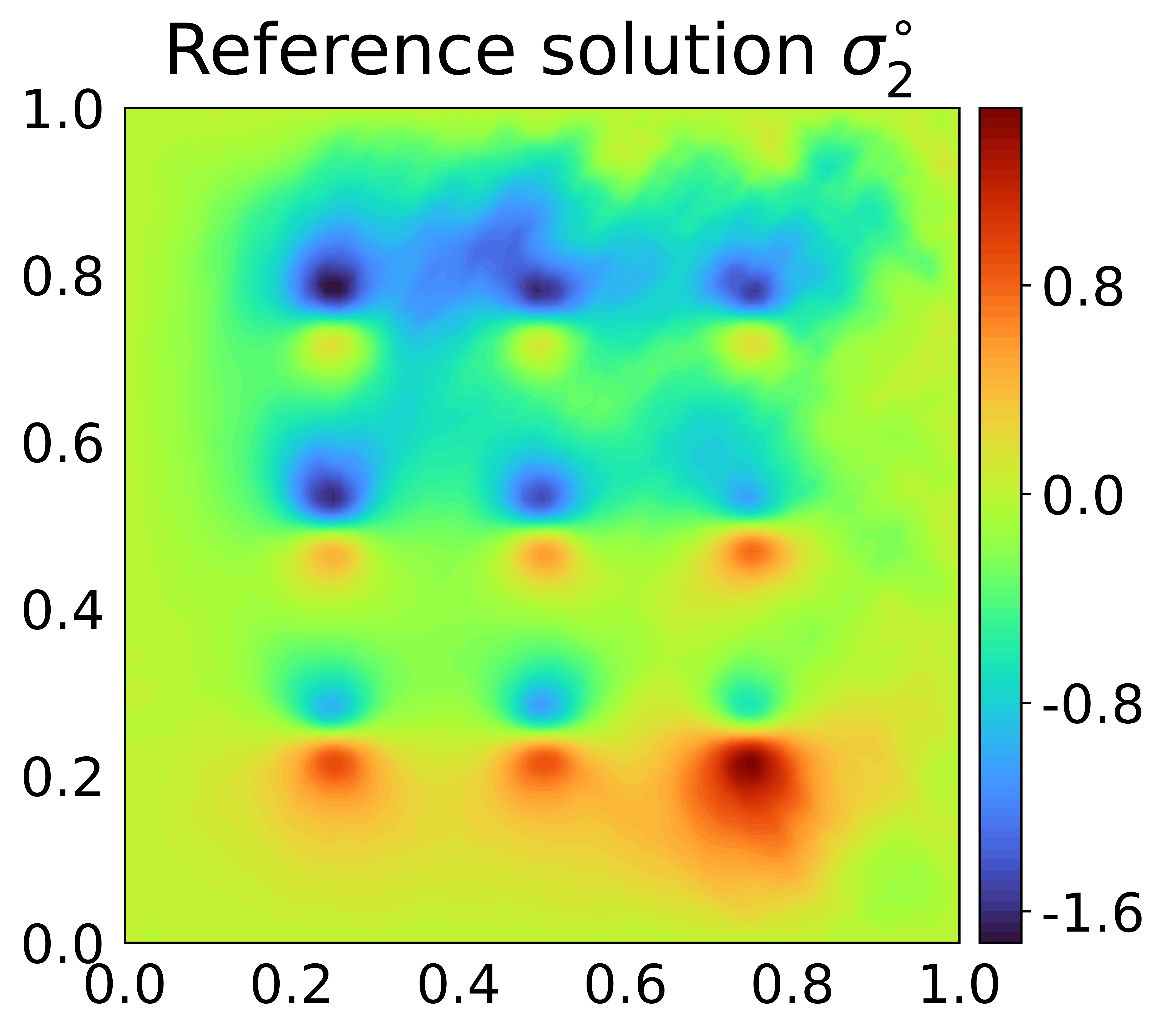}

\vspace{1em}
\includegraphics[width=.32\linewidth]{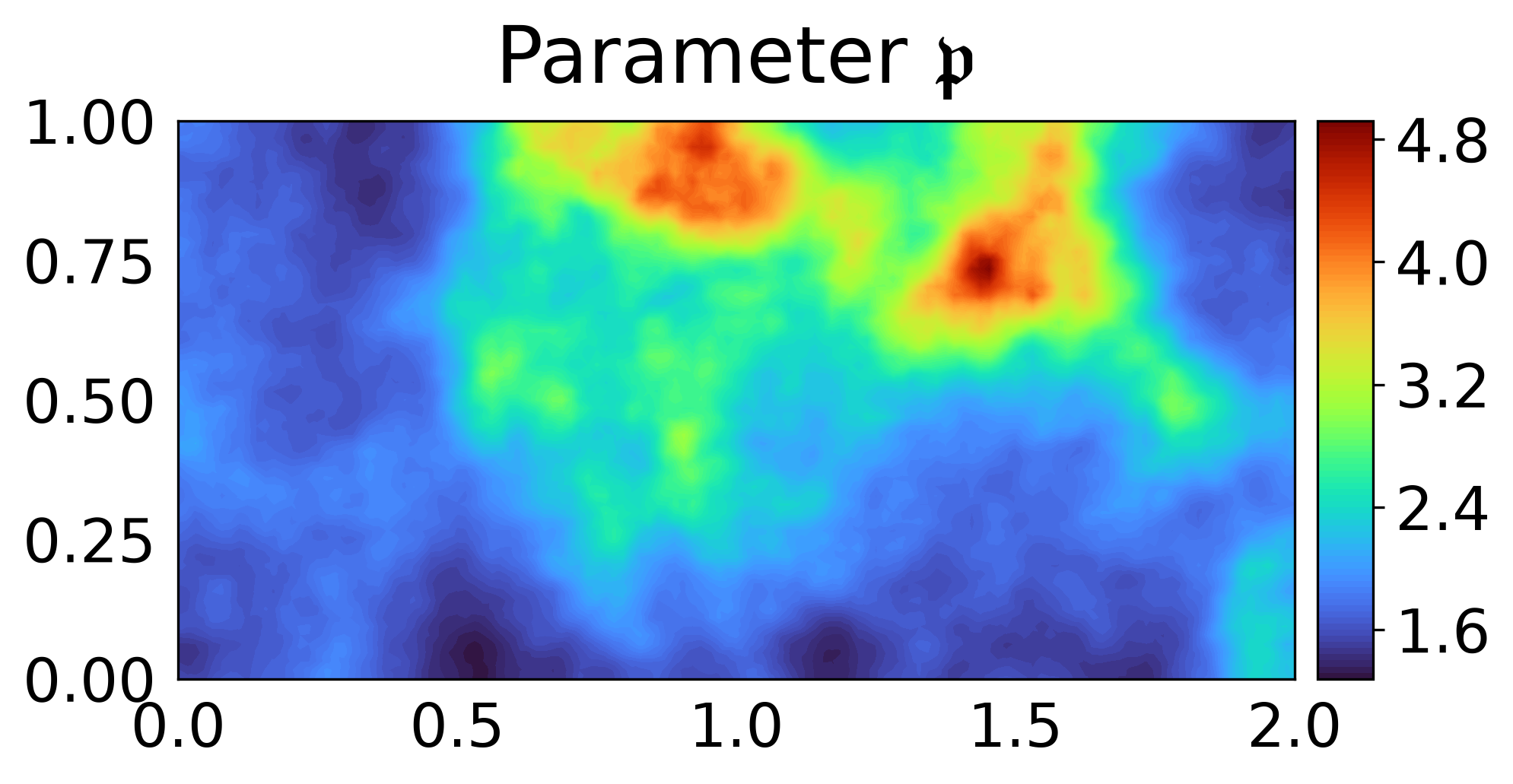}
\includegraphics[width=.32\linewidth]{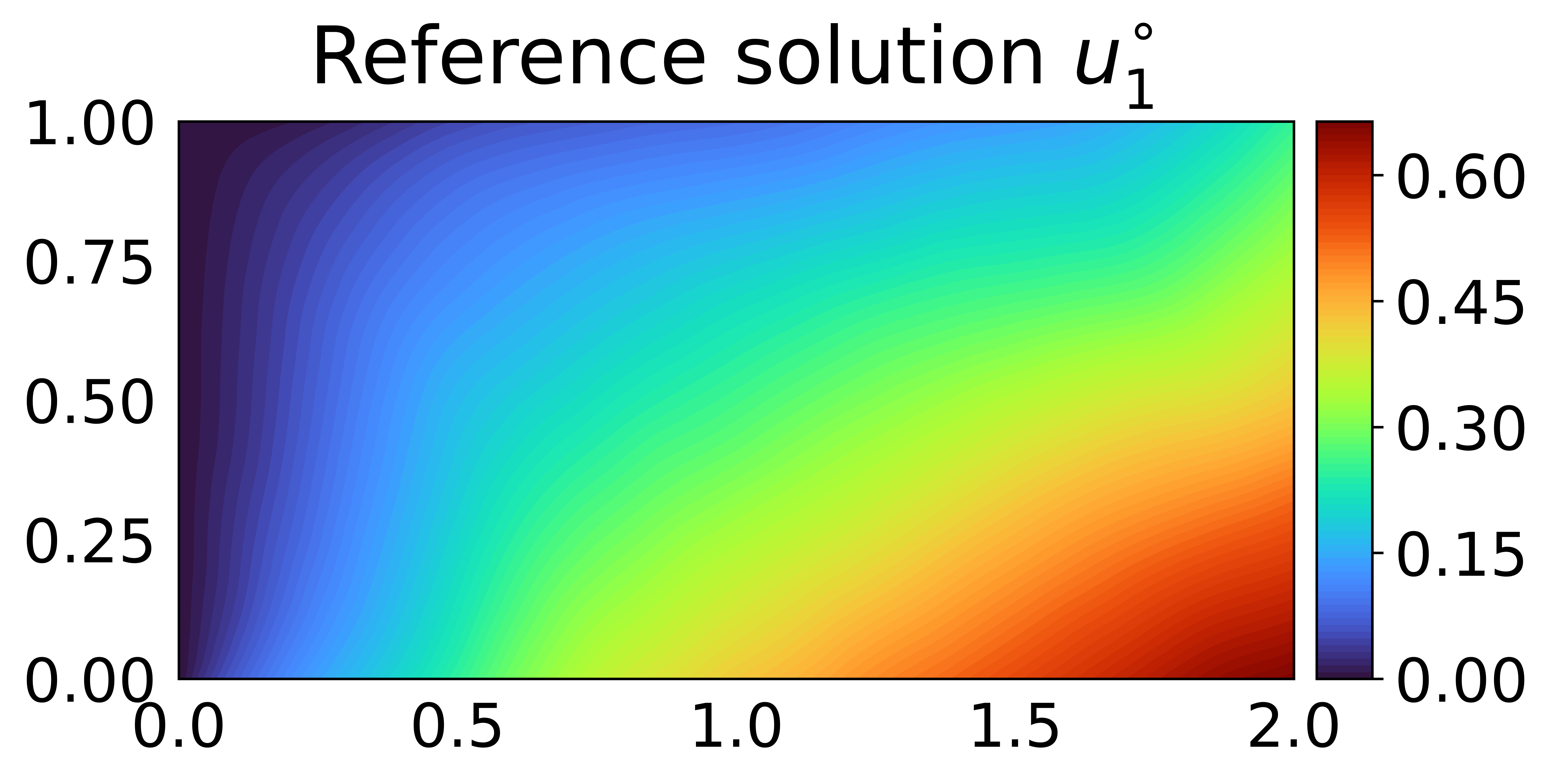}
\includegraphics[width=.32\linewidth]{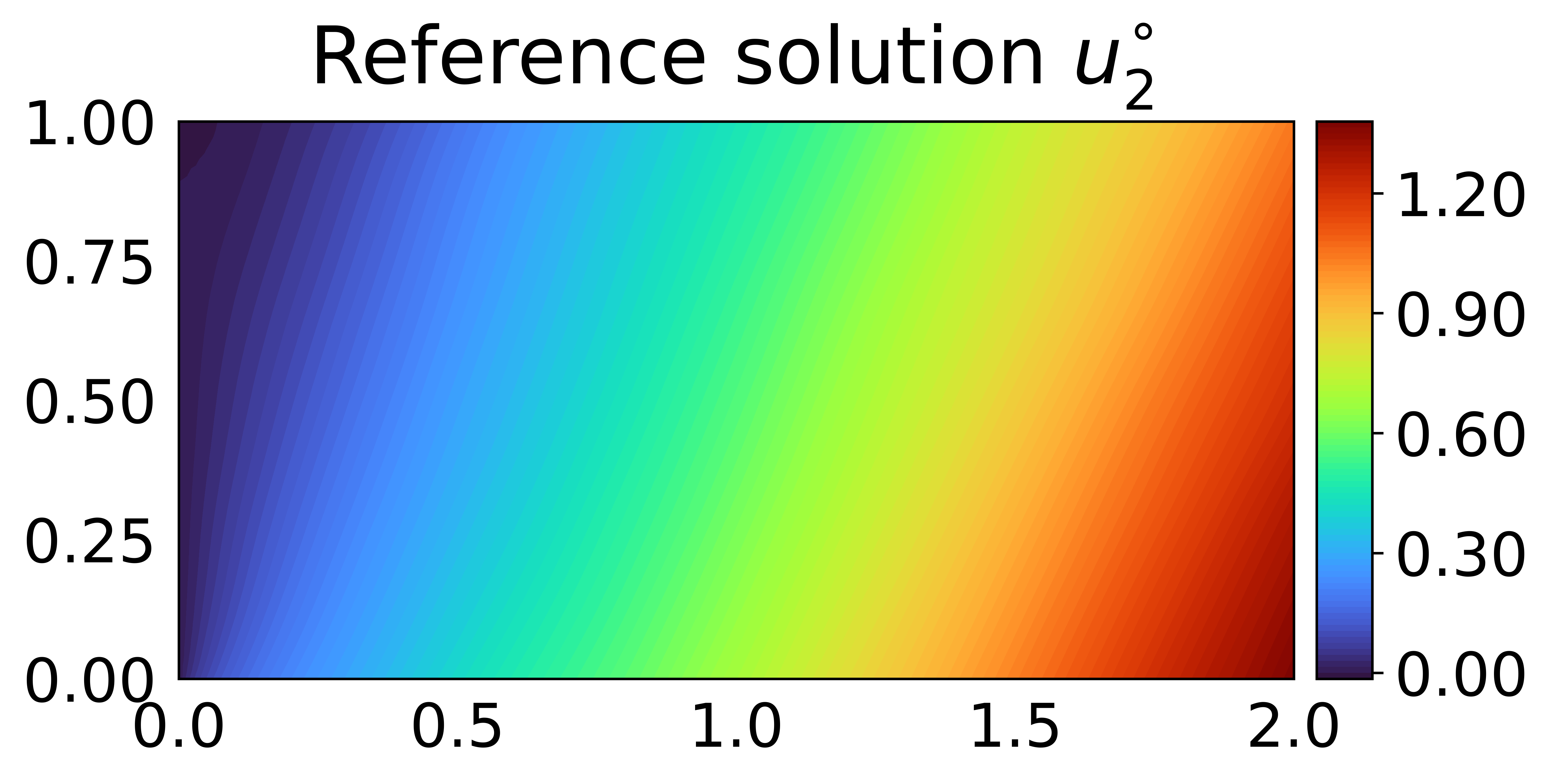}

\vspace{1em}
\includegraphics[width=.24\linewidth]{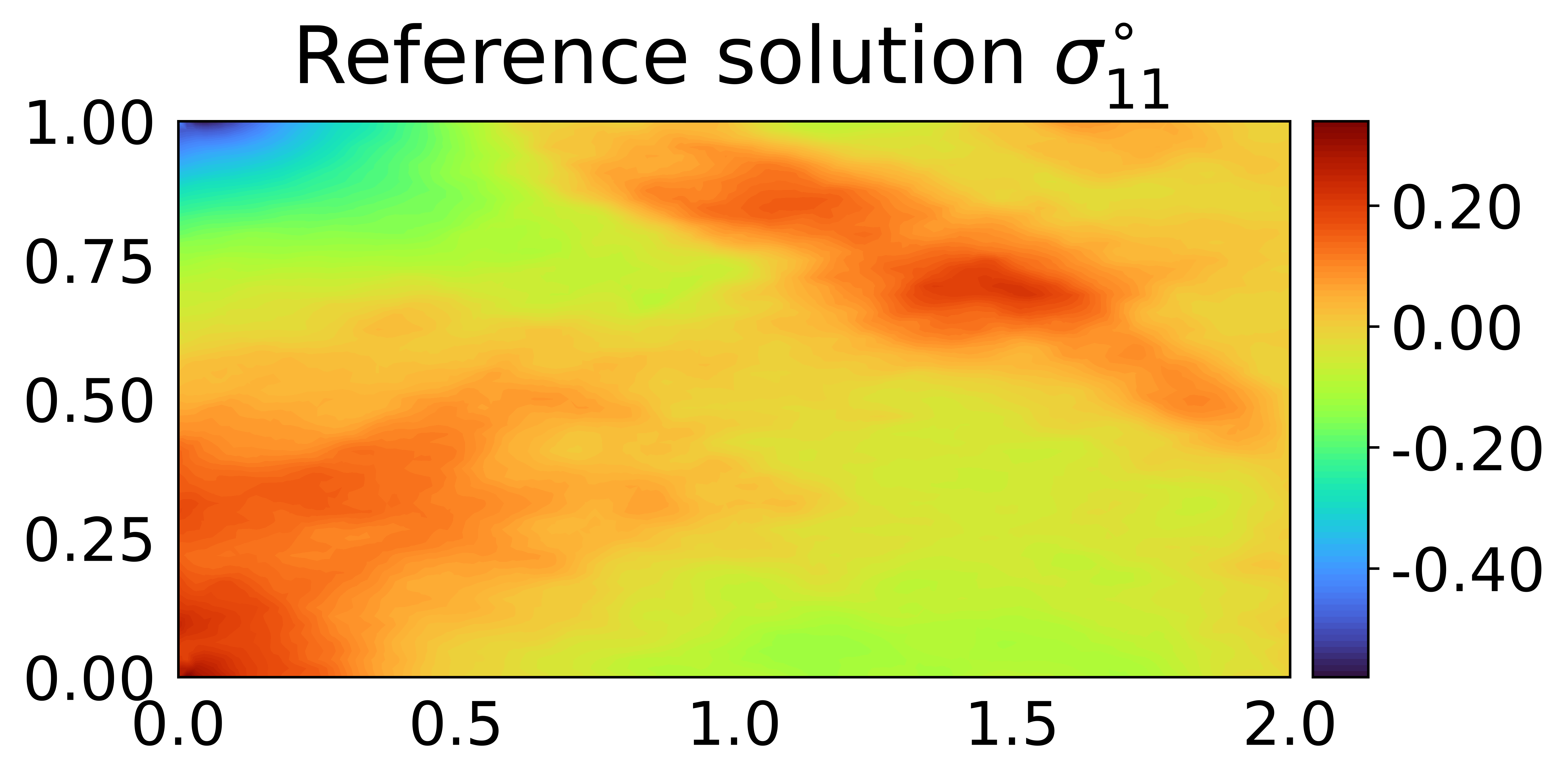}
\includegraphics[width=.24\linewidth]{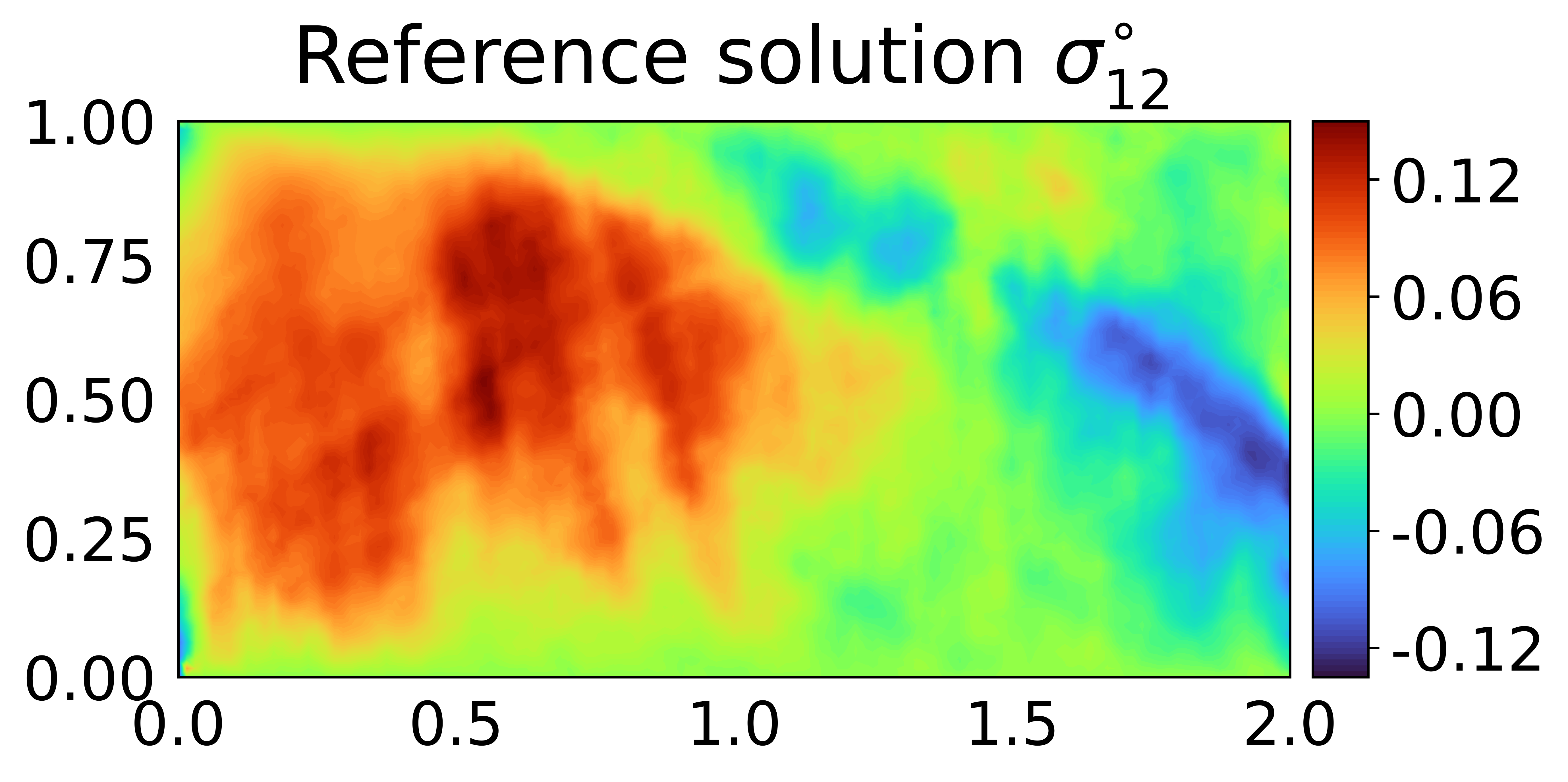}
\includegraphics[width=.24\linewidth]{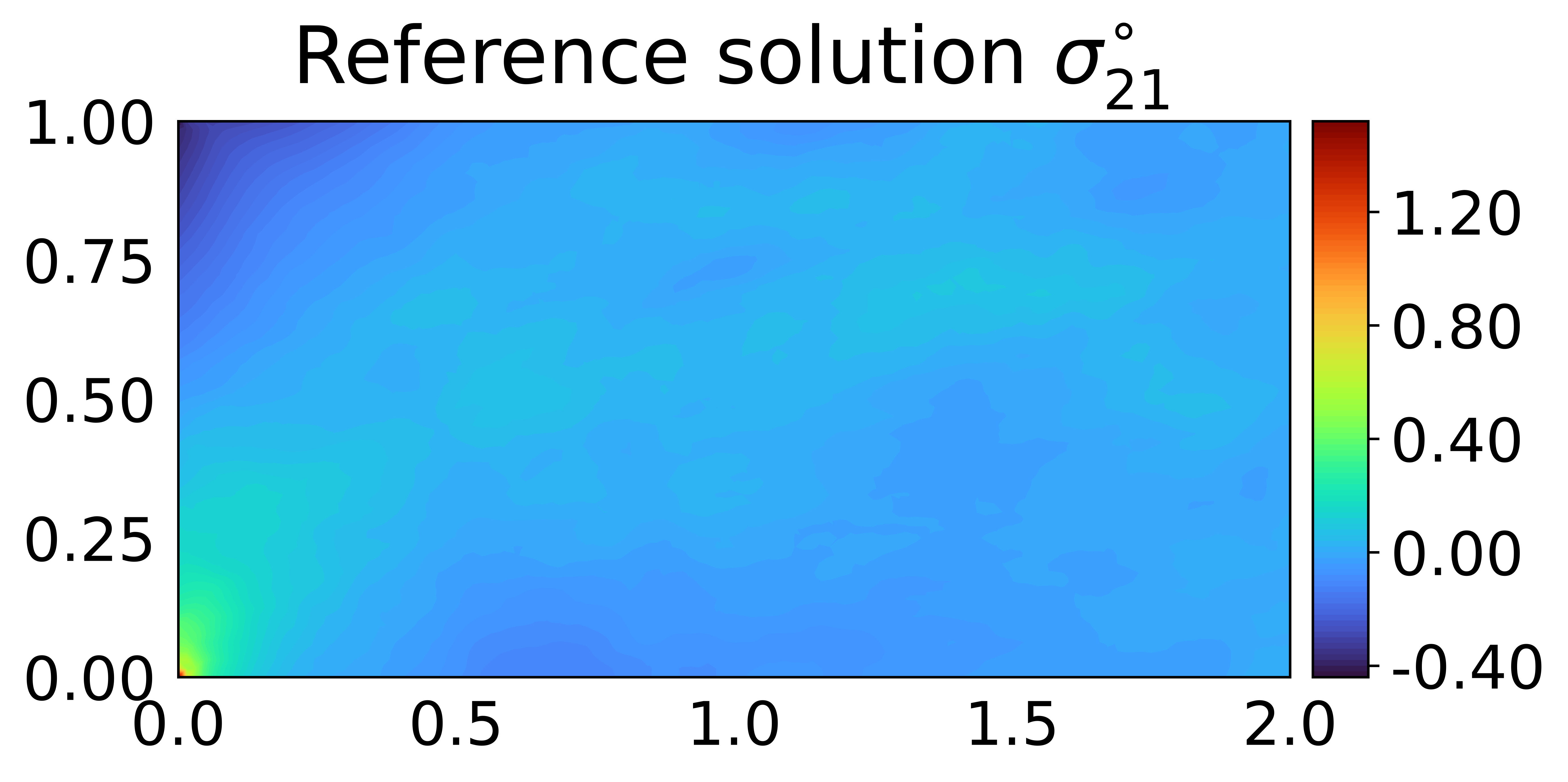}
\includegraphics[width=.24\linewidth]{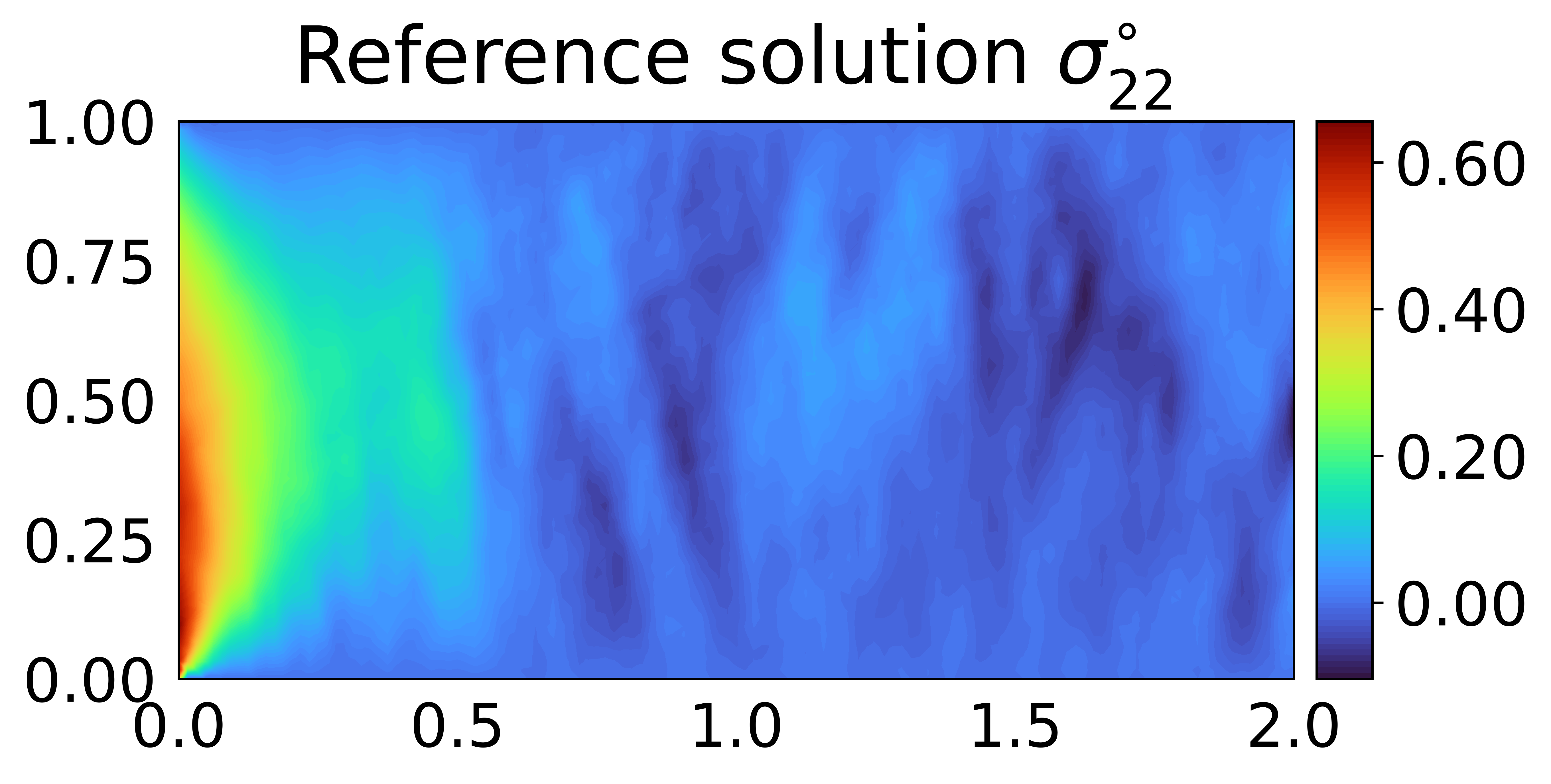}
\caption{One test parameter and the corresponding finite element solution for each benchmark. Top: semilinear diffusion--reaction benchmark ($\kappa=10^3$), with the conductivity $\pp$, the primal variable $u^\circ$, and the two components of the flux $\sigma^\circ$. Second row: quasilinear diffusion benchmark, in the same format. Bottom two rows: compressible hyperelasticity benchmark, with the Young's modulus, the two components of the displacement $u^\circ$, and the four components of the stress $\sigma^\circ$. The plots show the homogeneous parts $u^\circ$ and $\sigma^\circ$ of~\eqref{eq:lifts}.}\label{fig:benchmarks}
\end{figure}

We consider three benchmarks, one for each class of nonlinear PDEs in \cref{sec:examples}. We summarize their configurations in \cref{tab:benchmarks} and show a representative test parameter and its finite element solution in \cref{fig:benchmarks}. The parameters are sixteen inclusion conductivities with a contrast of up to one hundred for semilinear diffusion--reaction, a lognormal permeability field with a positive offset for quasilinear diffusion, and a lognormal stiffness field with a positive offset for compressible hyperelasticity. The corresponding solutions exhibit steep gradients near the inclusions, localized responses to the wells, and corner singularities at the clamped edge, respectively.

\begin{table}[!htb]
\centering
\footnotesize
\renewcommand{\arraystretch}{1.2}
\setlength{\tabcolsep}{4pt}
\newcommand{\grp}[1]{\multicolumn{4}{@{}l}{\itshape #1}\\[1pt]}
\begin{tabularx}{\linewidth}{@{}
  >{\raggedright\arraybackslash}p{0.16\linewidth}
  >{\raggedright\arraybackslash\hsize=1.05\hsize}X
  >{\raggedright\arraybackslash\hsize=0.85\hsize}X
  >{\raggedright\arraybackslash\hsize=1.10\hsize}X
@{}}
\toprule
 & \textbf{Semilinear diffusion--reaction} & \textbf{Quasilinear diffusion} & \textbf{Compressible hyperelasticity}\\
\midrule
\grp{Problem}
Nonlinearity
 & reaction term $\kappa u^3$, $\kappa=10^3$
 & diffusion coefficient $1+u^2$
 & stress law (Neo-Hookean)\\
Equation
 & $-\rdiv(\pp\nabla u)+\kappa u^3=f$
 & $-\rdiv\bigl(\pp(1+u^2)\nabla u\bigr)=f$
 & $\rdiv P(I+\nabla u;\pp)+f=0$\\
Domain
 & $(0,1)^2$ & $(0,1)^2$ & $(0,2)\times(0,1)$\\
Linear counterpart~\cite{QCD}
 & heat conduction & Darcy flow & linear elasticity\\
\addlinespace[3pt]
\grp{Parameter}
Definition
 & $\pp=10^{y_i}$ on $16$ disjoint mini-squares, $\pp=1$ elsewhere
 & $\pp=0.01+e^{m}$
 & $\pp=1+e^{m}$\\
Random input
 & $y_i\sim\mathcal U(-1,1)$, $i=1,\dots,16$
 & $m\sim\mathcal N\bigl(0,(1.5I-0.15\Delta)^{-2}\bigr)$
 & $m\sim\mathcal N\bigl(0,(1.5I-0.15\Delta)^{-2}\bigr)$ \\
Bounds
 & $0.1\le\pp\le10$
 & $0.01\le\pp\le0.01+e^{\hat m}$
 & $1\le\pp\le1+e^{\hat m}$\\
\addlinespace[3pt]
\grp{Boundary conditions and source term}
Dirichlet
 & $u_D=0.1(1-x_1)\sin(4\pi x_2)$\newline left/right
 & $u_D=1-x_1$\newline left/right
 & $u_D=(0,0)^\top$\newline left\\
\addlinespace[2pt]
Neumann
 & $g=0.1(1-x_2)\cos(2\pi x_1)$\newline top/bottom
 & $g=0$\newline top/bottom
 & $g=\begin{pmatrix}0.6\,e^{-(x_2-0.5)^2/4}\\[1pt] 0.3(1+x_2/10)\end{pmatrix}$ right;\newline
   $g=(0,0)^\top$ top/bottom\\
\addlinespace[2pt]
Source term
 & $f_2=1$ whole domain;\newline
   $f_1=(0.5,-0.5)^\top$ mini-squares
 & $f=\sum_{i=1}^{9}100\,e^{-\|x-c_i\|_2^2/\ell^2}$
 & $f=(0,0)^\top$\\
\addlinespace[3pt]
\grp{Discretization and data}
Mesh size
 & $128\times128$ & $128\times128$ & $128\times64$\\
\# DoFs
 & $229\,375$ & $229\,375$ & $229\,376$\\
Quadrature degree $q$
 & $6$ & $8$ & $8$\\
$N_{\mathrm{POD}}$ / $r$
 & $500$ / $128$ & $1000$ / $512$ & $1000$ / $256$\\
Train / valid / test
 & $3000$ / $500$ / $500$
 & $8000$ / $800$ / $500$
 & $5000$ / $500$ / $500$\\
\bottomrule
\end{tabularx}
\caption{The three benchmarks. The linear counterpart is the benchmark
of~\cite{QCD} with the same data and discretization. The mini-squares of
the semilinear benchmark are
$[\tfrac j8-\tfrac1{16},\tfrac j8+\tfrac1{16}]\times
 [\tfrac k8-\tfrac1{16},\tfrac k8+\tfrac1{16}]$, $j,k\in\{1,3,5,7\}$,
and are aligned with the mesh. The Gaussian bumps of the quasilinear
benchmark have centers $c_i=(j/4,k/4)$, $j,k\in\{1,2,3\}$, and width
$\ell=1/32$. The upper bound $\hat m$ is the largest grid-point value of
$m$ over $10^4$ samples: $\hat m=2.85$ (quasilinear diffusion) and $\hat m=2.98$
(compressible hyperelasticity).}
\label{tab:benchmarks}
\end{table}

For all benchmarks, we use a uniform triangulation and the conforming pair $\mathrm{RT}_1\times\mathrm{CG}_2$ for the state $s=[\sigma,u]$. This pair consists of Raviart--Thomas elements of order one, corresponding to degree two in the DOLFINx convention, for the flux and continuous Lagrange elements of degree two for the primal variable. For compressible hyperelasticity, the pair is used row-wise for the stress and component-wise for the displacement. We specify the quadrature degree for each benchmark in \cref{tab:benchmarks}, since the automatic degree estimates in DOLFINx lead to unnecessarily expensive integration of the nonlinear terms. For semilinear diffusion--reaction, all integrands are piecewise polynomial, and quadrature of degree twelve integrates the loss exactly. In our tests, degree six yields negligible relative changes in the loss, its gradient, and the finite element solution, while reducing the cost of loss and gradient evaluation by more than half.

We compute the finite element solutions by applying the Gauss--Newton method described in \cref{subsec:loss_excess} to the stationarity condition~\eqref{eq:discrete_normal}. This amounts to omitting the term $(\cR(s),D^2\cR(s)[\cdot,\cdot])$ from the Newton Jacobian while retaining the full stationarity residual. The stationary points therefore remain unchanged. In our benchmarks, this method is between $1.3$ and $7$ times faster than Newton's method. At each iteration, we solve a linear system using the direct solver MUMPS and select the step size by backtracking line search. We set both relative and absolute stopping tolerances to $10^{-8}$.

For semilinear diffusion--reaction, we initialize the iteration with the solution at $\kappa=0$ for the same parameter. For quasilinear diffusion, we use the solution corresponding to the mean of the parameter samples. For compressible hyperelasticity, we apply load continuation with factors $0.4$, $0.7$, and $1$. We initialize the first load step with the solution for the mean parameter at load factor $0.4$, and each subsequent step with the solution from the preceding step. All solves converge without further subdivision of the load steps. To keep the logarithm in the Neo-Hookean energy well defined during iteration, we replace $\det F$ by $\max(\det F,10^{-3})$ in both the solver and the loss evaluation. This modification leaves the constitutive law unchanged wherever $\det F>10^{-3}$.

At $\kappa=0$, the semilinear diffusion--reaction benchmark reduces to the heat-conduction problem of~\cite{QCD}, whose solutions satisfy $|u|\lesssim0.2$. Thus, the magnitude of $\kappa u^3$ becomes comparable to $f_2=1$ near the largest values of $|u|$ when $\kappa$ is of order $10^2$, and affects a larger portion of the domain as $\kappa$ increases. We choose $\kappa=10^3$, for which the reaction magnitude exceeds the source on about one third of the domain. Relative to the solution at $\kappa=0$, the finite element solution changes by a median of $42\%$ in the $\HH$ norm, with a range of $27$--$61\%$. Starting from the solution at $\kappa=0$, Gauss--Newton iteration converges in four or five steps.

\subsection{Residual and error metrics}\label{subsec:numerical_metrics}
For a parameter sample $\pp_i$, let $s_{h,i}=s_h(\pp_i)$ denote the finite element solution and $s_i=s_r(\pp_i;\hat\theta)$ the prediction of a neural operator with trained weights $\hat\theta$. For each test sample, we record the residual of the prediction and of the finite element solution, the error of the prediction, and the floor-corrected ratio of error and residual,
\begin{equation}\label{eq:numerical_metrics}
 \eta_i=\sqrt{J_{h,q}(s_i;\pp_i)},\qquad
 \eta_{h,i}=\sqrt{J_{h,q}(s_{h,i};\pp_i)},\qquad
 E_{h,i}=\norm{s_i-s_{h,i}}_\HH,\qquad
 \chi_i=\frac{E_{h,i}}{\bigl(\eta_i^2-\eta_{h,i}^2\bigr)^{1/2}} .
\end{equation}
By \cref{lem:loss_excess}, the loss excess $\eta_i^2-\eta_{h,i}^2$ is equivalent to the squared error $E_{h,i}^2$, with the constants of the linearized problem and up to the perturbation terms in~\eqref{eq:loss_excess}. The ratio $\chi_i$ therefore measures how well this equivalence holds in the finite element space. For a linear problem it satisfies $\Clin^{-1/2}\le\chi_i\le\clin^{-1/2}$ with the constants on $\HH_h$. We report the plain ratio $E_{h,i}/\eta_i$ only for comparison: it compares an error that vanishes at $s_{h,i}$ with a residual that does not, and hence tends to zero for accurate predictions. Both ratios coincide when $\eta_{h,i}\ll\eta_i$. In addition, we report the relative error $E_{h,i}/\norm{s_{h,i}}_\HH$ and, as in~\cite{QCD}, the relative error in the norm of $\LL^2=L^2\times L^2$.

\subsection{Accuracy and computational cost of RBNO}\label{subsec:exp_accuracy}
For each benchmark, we construct the POD basis described in \cref{subsec:pod} from the $N_{\mathrm{POD}}$ snapshots specified in \cref{tab:benchmarks}. We report the POD eigenvalues and projection errors on the test set in \cref{fig:pod_record} of \ref{app:numerical_results}. RBNO predicts the $r$ reduced basis coefficients and is trained on the parameter samples listed in \cref{tab:benchmarks}, using the loss and gradient assembly described in \cref{subsec:gradient}. For semilinear diffusion--reaction, we use the heat-conduction network architecture of~\cite{QCD} and train separate models for $\kappa=0$ and $\kappa=10^3$ on the same parameter samples. We provide the architectures and training settings in \ref{app:implementation} and retain the checkpoint with the lowest validation loss for each model.

\begin{table}[!htp]
\centering\footnotesize
\setlength{\tabcolsep}{4pt}
\begin{subtable}{\linewidth}
\centering
\begin{tabular}{@{}lccccc@{}}\toprule
 & & \multicolumn{4}{c}{Relative $\HH$ error $\norm{\,\cdot\,}_\HH/\norm{s_h}_\HH$ of}\\\cmidrule(l){3-6}
Benchmark & $r$ & \makecell{$s-s_h$\\(RBNO)} & \makecell{$s-P_rs_h$\\(coefficient part)} & \makecell{$P_rs_h-s_h$\\(basis part)} & \makecell{$s_r^\dagger-s_h$\\(reduced minimizer)}\\\midrule
Semilinear diffusion--reaction, $\kappa=0$    & $128$ & \num{1.09e-02} & \num{1.09e-02} & \num{5.69e-04} & \num{1.15e-03}\\
Semilinear diffusion--reaction, $\kappa=10^3$ & $128$ & \num{8.43e-03} & \num{8.30e-03} & \num{1.44e-03} & \num{2.57e-03}\\
Quasilinear diffusion                      & $512$ & \num{8.33e-03} & \num{7.14e-03} & \num{4.23e-03} & \num{4.58e-03}\\
Compressible hyperelasticity               & $256$ & \num{2.41e-02} & \num{1.49e-02} & \num{1.88e-02} & \num{2.18e-02}\\\bottomrule
\end{tabular}
\caption{Mean relative errors in the $\HH$ norm.}\label{tab:accuracy_error}
\end{subtable}

\vspace{0.8em}
\begin{subtable}{\linewidth}
\centering
\begin{tabular}{@{}lccccc@{}}\toprule
 & \multicolumn{5}{c}{Mean loss $J_{h,q}$ of}\\\cmidrule(l){2-6}
Benchmark & \makecell{$s$\\(RBNO, test)} & \makecell{$s$\\(RBNO, training)} & \makecell{$P_rs_h$\\(projection)} & \makecell{$s_r^\dagger$\\(reduced minimizer)} & \makecell{$s_h$\\(floor)}\\\midrule
Semilinear diffusion--reaction, $\kappa=0$    & \num{1.95e-04} & \num{1.49e-04} & \num{1.09e-04} & \num{1.08e-04} & \num{1.06e-04}\\
Semilinear diffusion--reaction, $\kappa=10^3$ & \num{1.35e-04} & \num{1.16e-04} & \num{1.01e-04} & \num{9.67e-05} & \num{9.32e-05}\\
Quasilinear diffusion                      & \num{2.12e-02} & \num{8.83e-03} & \num{1.09e-02} & \num{8.41e-03} & \num{4.00e-03}\\
Compressible hyperelasticity               & \num{2.15e-03} & \num{1.41e-03} & \num{2.75e-03} & \num{1.34e-03} & \num{2.19e-05}\\\bottomrule
\end{tabular}
\caption{Mean losses.}\label{tab:accuracy_loss}
\end{subtable}
\caption{Accuracy of RBNO over $500$ test parameters for each benchmark. We denote the RBNO prediction by $s=s_r(\pp;\hat\theta)$, the $\HH$-orthogonal projection of the finite element solution onto $\HH_r$ by $P_rs_h$, and the reduced minimizer of \cref{prop:rbno_error} by $s_r^\dagger$. In (a), we report mean relative errors. For each prediction, the squared error is the sum of the squared coefficient and basis contributions by~\eqref{eq:state_pythagoras}. In (b), we report mean test losses, together with the final mean RBNO training loss. The mean loss of the finite element solutions gives the finite element loss floor.}\label{tab:accuracy}
\end{table}

We summarize the test accuracy in \cref{tab:accuracy} and compare RBNO predictions with the finite element solutions for the parameters shown in \cref{fig:benchmarks} in \cref{fig:ref_pred_diff_semilinear,fig:ref_pred_diff_quasilinear,fig:ref_pred_diff_hyperelastic}. By~\eqref{eq:state_pythagoras}, the squared prediction error splits into the basis contribution $\norm{P_rs_h-s_h}_\HH^2$, which cannot be reduced by training with a fixed basis, and the coefficient contribution $\norm{s-P_rs_h}_\HH^2$. The coefficient contribution dominates for semilinear diffusion--reaction and quasilinear diffusion, whereas the basis contribution is larger for compressible hyperelasticity.

To assess the loss attainable in the reduced space, we approximate $s_r^\dagger$ for each test parameter by applying Gauss--Newton iteration to $J_{h,q}$ in the $r$ reduced basis coefficients, initialized at $P_rs_h$. The iteration converges for every test parameter in at most seven steps. The projection $P_rs_h$ has a smaller error but a larger loss than the computed reduced minimizer, consistent with its definition as the best approximation in the $\HH$ norm.

For quasilinear diffusion and compressible hyperelasticity, the final mean RBNO training loss is about $5\%$ above the mean loss of the reduced minimizers on the test set. This comparison suggests that training reaches losses close to those attainable in the reduced space, although the averages are taken over different parameter sets. The larger RBNO test losses indicate a generalization gap; we show the training and validation histories in \cref{fig:training_record}. For semilinear diffusion--reaction, the corresponding differences are about $20\%$ for $\kappa=10^3$ and $38\%$ for $\kappa=0$. Here the finite element floor accounts for much of the total loss, so the differences are more pronounced after subtracting this floor: the mean RBNO error is approximately three times that of the reduced minimizer for $\kappa=10^3$ and ten times that for $\kappa=0$.

\begin{table}[!htp]
\centering\small
\begin{tabular}{@{}lcccc@{}}\toprule
Benchmark & Training (h) & Loss and gradient (s) & Prediction (s) & FE solve (s)\\\midrule
Semilinear diffusion--reaction, $\kappa=10^3$ & $0.5$ & $0.080$ & $0.0002$ & $7.9$\\
Quasilinear diffusion & $4.1$ & $0.145$ & $0.0012$ & $4.9$\\
Compressible hyperelasticity & $1.2$ & $0.079$ & $0.0006$ & $37.9$\\\bottomrule
\end{tabular}
\caption{Computational cost of RBNO. Training: total wall-clock time using one GPU and $60$ CPU workers for loss and gradient assembly. Loss and gradient: mean time per evaluation in one CPU process over $20$ parameters. Prediction: time per parameter for batched GPU inference, including the reconstruction $\Pi_r\bss_r$. FE solve: mean time per solve in one CPU process over $5$ parameters, using the initialization described in \cref{subsec:problem_configs}. For semilinear diffusion--reaction, the solve time includes the linear solve at $\kappa=0$, which takes $2.2$~s. For compressible hyperelasticity, it includes all continuation load steps.}\label{tab:cost}
\end{table}

We report the computational costs in \cref{tab:cost}. Evaluating the loss and its gradient takes approximately $1/30$ to $1/500$ of the time required for a finite element solve, and residual-based training requires no additional finite element solves after constructing the reduced basis. Once trained, RBNO achieves a per-parameter prediction time of $0.2$--$1.2$ milliseconds in batched GPU inference, including reconstruction of the finite element state.

\subsection{Residual--error relation}\label{subsec:exp_ratio}
\begin{figure}[!htp]
\centering
\includegraphics[width=.32\linewidth]{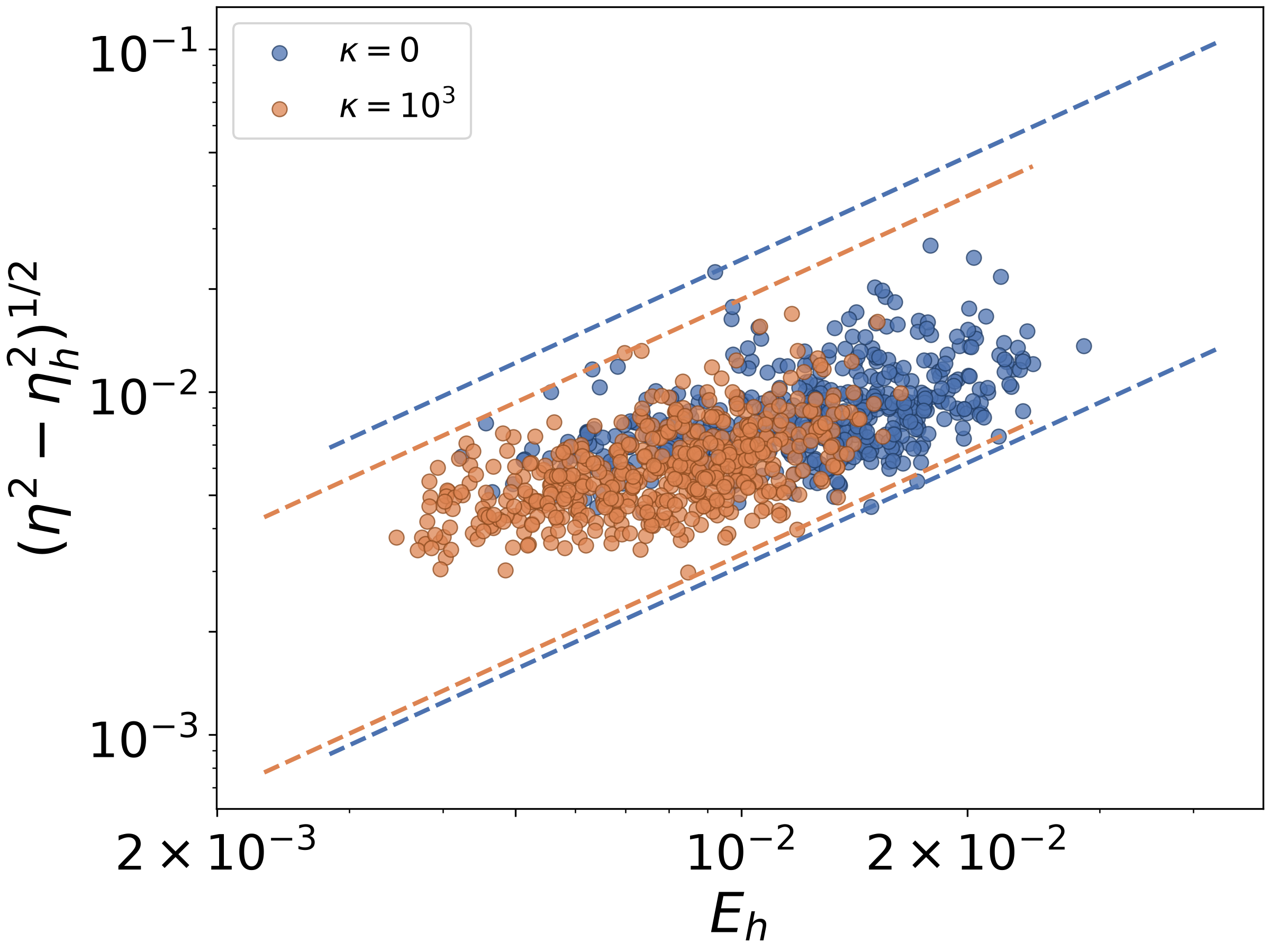}\hfill
\includegraphics[width=.32\linewidth]{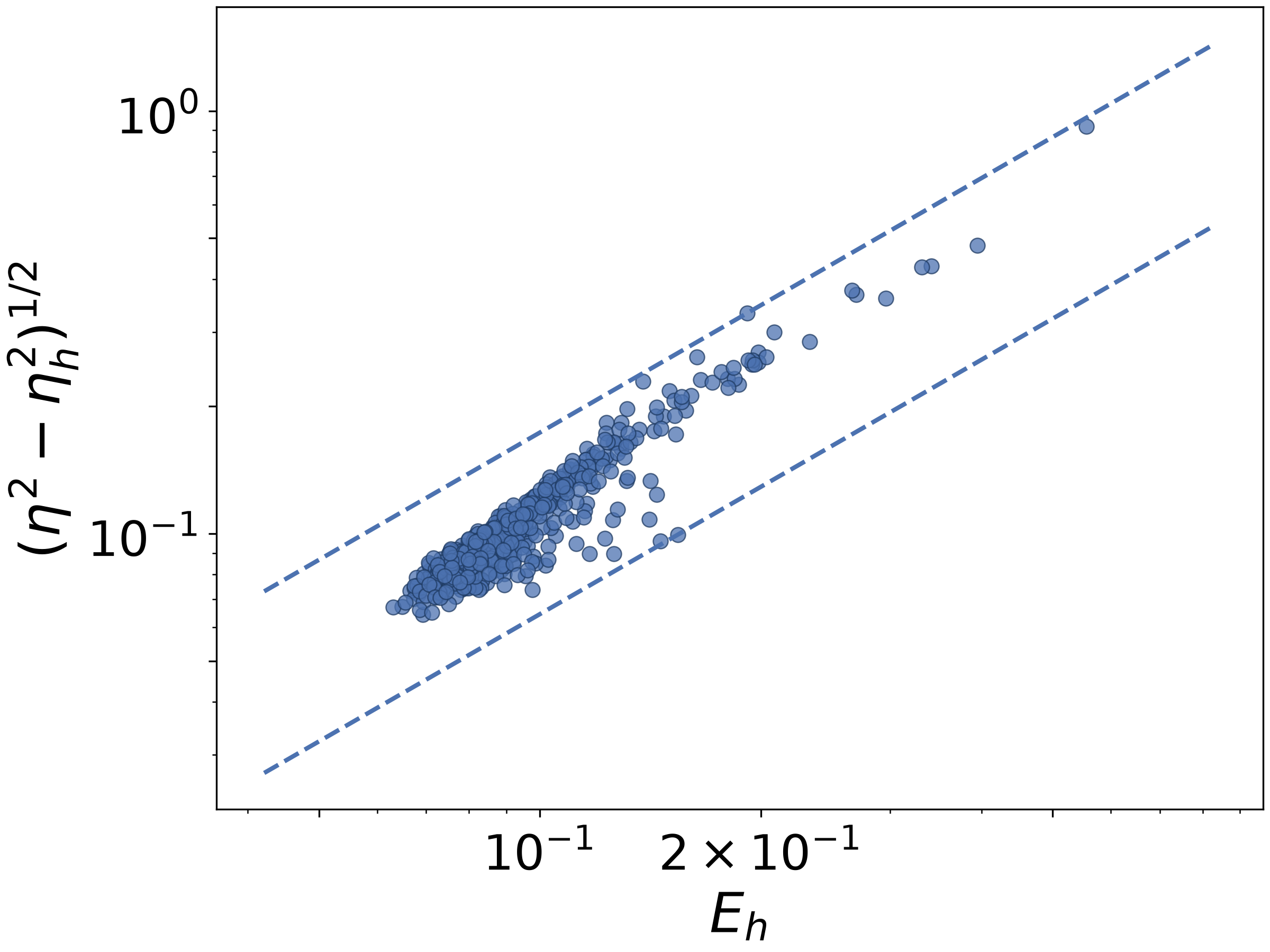}\hfill
\includegraphics[width=.32\linewidth]{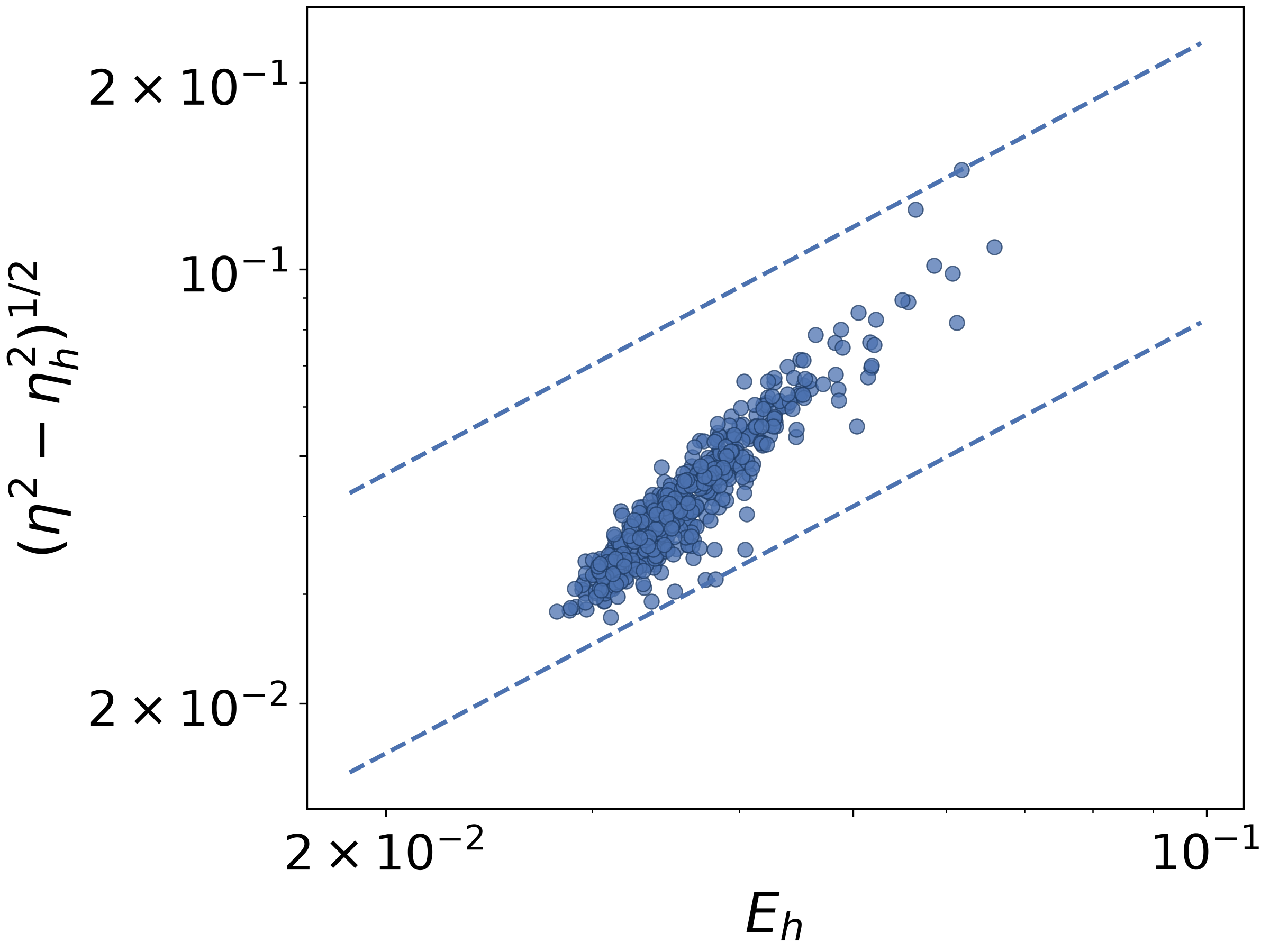}\\[4pt]
\includegraphics[width=.32\linewidth]{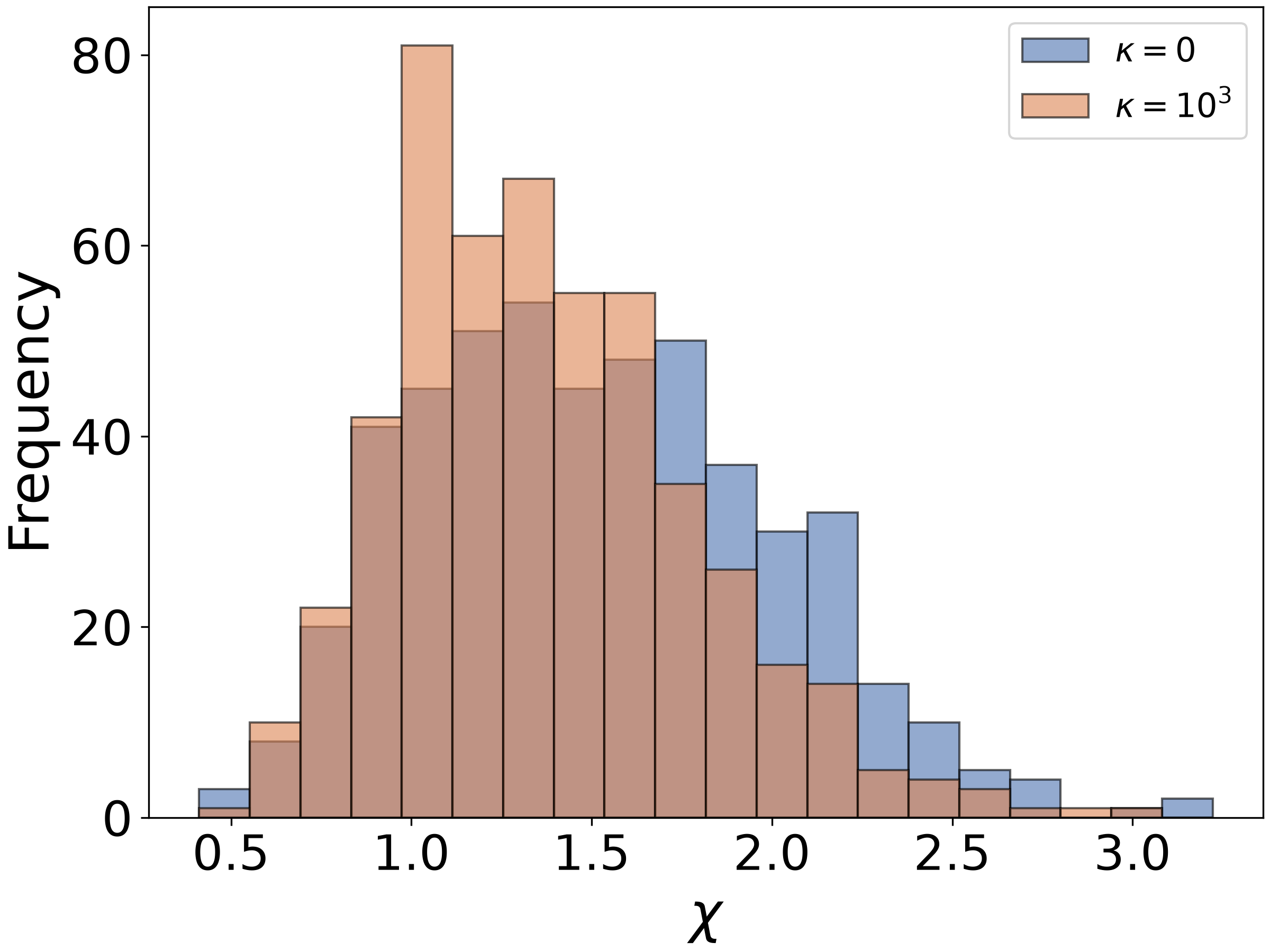}\hfill
\includegraphics[width=.32\linewidth]{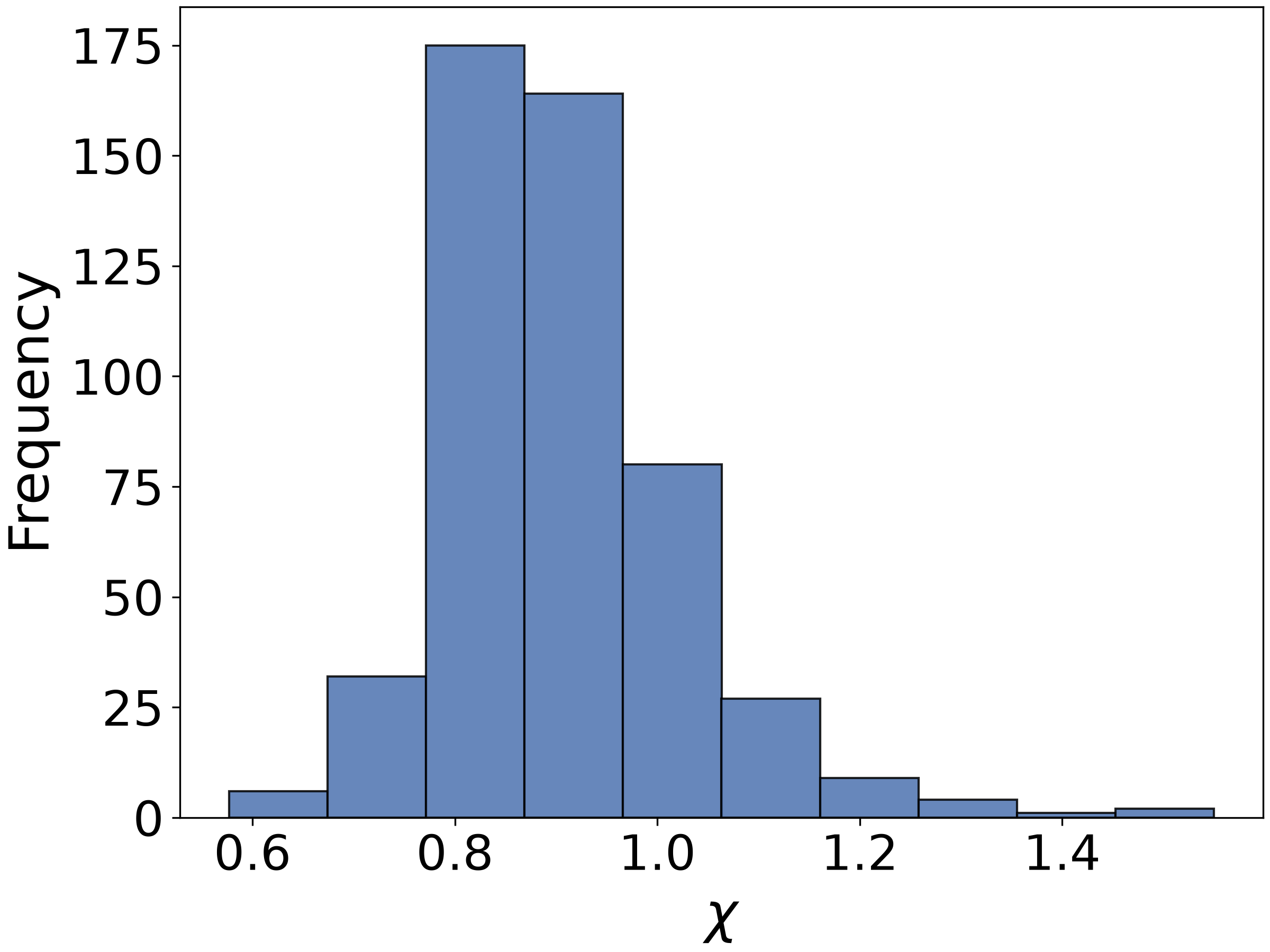}\hfill
\includegraphics[width=.32\linewidth]{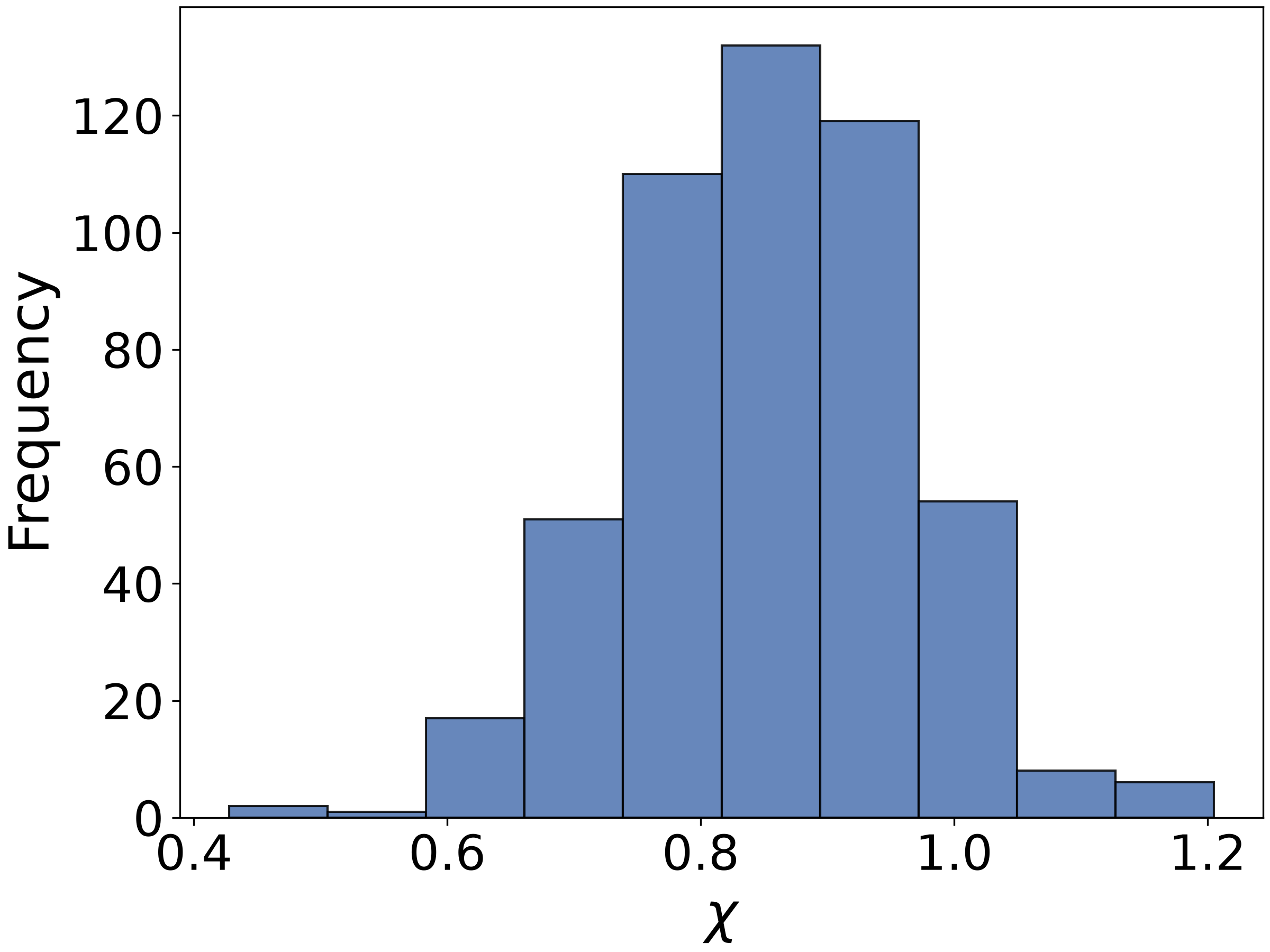}
\caption{Residual--error relation over $500$ test parameters for each benchmark. From left to right: semilinear diffusion--reaction with $\kappa=0$ (blue) and $\kappa=10^3$ (orange), quasilinear diffusion, and compressible hyperelasticity. Top row: we plot the square root of the loss excess, $(\eta_i^2-\eta_{h,i}^2)^{1/2}$, against the error $E_{h,i}$ on log--log axes. The dashed lines correspond to the minimum and maximum observed values of $\chi_i$ and indicate the empirical spread rather than the theoretical constants in \cref{thm:local_vc}. Bottom row: histograms of $\chi_i$.}
\label{fig:ratios}
\end{figure}

We examine the relation between the residual and the error on the test sets in \cref{tab:ratios,fig:ratios}. For every test sample, the prediction has a larger loss than the finite element solution, so $\eta_i>\eta_{h,i}$ and the ratio $\chi_i$ is well defined.

\begin{table}[!htp]
\centering\small
\begin{tabular}{@{}lccccc@{}}\toprule
 & \multicolumn{2}{c}{$E_h/\eta$} & \multicolumn{2}{c}{$\chi$} & \\\cmidrule(lr){2-3}\cmidrule(lr){4-5}
Benchmark & median & range & median & range & Floor share $\overline{\eta_h^2}/\overline{\eta^2}$\\\midrule
Semilinear diffusion--reaction, $\kappa=0$ & $0.94$ & $0.38$--$1.66$ & $1.48$ & $0.41$--$3.22$ & $54\%$\\
Semilinear diffusion--reaction, $\kappa=10^3$ & $0.68$ & $0.35$--$1.22$ & $1.31$ & $0.54$--$2.98$ & $69\%$\\
Quasilinear diffusion & $0.74$ & $0.57$--$1.31$ & $0.89$ & $0.58$--$1.55$ & $19\%$\\
Compressible hyperelasticity & $0.85$ & $0.43$--$1.20$ & $0.86$ & $0.43$--$1.21$ & $1\%$\\\bottomrule
\end{tabular}
\caption{Residual--error ratios over $500$ test parameters for each benchmark, using the notation of~\eqref{eq:numerical_metrics}. We report the median and range of the plain ratio $E_h/\eta$ and the floor-corrected ratio $\chi$. The floor share $\overline{\eta_h^2}/\overline{\eta^2}$ is the mean loss of the finite element solutions divided by the mean loss of the RBNO predictions, as reported in \cref{tab:accuracy_loss}.}\label{tab:ratios}
\end{table}

The floor-corrected ratio $\chi$ remains within a moderate range for all benchmarks, with medians between $0.86$ and $1.48$. The ratio of its maximum to its minimum ranges from $2.7$ for quasilinear diffusion to $7.8$ for semilinear diffusion--reaction with $\kappa=0$. These results show that the square root of the loss excess tracks the error relative to the finite element solution, with the observed proportionality factors given by the ranges of $\chi$ in \cref{tab:ratios}.

The identity $E_h/\eta=\chi\,(1-\eta_h^2/\eta^2)^{1/2}$
shows that the plain ratio cannot exceed $\chi$ and that the difference depends on the contribution of the finite element loss floor. For compressible hyperelasticity, the floor accounts for about $1\%$ of the mean prediction loss, and the two ratios nearly coincide. For quasilinear diffusion, this share is about one fifth. The floor arises almost entirely from the equilibrium component because the wells of width $1/32$ are only partly resolved on the $128\times128$ mesh; the same effect was reported in~\cite{QCD} for the linear Darcy flow problem. For semilinear diffusion--reaction, the floor accounts for more than half of the mean prediction loss. Subtracting this contribution makes $\chi$ more sensitive to errors in the evaluated losses, particularly for predictions close to the finite element solution.

In these experiments, adding the cubic reaction term does not weaken the residual--error relation observed for the linear problem. For $\kappa=0$, the finite element stability constants give $\Clin^{-1/2}\le\chi_i\le\clin^{-1/2}$, and the observed values range from $0.41$ to $3.22$. For $\kappa=10^3$, the range is slightly narrower, from $0.54$ to $2.98$. For comparison, \cite{QCD} reports plain ratios between approximately $0.47$ and $1.88$ for the linear problem, although the error there is measured against a higher-order finite element reference solution. The variation in the observed ratios across benchmarks also shows that a fixed residual threshold need not correspond to the same error level for different problems.

\subsection{Stability test}\label{subsec:exp_stability}
For semilinear diffusion--reaction and quasilinear diffusion, we establish lower stability bounds in \cref{prop:semilinear,prop:quasilinear}. For compressible hyperelasticity, the corresponding bound relies on the stability assumption on the branch tangent in \cref{prop:hyperelastic}. We examine \cref{ass:local_stability} numerically by computing the smallest generalized eigenvalue $\lambda_{\min,h}(\bar u,\pp)$ defined in~\eqref{eq:discrete_eigen_check} for the first $50$ test parameters at three groups of states. The first consists of finite element solutions on the mesh specified in \cref{tab:benchmarks} and on two coarser meshes with half and one quarter as many cells in each direction. For semilinear diffusion--reaction and quasilinear diffusion, we use $32\times32$, $64\times64$, and $128\times128$ cells on the unit square. For compressible hyperelasticity, we use $32\times16$, $64\times32$, and $128\times64$ cells on $(0,2)\times(0,1)$. Within each benchmark, refinement preserves the cell shape and halves the mesh size. The second and third groups consist of Gauss--Newton iterates and RBNO predictions on the finest mesh, respectively. For compressible hyperelasticity, we include only the iterates from the final load step.

On each mesh, we linearize at the finite element solution computed on that mesh. This comparison therefore differs from the nested approximation of a fixed operator considered in \cref{sec:stability_check}: changes in $\lambda_{\min,h}$ reflect both mesh refinement and changes in the linearization point. We assemble $G_h$ without projecting the residual, using the same quadrature rule as for the loss, and take $X_h$ to be the Gram matrix of the $\HH$ inner product. After restricting both matrices to the unconstrained degrees of freedom, we approximate $\lambda_{\min,h}$ by the reciprocal of the largest eigenvalue $\mu$ of $X_h\bst=\mu G_h\bst$, computed using a randomized double-pass method~\cite{villa2021hippylib}.

\begin{figure}[!htp]
\centering
\includegraphics[width=\linewidth]{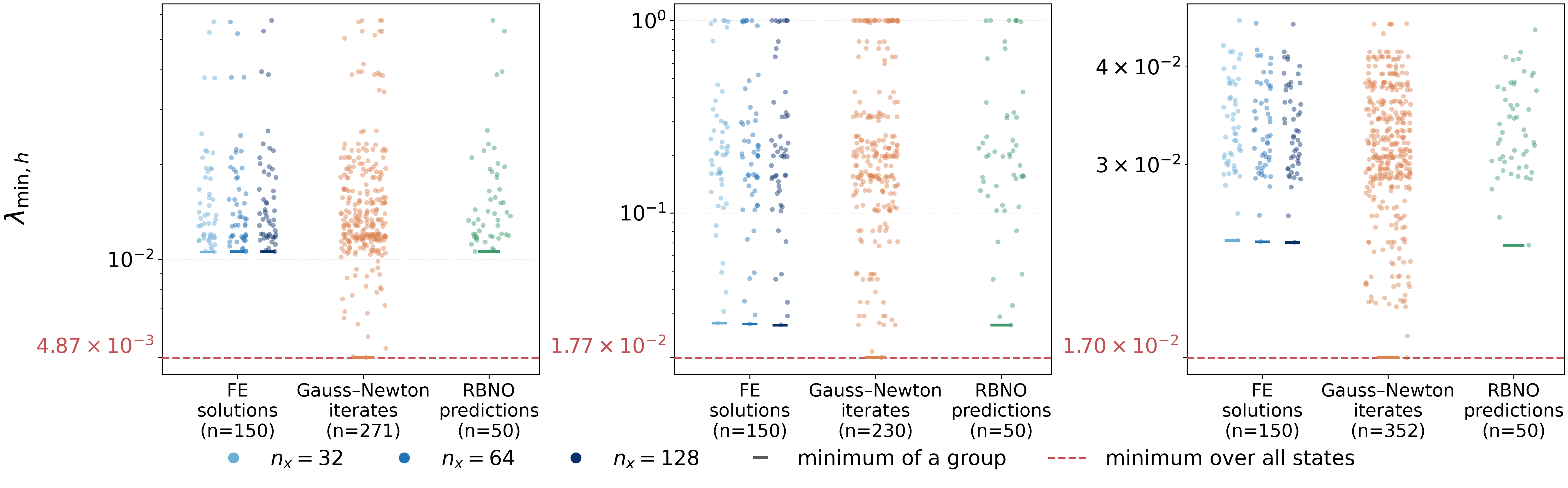}
\caption{Discrete stability for the first $50$ test parameters. We plot the smallest generalized eigenvalue $\lambda_{\min,h}$ defined in~\eqref{eq:discrete_eigen_check} at finite element solutions on three meshes and at Gauss--Newton iterates and RBNO predictions on the finest mesh. The meshes have $n_x=32$, $64$, and $128$ cells in the $x$-direction, with $n_y=n_x$ for semilinear diffusion--reaction and quasilinear diffusion, and $n_y=n_x/2$ for compressible hyperelasticity. For hyperelasticity, we include only the iterates from the final load step. Left: semilinear diffusion--reaction with $\kappa=10^3$. Middle: quasilinear diffusion. Right: compressible hyperelasticity. We denote the number of states in each group by $n$. Horizontal bars indicate the group minima, and the dashed line indicates the minimum over all states, with its value marked on the vertical axis.}\label{fig:stability}
\end{figure}

\begin{table}[!htp]
\centering\small
\begin{tabular}{@{}lcccc@{}}\toprule
Benchmark & States & \makecell{$\min\lambda_{\min,h}$\\(all states)} & \makecell{$\max\lambda_{\min,h}^{-1/2}$\\(at $s_h$)} & \makecell{$\max\chi_i\lambda_{\min,h}^{1/2}$\\(at $s_h$)}\\\midrule
Semilinear diffusion--reaction, $\kappa=10^3$ & $471$ & $4.87\times10^{-3}$ & $9.7$ & $0.34$\\
Quasilinear diffusion & $430$ & $1.77\times10^{-2}$ & $6.2$ & $0.90$\\
Compressible hyperelasticity & $552$ & $1.70\times10^{-2}$ & $6.5$ & $0.18$\\\bottomrule
\end{tabular}
\caption{Summary of the discrete stability tests for the first $50$ test parameters. We report the total number of linearization states and the minimum of $\lambda_{\min,h}$ over these states. For the last two columns, we evaluate $\lambda_{\min,h}$ at the finite element solutions $s_h$ on the finest mesh and report the maximum of $\lambda_{\min,h}^{-1/2}$ and of $\chi_i\lambda_{\min,h}^{1/2}$, respectively. The former gives the leading-order bound for $\chi_i$ from \cref{lem:loss_excess}; values below one in the last column are consistent with this bound.}\label{tab:stability}
\end{table}

We present the eigenvalues at the sampled states in \cref{fig:stability} and summarize the results in \cref{tab:stability}. For each benchmark, the smallest value occurs at an initial Gauss--Newton state. For semilinear diffusion--reaction, this state is the finite element solution for $\kappa=0$, at which we linearize the residual with $\kappa=10^3$. The minima over the finite element solutions are nearly unchanged across the three meshes, and the eigenvalues at RBNO predictions lie in the same range as those at the finite element solutions. Thus, we observe no deterioration of discrete stability under mesh refinement, and the eigenvalues remain positive throughout the sampled solver iterates and predictions.

For semilinear diffusion--reaction, we also compare the computed eigenvalues with the analytic lower bound. At a state with physical primal component $u$, the linearized reaction coefficient is $3\kappa u^2\ge0$. Applying the proof of \cref{prop:semilinear} to the linearized operator gives the lower constant $c_B$ in~\eqref{eq:dr_constants}, with $\zeta=K\norm{3\kappa u^2}_{L^\infty}^{1/2}$. Using the observed grid point evaluations $|u|\le0.13$ at the finite element solutions gives an estimate of $c_B^{-1/2}\approx29$, compared with the computed values $\lambda_{\min,h}^{-1/2}\le9.7$. At the initial states, the observed values $|u|\le0.22$ give $c_B^{-1/2}\approx44$, compared with $\lambda_{\min,h}^{-1/2}\le14.3$. This comparison suggests that the analytic estimate is conservative, as it accounts for the reaction coefficient only through its maximum.

Finally, we relate the eigenvalues at the finite element solutions to the residual--error ratios in \cref{subsec:exp_ratio}. Both the RBNO prediction and the finite element solution belong to $\HH_h$, and $\lambda_{\min,h}(u_h,\pp)$ gives the squared lower stability constant of $\cA_{s_h,\pp}$ restricted to this space. For sufficiently small errors, \cref{lem:loss_excess} therefore yields $\chi_i\lesssim\lambda_{\min,h}(u_{h,i},\pp_i)^{-1/2}$
up to the perturbation terms in~\eqref{eq:loss_excess}. The largest ratios over the full test sets reported in \cref{tab:ratios} are below the largest eigenvalue-based estimates computed for the first $50$ parameters. More directly, for each of these $50$ parameters, we find that $\chi_i\lambda_{\min,h}(u_{h,i},\pp_i)^{1/2}<1$. The largest product is $0.90$ for quasilinear diffusion, where the leading-order estimate is nearly attained, as shown in \cref{tab:stability}. These results provide numerical support for \cref{ass:local_stability} at the sampled states, including those of compressible hyperelasticity, for which we do not have an analytic lower bound on the tangent operator.

\subsection{Comparison with other neural operators}\label{subsec:exp_baselines}
We compare RBNO with three grid-based neural operators, PCA-Net~\cite{bhattacharya2021model}, FNO~\cite{li2020fourier}, and FC-PINO~\cite{ganeshram2022fc}, as well as two data-driven RBNO variants. These variants use the same network, reduced basis, and optimizer as RBNO but replace the residual loss with the relative squared error to the finite element solution, measured in either the $\LL^2$ norm or the $\HH$ norm.

The grid-based operators map a parameter field to the full first-order state, including both the primal variable and the flux or stress. We represent their inputs and outputs as images on a uniform grid whose points coincide with the vertices of the finite element mesh, and construct the training targets by evaluating the finite element solutions at these points. PCA-Net learns a map between the PCA coefficients of the input and output images, whereas FNO learns the map between the images directly. We train both models using a data loss. FC-PINO is a physics-informed neural operator~\cite{li2024physics} with the FNO architecture. It evaluates the same lifted FOSLS residual as RBNO on the grid, through Fourier continuation and spectral differentiation. We first train FC-PINO using the data loss, then continue training with the sum of the data loss, the residual loss evaluated over a larger set of unlabeled parameters, and a penalty for homogeneous boundary conditions.

For PCA-Net, FNO, and the two data-driven RBNO variants, we use the $500$ snapshots employed to construct the reduced basis for the semilinear diffusion--reaction benchmark and $2000$ labeled pairs for each of the other two benchmarks. For FC-PINO, we use $1000$ labeled pairs and the same number of unlabeled parameter samples as in RBNO training; these unlabeled samples require only residual evaluation. We summarize the architectures and training settings of PCA-Net, FNO, and FC-PINO in \cref{tab:architecture_training_details_pcanet,tab:architecture_training_details_fno,tab:architecture_training_details_fc_pino} in the appendix, and present the training histories of the data-driven RBNO variants and the grid-based neural operators in \cref{fig:baseline_training_rbno,fig:baseline_training}, respectively.

To compare grid-based predictions with RBNO, we follow~\cite{QCD} and use two evaluations, labeled $[\mathrm{CG}_1\times\mathrm{CG}_1]$ and $[\mathrm{RT}_1\times\mathrm{CG}_2]$ in \cref{tab:baseline_fe}. In the first, we interpret the predicted nodal values as the degrees of freedom of a $\mathrm{CG}_1\times\mathrm{CG}_1$ function and compare it with the corresponding nodal interpolant of the finite element solution. This reference is the representation that the grid-based neural operators are trained to reproduce. In the second, we compare the prediction with the finite element solution itself. For this comparison, we project the $\mathrm{CG}_1\times\mathrm{CG}_1$ prediction $(\hat{\sigma}^{\circ},\hat{u}^{\circ})$ onto $\HH_h$ using the inner product of $\HH=H(\rdiv)\times H^1$. Specifically, we seek $(\sigma^{\circ},u^{\circ})\in\HH_h$ such that
\begin{align}
\label{eq:hdiv_h1_projection}
    &(\sigma^{\circ},\tau)
    +(\rdiv{\sigma^{\circ}},\rdiv{\tau})
    +(u^{\circ},v)
    +(\nabla{u^{\circ}},\nabla{v}) \notag\\
    &\qquad =
    (\hat{\sigma}^{\circ},\tau)
    +(\rdiv{\hat{\sigma}^{\circ}},\rdiv{\tau})
    +(\hat{u}^{\circ},v)
    +(\nabla{\hat{u}^{\circ}},\nabla{v}),
    \qquad \forall\,(\tau,v)\in\HH_h,
\end{align}
where $\HH_h$ incorporates the homogeneous boundary conditions $\sigma^{\circ}\cdot n=0$ on $\Gamma_N$ and $u^{\circ}=0$ on $\Gamma_D$. Since continuous piecewise linear fields belong to $\mathrm{RT}_1\times\mathrm{CG}_2$, this projection leaves predictions satisfying these boundary conditions unchanged.

The two evaluations thus differ primarily in the reference solution. In both cases, we compute the residual loss using the same finite element procedure as for RBNO, with spatial derivatives evaluated from the finite element representation and integration performed by Gaussian quadrature on the mesh. The exception is the row labeled FC-PINO $[\mathrm{CG}_1\times\mathrm{CG}_1]^\dagger$, where we report the residual loss evaluated using the grid-based procedure employed during FC-PINO training. This procedure uses spectral differentiation with Fourier continuation for derivatives and the rectangle rule for integration on the uniform grid. These differences in differentiation and integration can lead to substantial discrepancies between the grid-based and finite element residual losses.

To isolate the error introduced by the grid representation, we also evaluate an idealized operator that reproduces the finite element solution exactly at every grid point. We report this result in the last row of each block of \cref{tab:baseline_fe}, labeled pointwise reconstruction. This comparison identifies the representation error that remains even in the absence of training and generalization errors.

\begin{table}[!tbp]
\centering\small
\sisetup{print-zero-exponent=true, exponent-mode=scientific, table-text-alignment=left, detect-weight=true}
\setlength{\tabcolsep}{4pt}
\begin{tabular}{@{}l
  S[table-format=1.2e-1]@{\,(}S[table-format=1.2e-1]@{)\quad}
  S[table-format=1.2e-1]@{\,(}S[table-format=1.2e-1]@{)\quad}
  S[table-format=1.2e-1]@{\,(}S[table-format=1.2e-1]@{)}@{}}\toprule
Method & \multicolumn{2}{l}{Relative $\LL^2$ error} & \multicolumn{2}{l}{Relative $\HH$ error} & \multicolumn{2}{l}{Residual loss}\\\midrule
\multicolumn{7}{@{}l}{\textit{Semilinear diffusion--reaction, $\kappa=10^3$}}\\
PCA-Net [CG$_1\times$CG$_1$] & 1.13e-01 & 2.31e-02 & 2.99e-01 & 5.51e-02 & 2.63e+01 & 3.99e+00\\
PCA-Net [RT$_1\times$CG$_2$] & 2.11e-01 & 1.31e-02 & 5.31e+00 & 3.75e-01 & 2.63e+01 & 3.99e+00\\
FNO [CG$_1\times$CG$_1$] & \bmnum{6.26e-03} & \bmnum{7.05e-04} & 3.76e-02 & 4.39e-03 & 2.74e+01 & 5.56e+00\\
FNO [RT$_1\times$CG$_2$] & 1.79e-01 & 7.98e-03 & 5.40e+00 & 5.09e-01 & 2.74e+01 & 5.56e+00\\
FC-PINO [CG$_1\times$CG$_1$]$^\dagger$ & 8.22e-02 & 4.06e-03 & 9.20e-01 & 6.01e-03 & 2.61e-02 & 4.32e-03\\
FC-PINO [RT$_1\times$CG$_2$] & 1.66e-01 & 6.99e-03 & 1.46e+00 & 1.13e-01 & 2.00e+00 & 3.37e-01\\
RBNO & 7.18e-03 & 1.97e-03 & 8.43e-03 & 2.64e-03 & \uwnum{1.35e-04} & \uwnum{3.96e-05}\\
RBNO ($\LL^2$ data loss) & 1.03e-02 & 3.32e-03 & 1.20e-02 & 2.97e-03 & 7.07e-04 & 5.90e-04\\
RBNO ($\HH$ data loss) & 9.55e-03 & 2.56e-03 & \bmnum{8.28e-03} & \bmnum{2.15e-03} & 3.38e-04 & 2.14e-04\\
\cmidrule{1-7}
Pointwise reconstruction & 1.79e-01 & 7.99e-03 & 5.44e+00 & 5.14e-01 & 2.78e+01 & 5.66e+00\\
\cmidrule{1-7}\morecmidrules\cmidrule{1-7}
\multicolumn{7}{@{}l}{\textit{Quasilinear diffusion}}\\
PCA-Net [CG$_1\times$CG$_1$] & 6.03e-02 & 1.76e-02 & 8.86e-02 & 2.69e-02 & 1.69e+00 & 6.27e-01\\
PCA-Net [RT$_1\times$CG$_2$] & 6.04e-02 & 1.76e-02 & 1.05e-01 & 8.67e-03 & 1.69e+00 & 6.27e-01\\
FNO [CG$_1\times$CG$_1$] & \bmnum{1.11e-02} & \bmnum{3.85e-03} & 1.55e-01 & 5.36e-02 & 5.65e+00 & 3.95e+00\\
FNO [RT$_1\times$CG$_2$] & 1.13e-02 & 3.82e-03 & 1.90e-01 & 5.08e-02 & 5.65e+00 & 3.95e+00\\
FC-PINO [CG$_1\times$CG$_1$]$^\dagger$ & 2.45e-02 & 6.26e-03 & 7.46e-02 & 2.20e-02 & 1.84e-02 & 2.62e-02\\
FC-PINO [RT$_1\times$CG$_2$] & 2.47e-02 & 6.24e-03 & 9.57e-02 & 9.11e-04 & 1.32e+00 & 5.13e-02\\
RBNO & 4.31e-02 & 9.27e-03 & \bmnum{8.33e-03} & \bmnum{3.30e-03} & \uwnum{2.12e-02} & \uwnum{4.28e-02}\\
RBNO ($\LL^2$ data loss) & 5.56e-02 & 1.70e-02 & 1.60e-02 & 6.19e-03 & 1.44e-01 & 3.24e-01\\
RBNO ($\HH$ data loss) & 6.24e-02 & 1.71e-02 & 1.31e-02 & 5.55e-03 & 7.18e-02 & 1.66e-01\\
\cmidrule{1-7}
Pointwise reconstruction & 1.89e-03 & 2.46e-04 & 1.22e-01 & 1.56e-02 & 2.20e+00 & 6.83e-01\\
\cmidrule{1-7}\morecmidrules\cmidrule{1-7}
\multicolumn{7}{@{}l}{\textit{Compressible hyperelasticity}}\\
PCA-Net [CG$_1\times$CG$_1$] & 3.64e-02 & 1.24e-02 & 1.61e-01 & 3.32e-02 & 4.40e-01 & 2.48e-01\\
PCA-Net [RT$_1\times$CG$_2$] & 3.67e-02 & 1.25e-02 & 4.19e-01 & 1.31e-01 & 4.40e-01 & 2.48e-01\\
FNO [CG$_1\times$CG$_1$] & \bmnum{6.35e-03} & \bmnum{2.23e-03} & 2.41e-01 & 4.67e-02 & 8.82e-01 & 3.40e+00\\
FNO [RT$_1\times$CG$_2$] & 7.47e-03 & 2.43e-03 & 4.88e-01 & 1.37e-01 & 8.82e-01 & 3.40e+00\\
FC-PINO [CG$_1\times$CG$_1$]$^\dagger$ & 1.45e-02 & 4.06e-03 & 2.50e-01 & 6.86e-02 & 6.73e-03 & 3.87e-03\\
FC-PINO [RT$_1\times$CG$_2$] & 1.48e-02 & 4.07e-03 & 2.99e-01 & 7.59e-02 & 2.19e-01 & 8.51e-02\\
RBNO & 2.60e-02 & 7.67e-03 & \bmnum{2.41e-02} & \bmnum{5.45e-03} & \uwnum{2.15e-03} & \uwnum{1.67e-03}\\
RBNO ($\LL^2$ data loss) & 2.48e-02 & 8.78e-03 & 2.79e-02 & 8.26e-03 & 8.48e-03 & 1.53e-02\\
RBNO ($\HH$ data loss) & 2.59e-02 & 8.81e-03 & 2.49e-02 & 7.24e-03 & 4.94e-03 & 8.57e-03\\
\cmidrule{1-7}
Pointwise reconstruction & 3.88e-03 & 1.24e-03 & 4.36e-01 & 1.36e-01 & 4.63e-01 & 2.59e-01\\
\bottomrule\end{tabular}
\caption{Mean and standard deviation (in parentheses) over the $500$ test samples of each benchmark of the relative errors in $\LL^2=L^2\times L^2$ and in $\HH=H(\rdiv)\times H^1$ norm, and of the residual loss. Unmarked rows and rows marked $[\mathrm{RT}_1\times\mathrm{CG}_2]$ are evaluated in $\mathrm{RT}_1\times\mathrm{CG}_2$ against the finite element solution. Grid-based predictions are first represented as $\mathrm{CG}_1\times\mathrm{CG}_1$ functions. Rows marked $[\mathrm{CG}_1\times\mathrm{CG}_1]$ measure the grid-based predictions in their own representation, against the finite element solution evaluated at the same grid points. The residual loss is the finite element loss $J_{h,q}$, except in the rows marked $\dagger$, which give the loss of FC-PINO in its own grid evaluation, with spectral differentiation and Fourier continuation for the derivatives and the rectangle rule for the integrals. Pointwise reconstruction: the finite element solution evaluated at the grid points and treated like the grid-based predictions. The smallest mean errors are in bold, and the smallest finite element loss is underlined with a wave; the rows marked $\dagger$ and the pointwise reconstruction are not ranked.}
\label{tab:baseline_fe}
\end{table}

RBNO achieves the smallest residual loss under finite element evaluation for all three benchmarks and the smallest $\HH$ error for quasilinear diffusion and compressible hyperelasticity. For semilinear diffusion--reaction, the $\HH$-data variant achieves a slightly smaller $\HH$ error but a larger residual loss. Overall, the two data-loss variants yield $\HH$ errors comparable to those of RBNO, although the $\LL^2$-data loss variant is less accurate in this norm, and both produce larger residual losses. FNO achieves the smallest $\LL^2$ error against the finite element solution sampled at the grid points in all three benchmarks, but its $\HH$ error exceeds that of RBNO.

When we evaluate the grid-based predictions against the full finite element solution in $\mathrm{RT}_1\times\mathrm{CG}_2$, their $\HH$ errors increase, particularly for semilinear diffusion--reaction, and become comparable to the pointwise reconstruction error. This indicates that the grid representation is a major limitation on their $\HH$ accuracy. In particular, continuous piecewise linear interpolation of the flux or stress can introduce large divergence errors. This effect is most pronounced for semilinear diffusion--reaction, where the conductivity and $f_1$ jump across the inclusion boundaries and only the normal component of the flux is continuous across these interfaces. Even exact reproduction of the solution at the grid points gives a relative $\HH$ reconstruction error of about $5.4$. The finite element residual losses of the grid-based predictions are correspondingly at least sixty times larger than those of RBNO, with the equilibrium component accounting for most of the loss. Residual training in FC-PINO reduces this discrepancy for semilinear diffusion--reaction, but its residual loss remains substantially larger than that of RBNO.

In the rows marked $\dagger$ in \cref{tab:baseline_fe}, FC-PINO has a small loss under its own evaluation procedure; for quasilinear diffusion, this loss is even smaller than the finite element loss of RBNO. We compare the corresponding pointwise residuals in \cref{fig:residual_maps}. Although both procedures discretize the same lifted FOSLS functional, they yield substantially different loss values, illustrating the effect of residual evaluation errors described in \cref{prop:consistency}. Thus, a physics-informed loss inherits the error control of the residual--error equivalence only when the residual is evaluated with sufficient accuracy.

\FloatBarrier
\section{Conclusions}\label{sec:conclusions}
A physics-informed loss is a sound basis for operator learning when its residual is equivalent to the error. For nonlinear parametric PDEs written as first-order systems, we showed that the least-squares residual of a conforming prediction has this property near the selected solution branch if the linearized residual is stable and the residual derivative is Lipschitz continuous, possibly in a norm stronger than the error norm, and that errors in the evaluation of the residual enter the estimate additively. We verified these conditions for three classes of nonlinearity, which differ in the set of states on which the estimate holds. For a monotone cubic reaction and for Kirchhoff-type diffusion, it is every bounded set, and for the reaction we obtained an error bound computable from the data and the prediction. For compressible Neo-Hookean hyperelasticity, it is a neighborhood of the branch, and the stability of the branch tangent remains an assumption.

In a conforming finite element space, we showed that the loss of a state near the finite element solution is the loss of that solution plus an excess equivalent to the squared distance from it. We trained RBNO on this loss with the loss and its gradient assembled by the finite element library; nonlinear solves were needed only for the snapshots of the reduced basis.

On one benchmark for each class, we observed that the ratio of the error to the square root of the loss excess has medians between $0.86$ and $1.31$ and stays between $0.43$ and $2.98$, comparable to the linear counterpart of the first benchmark ($0.41$--$3.22$). The discrete stability eigenvalues showed no loss of stability across meshes, solver iterates, and predictions; for hyperelasticity, they are our only evidence for the stability assumption. Grid-based operators evaluated in the same finite element space have large errors in the divergence of the flux, which the grid evaluation of the residual does not reveal. For the data of the first benchmark, the pointwise residual of the second-order equation also vanishes at states that violate the transmission conditions at the coefficient jumps, whereas the FOSLS residual bounds the $H^1$ error of the primal component of every state.

Several questions remain open. The loss does not select a branch, and for hyperelasticity we have no computable error bound. The explicit bound for the semilinear model combines worst-case constants: the lower conductivity bound enters through $\pmin^{-1}$, and the reaction through a factor proportional to $\kappa\norm{u}_{L^\infty}^2/\pmin$. A weighted constitutive residual would reduce the dependence on the conductivity, at the price of changing the trained loss. Transport-dominated or strongly localized solution sets are poorly represented by a single linear reduced space, and localized or nonlinear reduced representations~\cite{amsallem2012local,lee2020nonlinear} would require a compatible residual evaluation. A natural next step is to use the residual at inference to decide whether a prediction can be accepted, how it should be corrected, and when the model should be enriched.

\section*{Acknowledgments}

This work has been supported in part by NSF grants DMS-2245111, DMS-2245097, DMS-2513112, and DMS-2038080, by the U.S. Department of Energy under Contract No. DE-AC05-00OR22725 through the Genesis Mission project, and by the German Research Foundation grant CRC 1481.

\appendix
\setcounter{figure}{0}\setcounter{table}{0}%
\renewcommand{\theHfigure}{app.\arabic{figure}}\renewcommand{\theHtable}{app.\arabic{table}}%
\renewcommand{\thetheorem}{\Alph{section}.\arabic{theorem}}%
\aliascntresetthe{lemma}\aliascntresetthe{proposition}\aliascntresetthe{corollary}\aliascntresetthe{remark}
\section{Supplementary derivations}\label{app:theory}
\subsection{Embedding constants for the semilinear benchmark}\label{app:embedding}
\begin{lemma}\label{lem:embedding}
Let $\Omega=(0,1)^2$, $\Gamma_D=\{0,1\}\times[0,1]$, and $\U=H^1_{0,\Gamma_D}(\Omega)$. For $v\in\U$,
\[
 \norm{v}\le C_P\norm{\nabla v},\qquad
 \norm{v}_{L^4}\le C_4\norm{\nabla v},\qquad
 \norm{v}_{L^6}\le C_S\norm{\nabla v},
\]
with $C_P=\pi^{-1}$, $C_4=(C_P^3+2C_P^2)^{1/4}<0.70$, and $C_S=\bigl(\tfrac32C_4^4(3+C_P)\bigr)^{1/6}<1.03$.
\end{lemma}
\begin{proof}
By density we may assume that $v$ is smooth on $\bar\Omega$ and vanishes for $x_1\in\{0,1\}$. The Poincar\'e bound follows from the one-dimensional inequality $\int_0^1\varphi^2\le\pi^{-2}\int_0^1\varphi'^2$ for $\varphi(0)=\varphi(1)=0$, applied in $x_1$. For a function $\varphi$ on $(0,1)$ we use $\sup|\varphi|\le\tfrac12\int_0^1|\varphi'|$ if $\varphi(0)=\varphi(1)=0$, and $\sup|\varphi|\le\int_0^1|\varphi|+\int_0^1|\varphi'|$ in general. Since $|v(x_1,x_2)|^{2k}\le\sup_{x_1}|v(\cdot,x_2)|^k\,\sup_{x_2}|v(x_1,\cdot)|^k$,
\[
 \norm{v}_{L^{2k}}^{2k}\le\Bigl(\int_0^1\sup_{x_1}|v|^k\,dx_2\Bigr)\Bigl(\int_0^1\sup_{x_2}|v|^k\,dx_1\Bigr).
\]
For $k=2$, the first factor is at most $\int_\Omega|v||\partial_1v|\le\norm{v}\norm{\partial_1v}$ and the second at most $\norm{v}^2+2\norm{v}\norm{\partial_2v}$, which gives $\norm{v}_{L^4}^4\le(C_P^3+2C_P^2)\norm{\nabla v}^4$. For $k=3$, the first factor is at most $\tfrac32\int_\Omega v^2|\partial_1v|\le\tfrac32\norm{v}_{L^4}^2\norm{\partial_1v}$, and the second at most $\norm{v}_{L^3}^3+3\norm{v}_{L^4}^2\norm{\partial_2v}$ with $\norm{v}_{L^3}^3\le\norm{v}\norm{v}_{L^4}^2$. Inserting the bounds for $\norm{v}$ and $\norm{v}_{L^4}$ gives $\norm{v}_{L^6}^6\le\tfrac32C_4^4(C_P+3)\norm{\nabla v}^6$. Numerically $C_4^4\approx0.2349$ and $C_S^6\approx1.169$.
\end{proof}

\subsection{Regularized \texorpdfstring{$p$}{p}-Laplacian}\label{app:plaplacian}
For~\eqref{eq:nonlinear_plaplacian} put $\beta=(p-2)/2$ and $\varrho=1+|F|^2$. The derivative is $D\Psi_p(F)[H]=\varrho^{\beta}H+(p-2)\varrho^{\beta-1}(F\cdot H)F$. Its eigenvalue in directions orthogonal to $F$ is $\varrho^\beta$, and in the direction of $F$ it is $\varrho^{\beta-1}(1+(p-1)|F|^2)$; at $F=0$ both equal one. For $p\ge2$ both are at least one. For $1<p<2$ and $|F|\le M$ both are at least $(p-1)(1+M^2)^{(p-2)/2}$. Hence $D\Psi_p(F)[H]\cdot H\ge\alpha_p(M)|H|^2$ for $|F|\le M$, with $\alpha_p(M)=1$ for $p\ge2$ and $\alpha_p(M)=(p-1)(1+M^2)^{(p-2)/2}$ for $1<p<2$, and the tangent weak operator is coercive on sets with bounded gradients. The second derivative is
\[
 D^2\Psi_p(F)[H,K]=2\beta\varrho^{\beta-1}\bigl((F\cdot K)H+(F\cdot H)K+(H\cdot K)F\bigr)
 +4\beta(\beta-1)\varrho^{\beta-2}(F\cdot H)(F\cdot K)F ,
\]
so that on $|F|\le M$
\[
 \norm{D^2\Psi_p(F)}\le C_p^{(2)}(M):=6|\beta|M\max\{1,(1+M^2)^{\beta-1}\}+4|\beta(\beta-1)|M^3\max\{1,(1+M^2)^{\beta-2}\}.
\]
The derivative-Lipschitz constant is $\MK=\pmax C_p^{(2)}(M)$, the tangent flux is bounded by $L_T=\pmax\max\{1,p-1\}\max\{1,(1+M^2)^{(p-2)/2}\}$, and the tangent weak operator has $C_I=(1+C_P^2)/(\pmin\alpha_p(M))$. \Cref{lem:fosls_reduction} then gives explicit linearized constants, uniformly for $p$ in a compact subinterval of $(1,\infty)$.

\subsection{Neo-Hookean derivative bounds and the determinant calibration}\label{app:hyperelastic_derivatives}
Write $Z=F^{-\top}$ and $a=\lambda_L\log\det F-\mu$, so that $P(F)=\mu F+aZ$, $DZ[K]=-ZK^\top Z$, and $Da[K]=\lambda_L\operatorname{tr}(F^{-1}K)$. Differentiating the tangent~\eqref{eq:neo_hookean_tangent_note} gives
\begin{align}\label{eq:hyperelastic_second_derivative}
 D^2P(F)[H,K]={}&-\lambda_L\operatorname{tr}(F^{-1}KF^{-1}H)\,Z
 -\lambda_L\operatorname{tr}(F^{-1}H)\,ZK^\top Z-\lambda_L\operatorname{tr}(F^{-1}K)\,ZH^\top Z\nonumber\\
 &+a\bigl(ZK^\top ZH^\top Z+ZH^\top ZK^\top Z\bigr).
\end{align}
If $\norm{F^{-1}}_{\mathrm{op}}\le K_J$, each of the first three terms is bounded in the Frobenius norm by $\lambda_{L,\max}\sqrt d\,K_J^3\norm{H}_{\mathrm F}\norm{K}_{\mathrm F}$; for instance $|\operatorname{tr}(F^{-1}KF^{-1}H)|\le\norm{F^{-1}K}_{\mathrm F}\norm{F^{-1}H}_{\mathrm F}\le K_J^2\norm{K}_{\mathrm F}\norm{H}_{\mathrm F}$ and $\norm{Z}_{\mathrm F}\le\sqrt d\,K_J$. The last two terms together are bounded by $2(\mu_{\max}+\lambda_{L,\max}L_J)K_J^3\norm{H}_{\mathrm F}\norm{K}_{\mathrm F}$. This proves~\eqref{eq:hyperelastic_D2P_bound}, and~\eqref{eq:hyperelastic_DP_op_bound} follows in the same way from the three terms of~\eqref{eq:neo_hookean_tangent_note}. The constants are explicit rather than sharp.

For the calibration let $F=F_*+H$ with $\norm{H}_{\mathrm F}\le R_0$ pointwise. The first condition in~\eqref{eq:hyperelastic_calibration} gives $\norm{F_*+tH}_{\mathrm{op}}\le M_F$ for $0\le t\le1$. The derivative of the determinant is $D\det(F)[H]=\operatorname{cof}(F):H$, and $\norm{\operatorname{cof}(F)}_{\mathrm F}\le\sqrt d\,\norm{\operatorname{cof}(F)}_{\mathrm{op}}\le\sqrt d\,M_F^{d-1}$, since the singular values of $\operatorname{cof}(F)$ are products of $d-1$ singular values of $F$. Hence
\[
 \det(F_*+tH)\ge\det F_*-t\sqrt d\,M_F^{d-1}\norm{H}_{\mathrm F}\ge j_*/2\ge j_0 ,
\]
so the whole segment from the branch to the candidate is admissible. A segment between two candidates lies in the same convex neighborhood and satisfies the same bound. This, rather than convexity of a set of matrices with bounded determinant, justifies the mean-value estimate.

\subsection{Convergence of the discrete stability constant}\label{app:diagnostic}
Fix a linearization point $\bar u$ and a parameter $\pp$, write $\cA=\cA_{\bar u,\pp}$, and let $\lambda_{\min}=\inf_{t\ne0}\norm{\cA t}_{\LL^2}^2/\norm{t}_\HH^2$. For exact, unprojected residual norms on conforming spaces, $\lambda_{\min,h}\ge\lambda_{\min}$. Given $\varepsilon>0$ choose $t\ne0$ whose quotient is at most $\lambda_{\min}+\varepsilon$. If every state can be approximated, there are $t_h\in\HH_h$ with $t_h\to t$ in $\HH$, and boundedness of $\cA$ gives
\[
 \limsup_{h\to0}\lambda_{\min,h}\le\lim_{h\to0}\frac{\norm{\cA t_h}_{\LL^2}^2}{\norm{t_h}_\HH^2}\le\lambda_{\min}+\varepsilon .
\]
Letting $\varepsilon\to0$ proves convergence. The argument is pointwise in the state and the parameter; uniform convergence over a neighborhood or a parameter set would need additional compactness.

\section{Supplementary numerical results}\label{app:numerical_results}
In this appendix, we present additional results on reduced basis approximation in \ref{app:numerical_record}, describe the RBNO architecture and training settings in \ref{app:implementation}, and visualize RBNO predictions in \ref{app:ref_pred_diff}. We also provide the settings and training histories of the other neural operators in \ref{app:baselines} and compare the pointwise residuals of FC-PINO and RBNO in \ref{app:residual_maps}.

\subsection{Reduced basis approximation}\label{app:numerical_record}
\begin{figure}[!htbp]
\centering
\includegraphics[width=.45\linewidth]{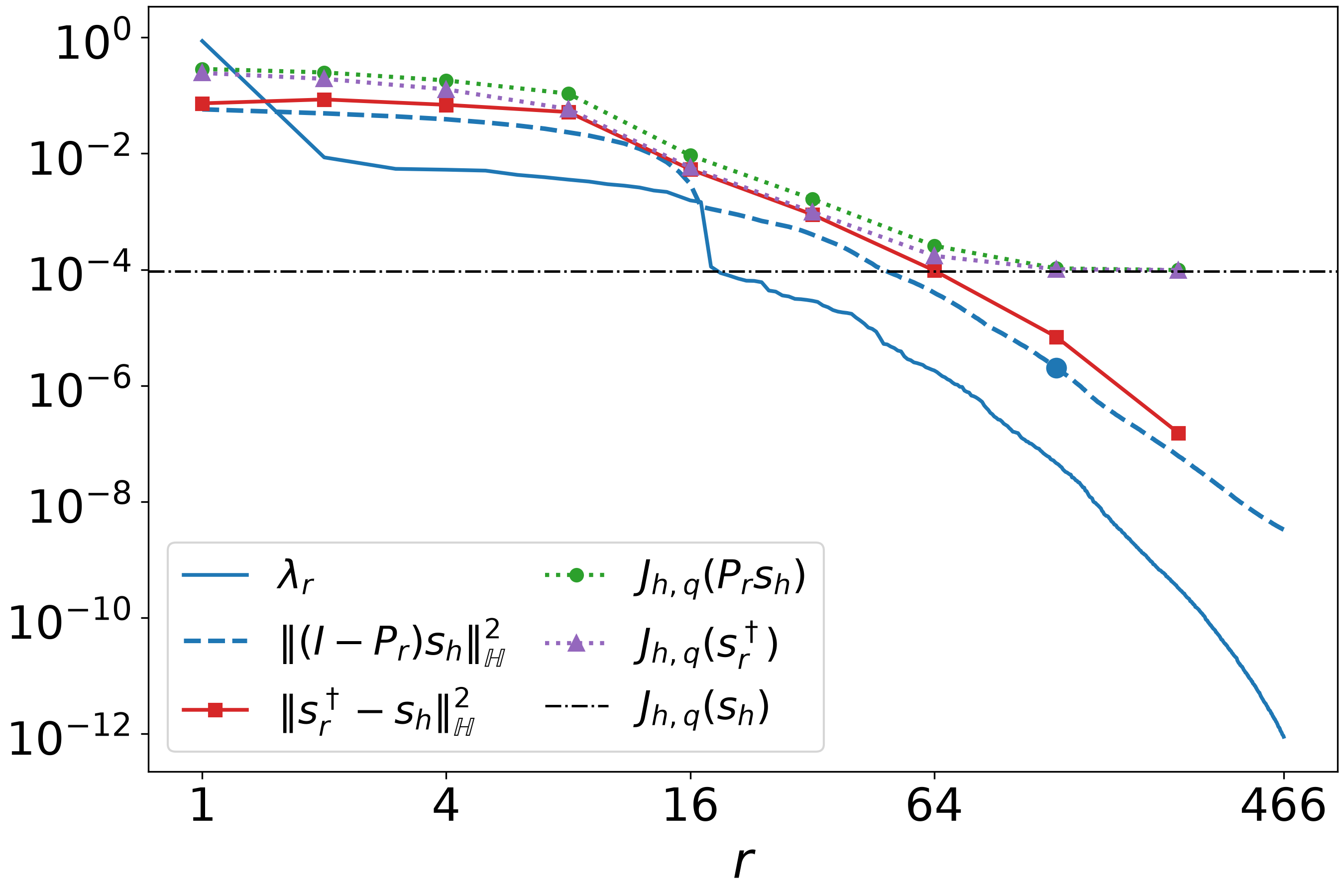}\hfill
\includegraphics[width=.45\linewidth]{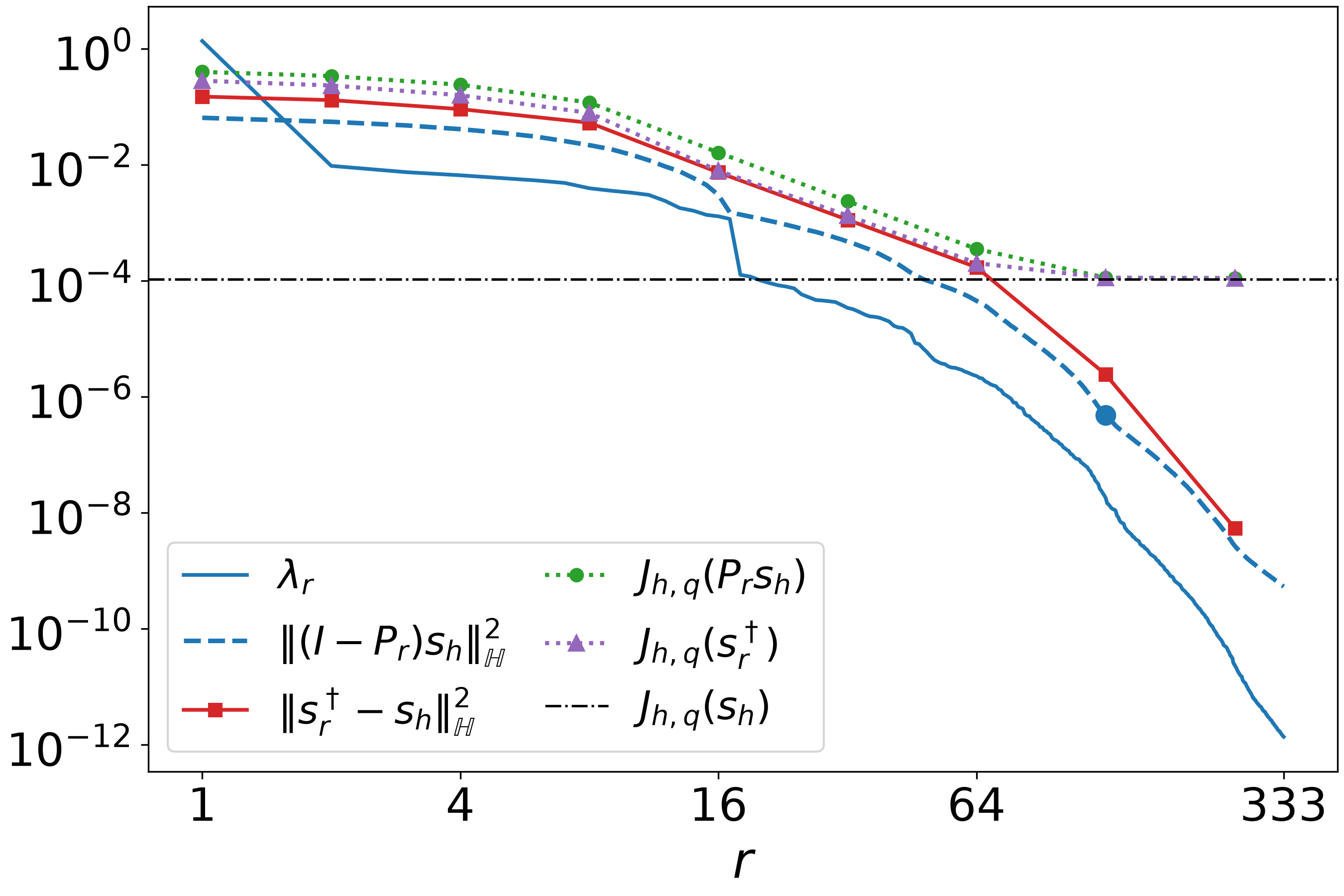}\\[4pt]
\includegraphics[width=.45\linewidth]{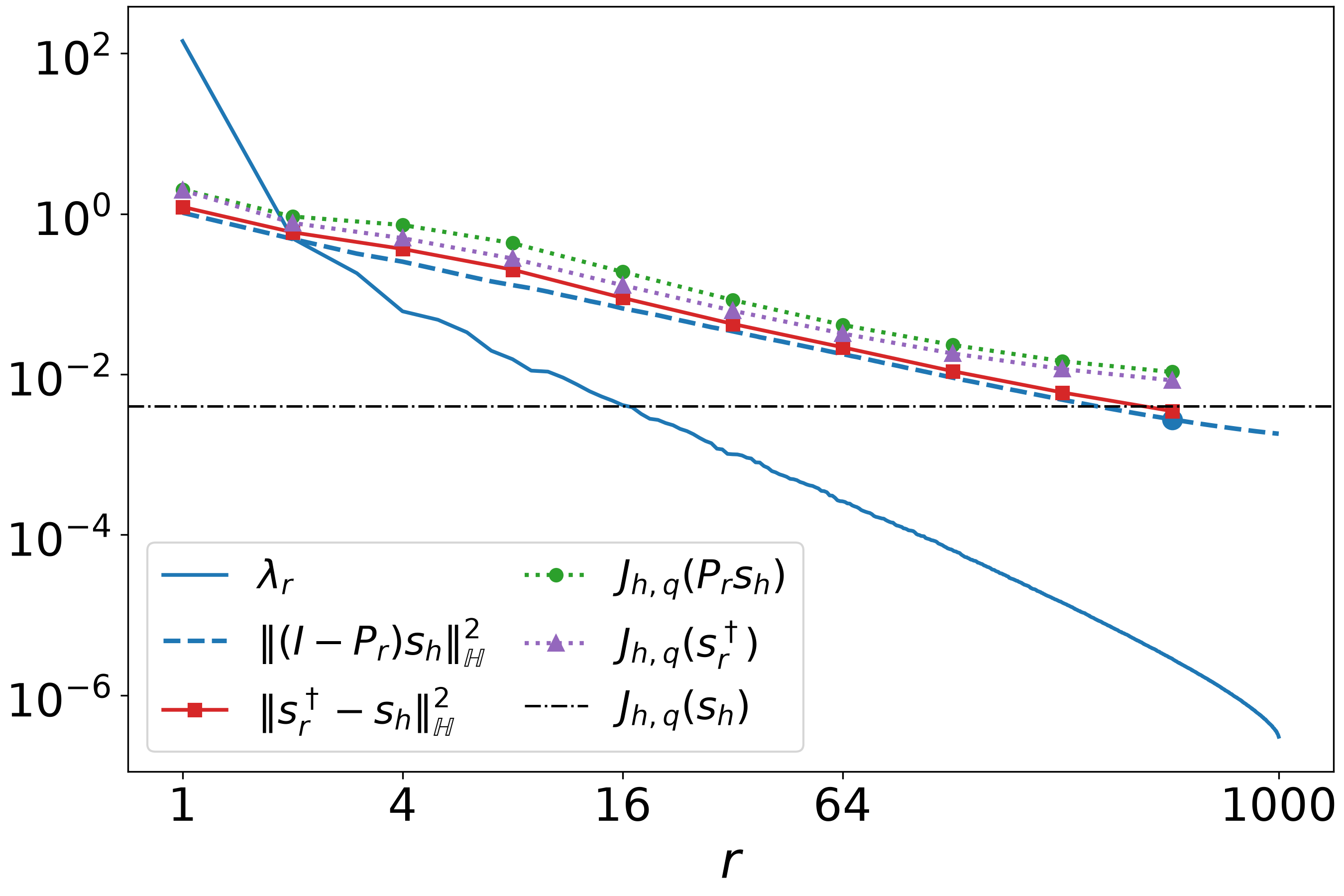}\hfill
\includegraphics[width=.45\linewidth]{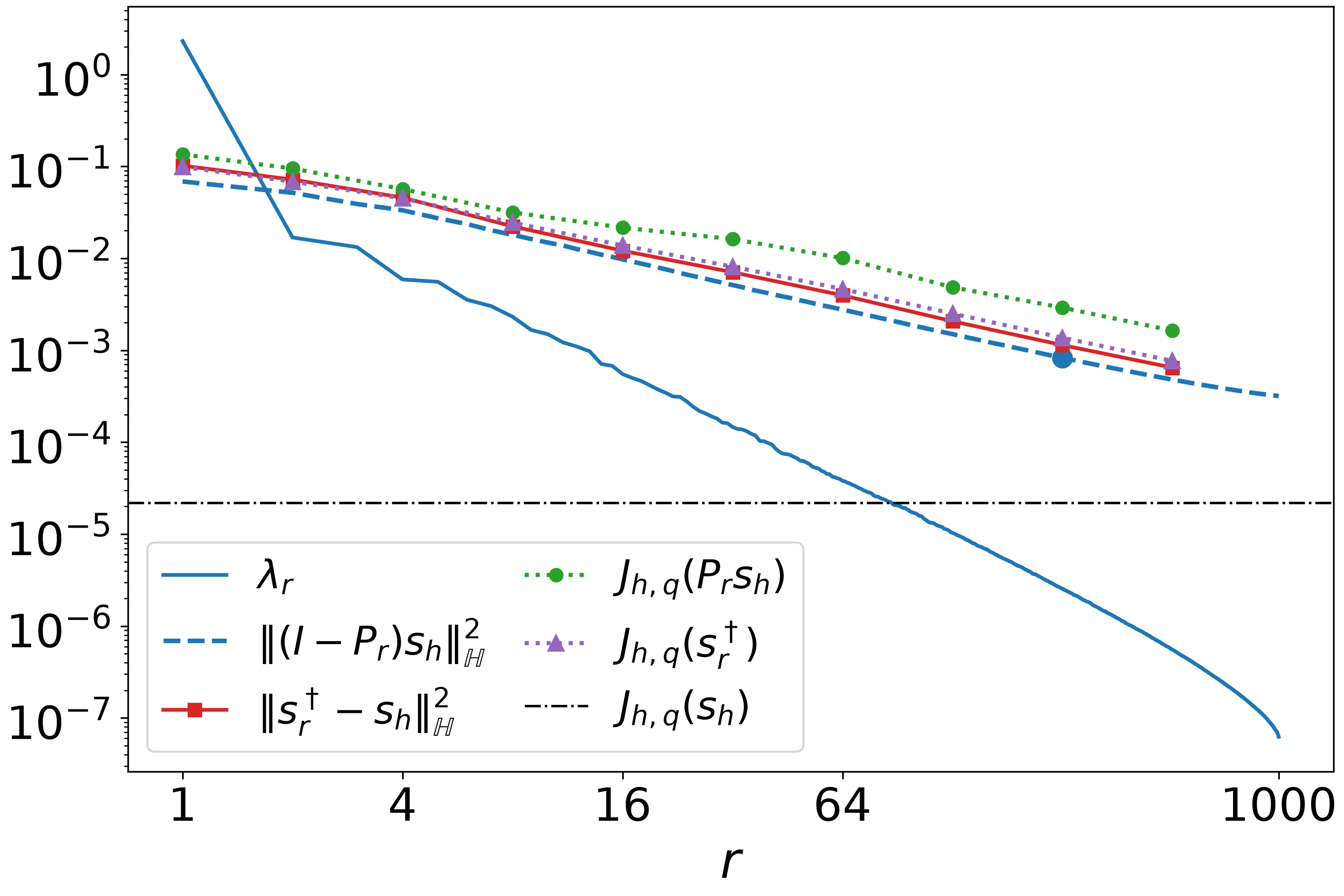}
\caption{Reduced basis approximation for $r=1,2,4,\dots$. We plot the POD eigenvalues $\lambda_r$ (solid), the squared projection error $\norm{(I-P_r)s_h}_\HH^2$ averaged over $500$ test solutions (dashed), and the squared error $\|s_r^\dagger-s_h\|_\HH^2$ of the reduced minimizer (squares) together with the losses $J_{h,q}(P_rs_h)$ (circles) and $J_{h,q}(s_r^\dagger)$ (triangles), each averaged over the first $100$ test parameters. The dash-dotted line shows the mean loss $J_{h,q}(s_h)$ of the finite element solutions. Top: semilinear diffusion--reaction with $\kappa=10^3$ (left) and $\kappa=0$ (right). Bottom: quasilinear diffusion (left) and compressible hyperelasticity (right).}\label{fig:pod_record}
\end{figure}

In \cref{fig:pod_record}, we show the decay of the POD eigenvalues, the projection error on the test set, and the error and the loss of the reduced minimizer $s_r^\dagger$ as functions of the reduced dimension.

\FloatBarrier
\subsection{RBNO architecture and training settings}\label{app:implementation}
\begin{figure}[!htbp]
\centering
\includegraphics[width=.38\linewidth]{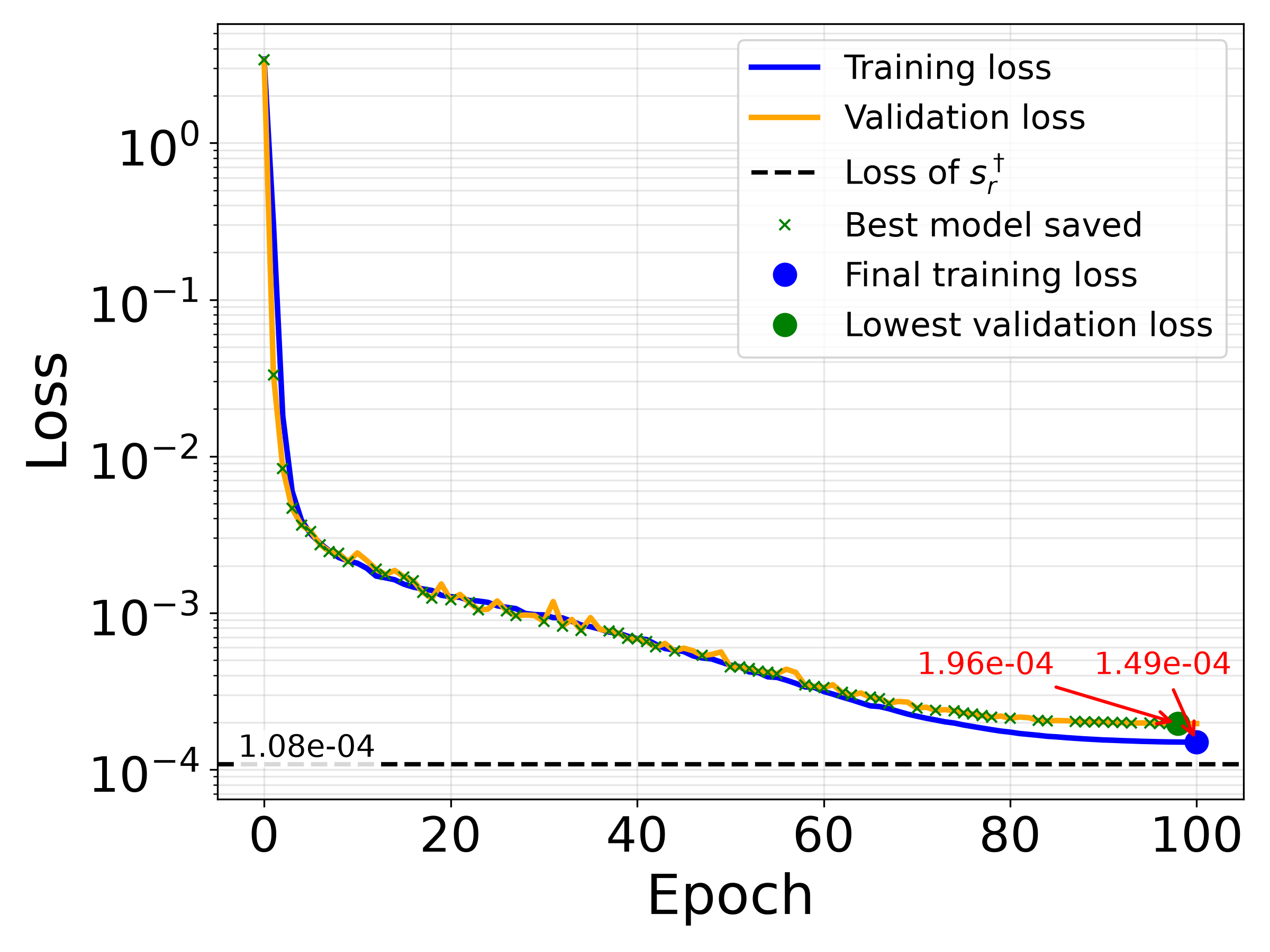}\hspace{2em}
\includegraphics[width=.38\linewidth]{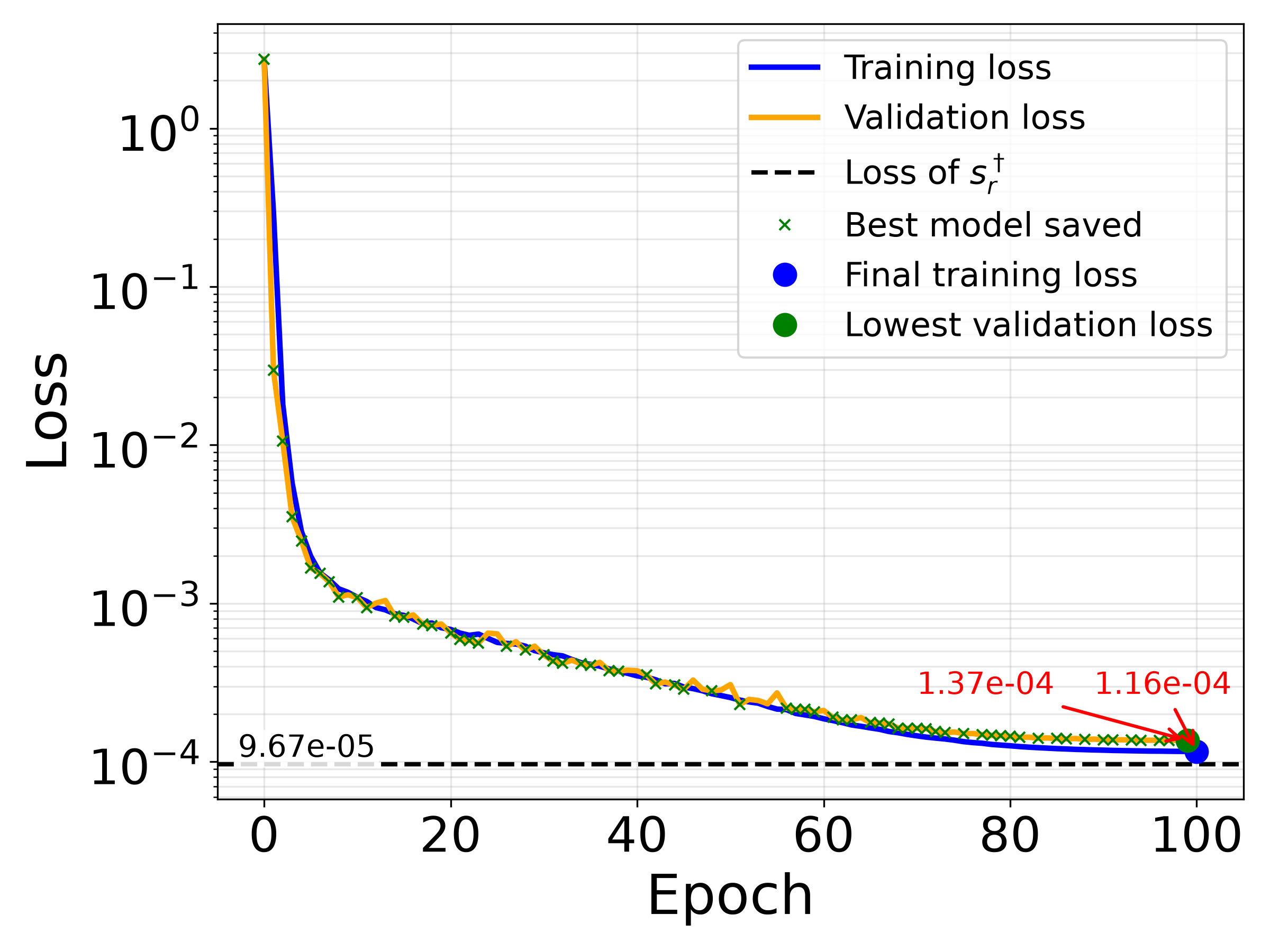}\\[4pt]
\includegraphics[width=.38\linewidth]{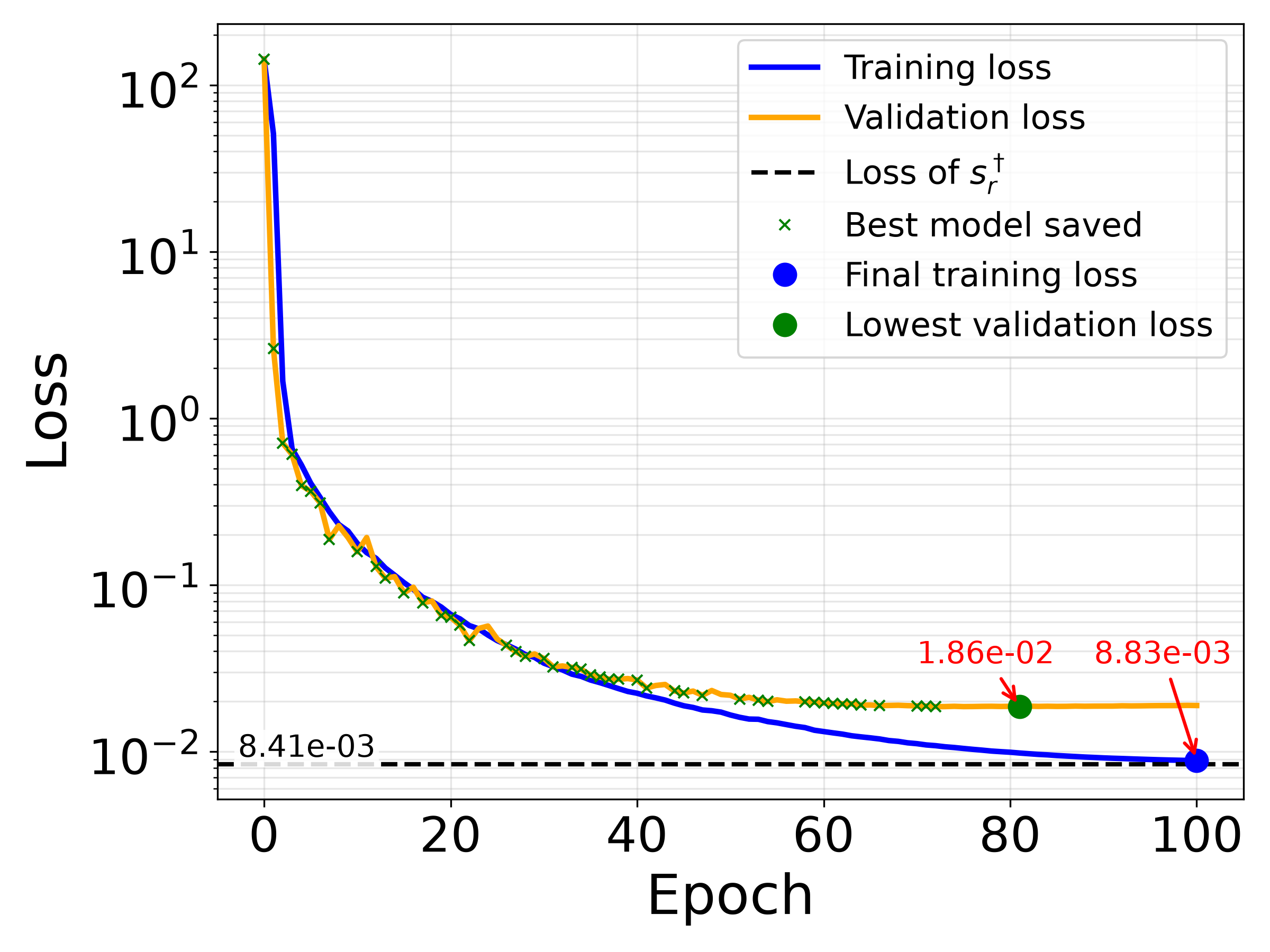}\hspace{2em}
\includegraphics[width=.38\linewidth]{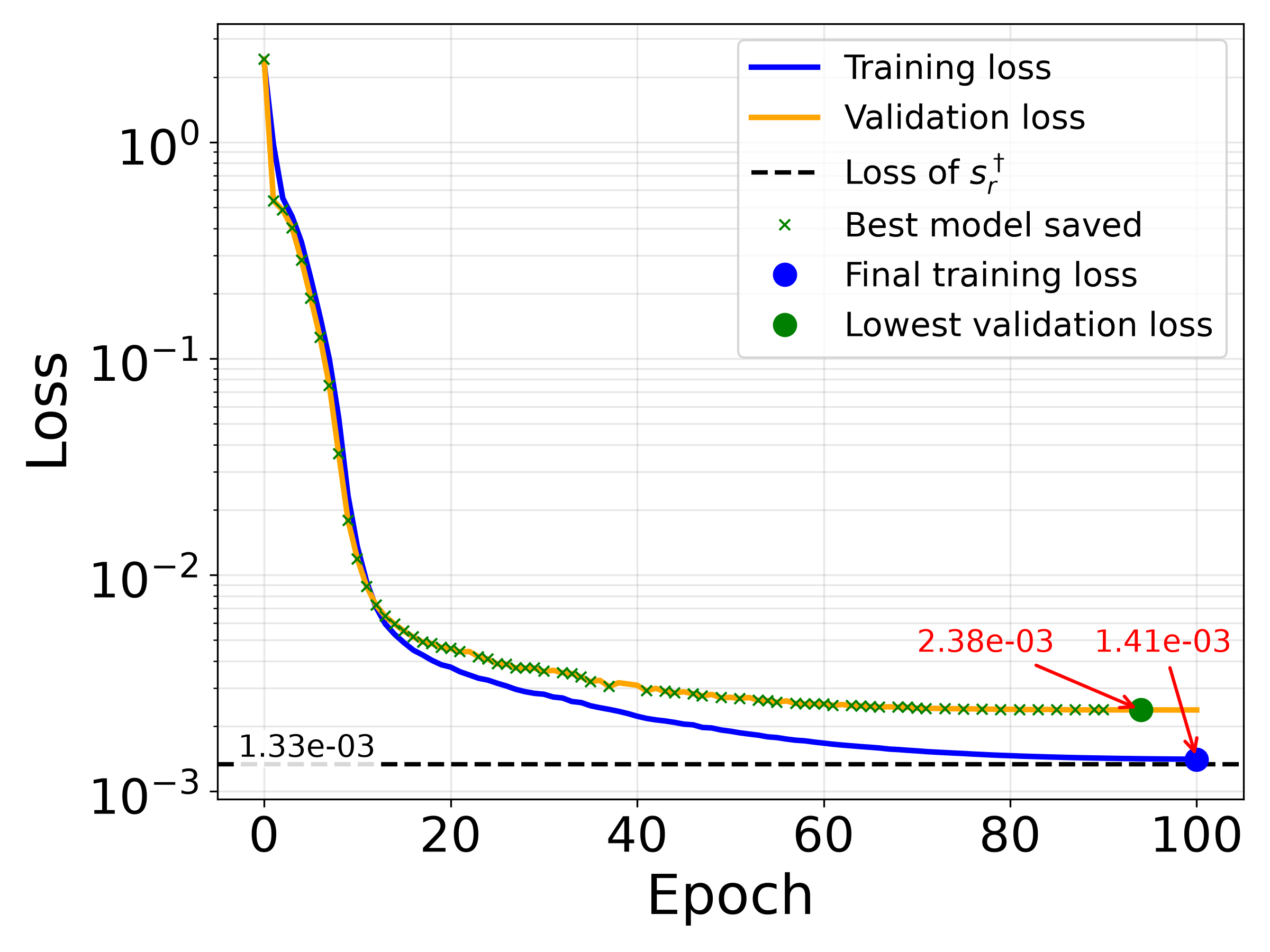}
\caption{Training and validation loss of RBNO per epoch: semilinear diffusion--reaction benchmark with $\kappa=0$ (top left) and $\kappa=10^3$ (top right), quasilinear diffusion benchmark (bottom left), and compressible hyperelasticity benchmark (bottom right). The dashed horizontal line is the mean loss of the reduced minimizers $s_r^\dagger$ of the test parameters.}\label{fig:training_record}
\end{figure}

We summarize the architectures and training settings of RBNO and its two data-loss variants in \cref{tab:architecture_training_details_rbno}, and show the training and validation losses of RBNO in \cref{fig:training_record}.

\begin{table}[!htbp]
\centering
\begin{scriptsize}
\begin{adjustbox}{max width=\linewidth}
\begin{tabular}{@{}ccc@{}}
\toprule
\textbf{Layer type} & \textbf{Configuration} & \textbf{Output shape (C, H, W)} \\
\midrule
\multicolumn{3}{c}{\textbf{Semilinear diffusion--reaction (SDR)}} \vspace{0.5em} \\
Input & 1 channel image ($\pp$) & (1, 129, 129) \\
Conv2d + act. & kernel=5, stride=2, out\_ch=64 & (64, 63, 63) \\
Conv2d + act. & kernel=5, stride=2, out\_ch=128 & (128, 30, 30) \\
Conv2d + act. & kernel=5, stride=3, out\_ch=256 & (256, 9, 9) \\
Conv2d + act. & kernel=3, stride=1, out\_ch=512 & (512, 7, 7) \\
Flatten & -- & (512$\times$7$\times$7) \\
Linear + act. & 512$\times$7$\times$7 $\rightarrow$ 512 & (512) \\
Linear & 512 $\rightarrow$ 128 & (128) \\
\# trainable parameters & \multicolumn{2}{c}{15.1 M} \\
\midrule
\multicolumn{3}{c}{\textbf{Quasilinear diffusion (QD)}} \vspace{0.5em} \\
Input & 1 channel image ($\pp$) & (1, 129, 129) \\
Conv2d + act. & kernel=3, stride=1, out\_ch=128 & (128, 129, 129) \\
Conv2d + act. & kernel=3, stride=1, out\_ch=256 & (256, 129, 129) \\
AvgPool2d & kernel=2 & (256, 64, 64) \\
Conv2d + act. & kernel=3, stride=1, out\_ch=512 & (512, 64, 64) \\
AvgPool2d & kernel=2 & (512, 32, 32) \\
Flatten & -- & (512$\times$32$\times$32) \\
Linear + act. & 512$\times$32$\times$32 $\rightarrow$ 512 & (512) \\
Linear & 512 $\rightarrow$ 512 & (512) \\
\# trainable parameters & \multicolumn{2}{c}{270.2 M} \\
\midrule
\multicolumn{3}{c}{\textbf{Compressible hyperelasticity (CH)}} \vspace{0.5em} \\
Input & 1 channel image ($\pp$) & (1, 65, 129) \\
Conv2d + act. & kernel=3, stride=1, out\_ch=128 & (128, 65, 129) \\
Conv2d + act. & kernel=3, stride=1, out\_ch=256 & (256, 65, 129) \\
AvgPool2d & kernel=2 & (256, 32, 64) \\
Conv2d + act. & kernel=3, stride=1, out\_ch=512 & (512, 32, 64) \\
AvgPool2d & kernel=2 & (512, 16, 32) \\
Flatten & -- & (512$\times$16$\times$32) \\
Linear + act. & 512$\times$16$\times$32 $\rightarrow$ 512 & (512) \\
Linear & 512 $\rightarrow$ 256 & (256) \\
\# trainable parameters & \multicolumn{2}{c}{135.8 M} \\
\midrule
initialization & \multicolumn{2}{c}{Xavier uniform (SDR), Kaiming uniform (QD, CH)} \\
activation & \multicolumn{2}{c}{LeakyReLU, slope 0.01 (SDR), GELU (QD, CH)} \\
input normalization & \multicolumn{2}{c}{none (SDR), z-score scaling (QD, CH)} \\
\midrule
\multicolumn{3}{c}{\textbf{Training}} \vspace{0.5em} \\
\# training / validation samples & \multicolumn{2}{c}{3000 / 500 (SDR), 8000 / 800 (QD), 5000 / 500 (CH)} \\
batch size & \multicolumn{2}{c}{50 (SDR), 200 (QD, CH)} \\
epochs & \multicolumn{2}{c}{100} \\
optimizer & \multicolumn{2}{c}{Muon~\cite{jordan2024muon} with auxiliary Adam} \\
 -- (Muon group) & \multicolumn{2}{c}{lr=$3{\times}10^{-3}$ (SDR; $5{\times}10^{-3}$ for $\kappa=0$), $10^{-2}$ (QD), $10^{-3}$ (CH)} \\
 & \multicolumn{2}{c}{momentum=0.95; weight\_decay=$10^{-4}$} \\
 -- (Adam group) & \multicolumn{2}{c}{lr as for the Muon group; betas=(0.9,0.95); weight\_decay=$10^{-4}$} \\
lr scheduler & \multicolumn{2}{c}{CosineAnnealingLR ($T_{\max}=$ \# epochs; $\eta_{\min}=$ lr$/100$ (SDR), $5{\times}10^{-4}$ (QD), $2.5{\times}10^{-5}$ (CH))} \\
numerical precision & \multicolumn{2}{c}{neural network: float32; POD reconstruction and residual loss: float64} \\
\midrule
\multicolumn{3}{c}{\textbf{Data-loss ablations (same network, optimizer, and scheduler)}} \vspace{0.5em} \\
loss & \multicolumn{2}{c}{$\norm{s-s_h}^2/\norm{s_h}^2$ in the $\LL^2$ or the $\HH$ norm} \\
\# training / validation pairs & \multicolumn{2}{c}{500 / 500 (SDR), 2000 / 200 (QD, CH)} \\
\bottomrule
\end{tabular}
\end{adjustbox}
\end{scriptsize}
\caption{Architecture and training settings of RBNO and of its two data-loss variants.}
\label{tab:architecture_training_details_rbno}
\end{table}

\FloatBarrier
\subsection{Visualization of the predictions of RBNO}\label{app:ref_pred_diff}
\begin{figure}[!htbp]
\centering
\includegraphics[width=\linewidth]{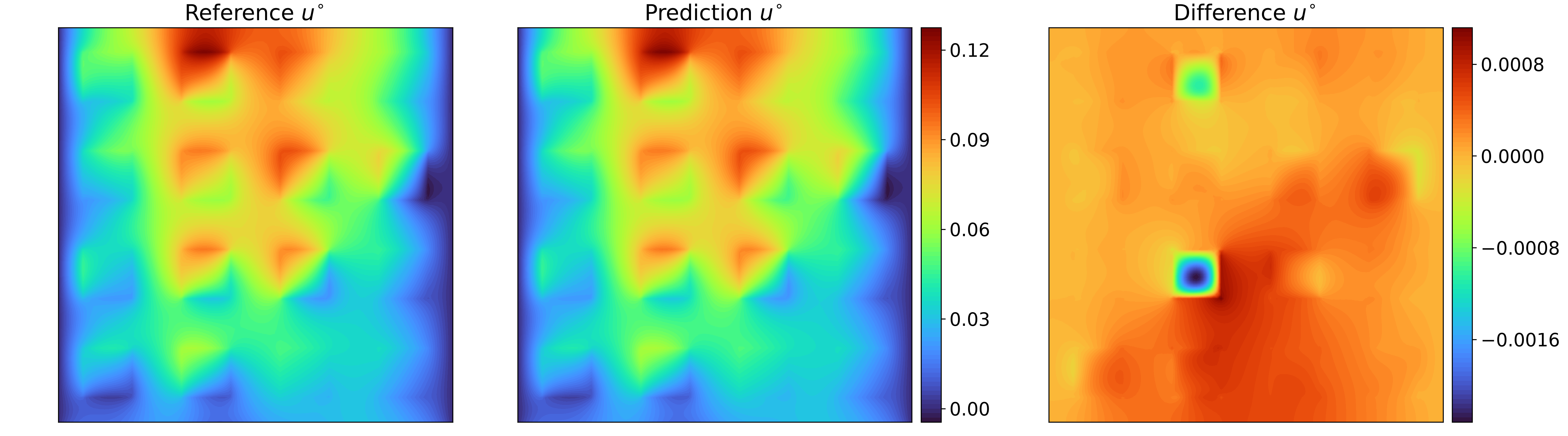}
\includegraphics[width=\linewidth]{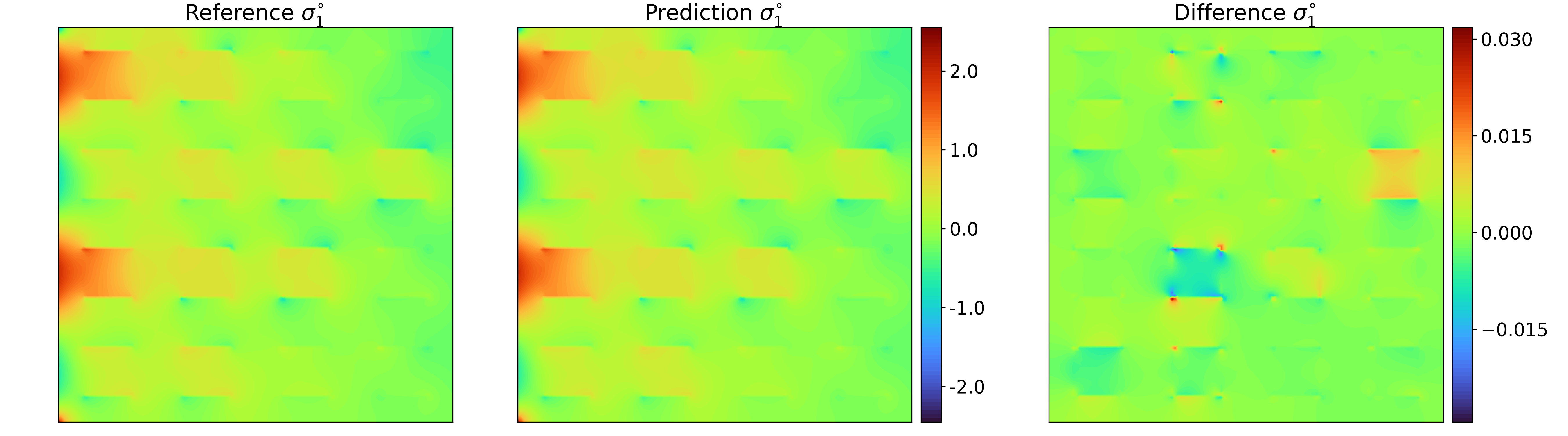}
\includegraphics[width=\linewidth]{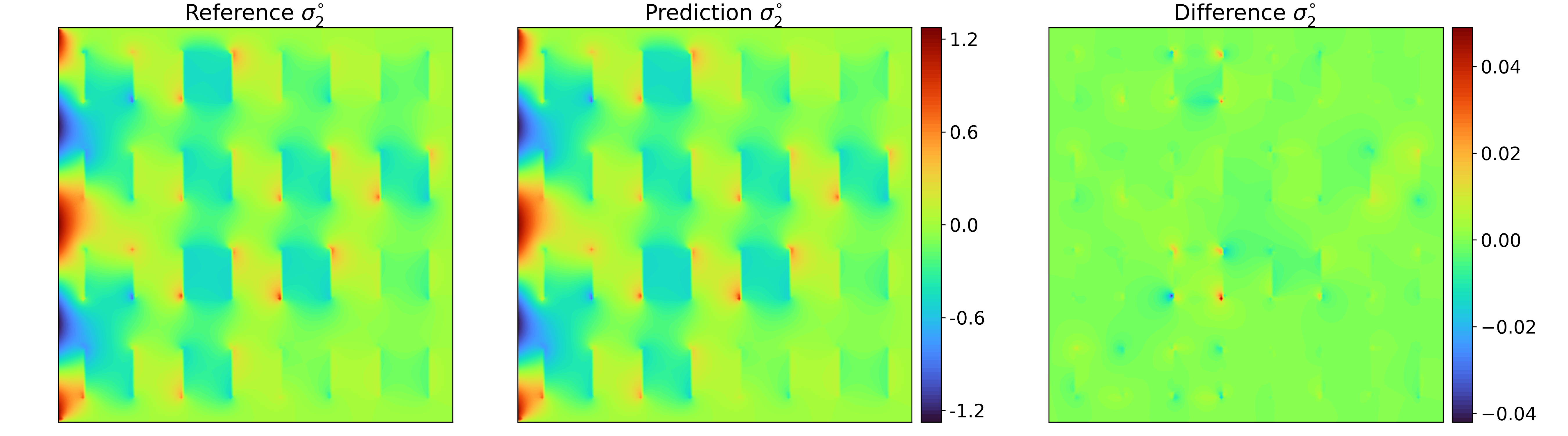}
\caption{Reference, RBNO prediction, and difference (left to right) of $u^{\circ}$ (top) and of the flux components $\sigma_1^{\circ}$ and $\sigma_2^{\circ}$ (middle, bottom) for the semilinear diffusion--reaction benchmark ($\kappa=10^3$).}\label{fig:ref_pred_diff_semilinear}
\end{figure}

\begin{figure}[!htbp]
\centering
\includegraphics[width=\linewidth]{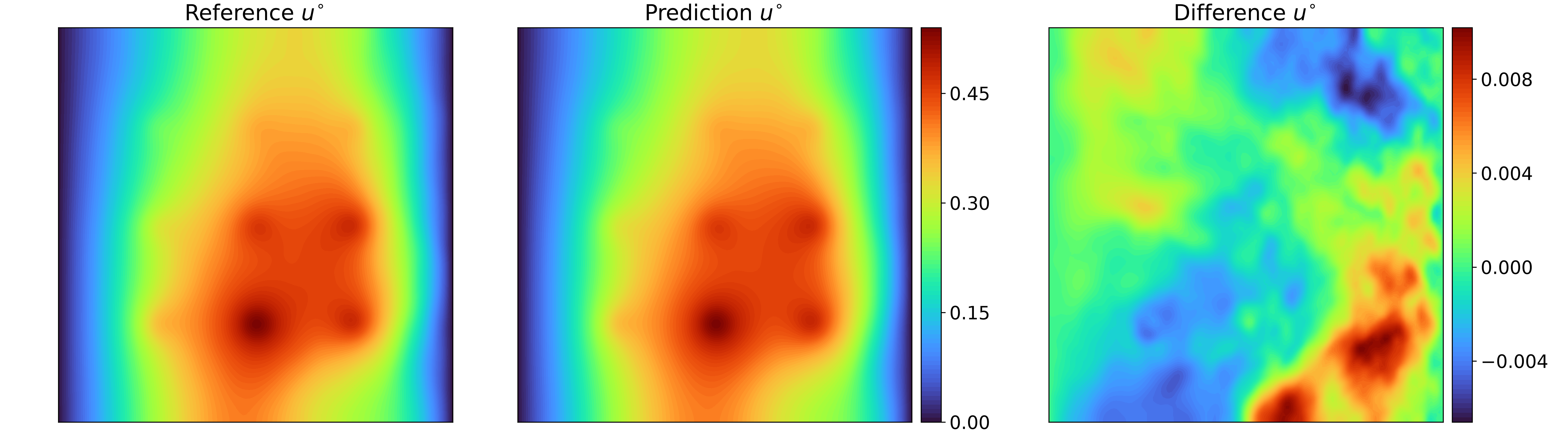}
\includegraphics[width=\linewidth]{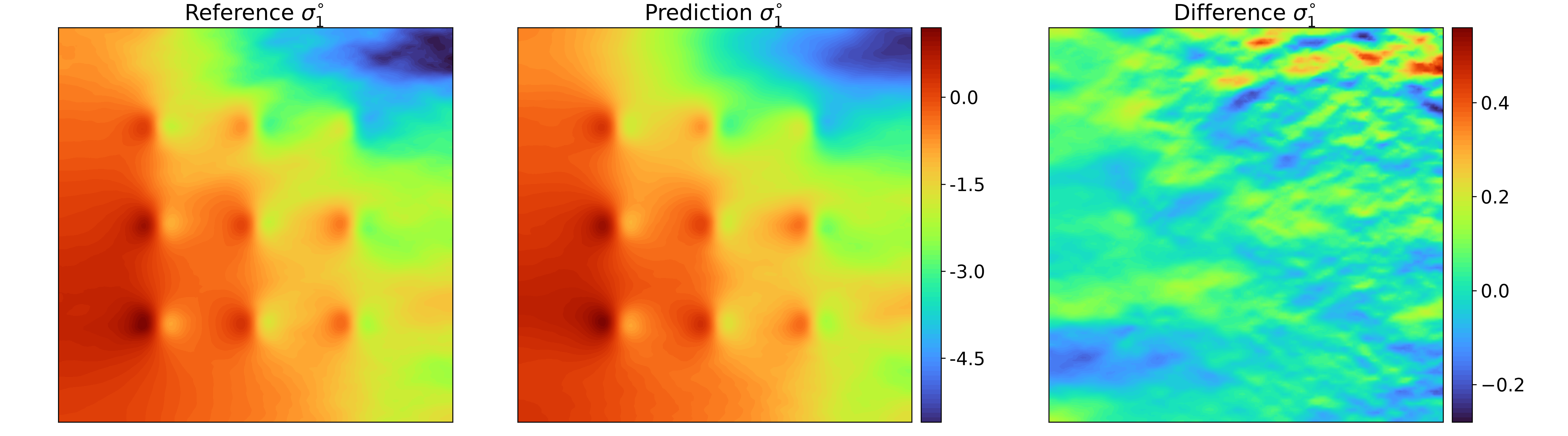}
\includegraphics[width=\linewidth]{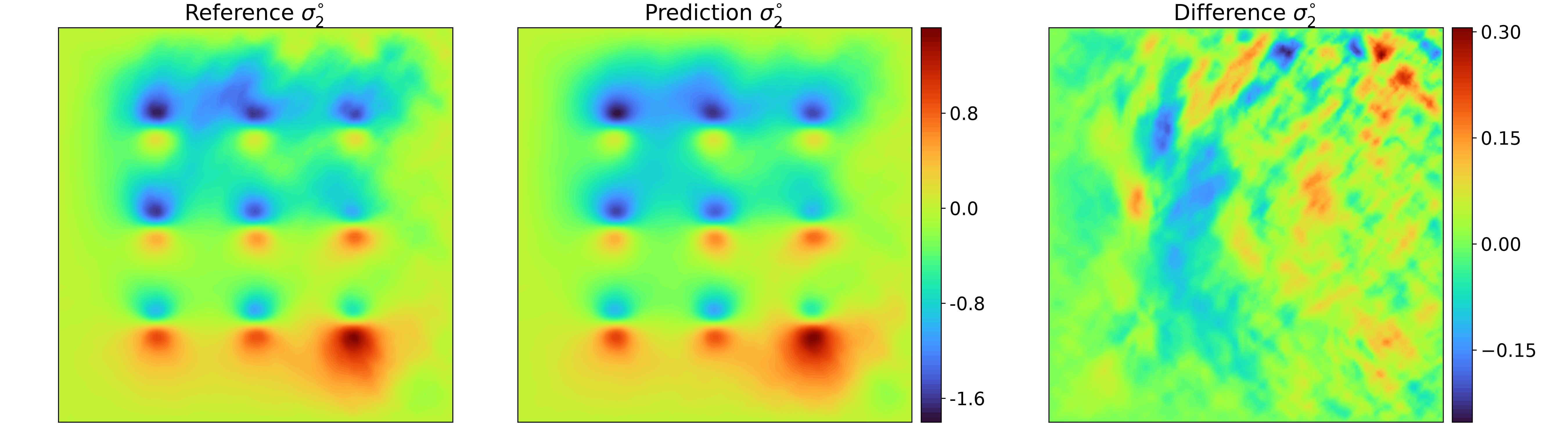}
\caption{Reference, RBNO prediction, and difference (left to right) of $u^{\circ}$ (top) and of the flux components $\sigma_1^{\circ}$ and $\sigma_2^{\circ}$ (middle, bottom) for the quasilinear diffusion benchmark.}\label{fig:ref_pred_diff_quasilinear}
\end{figure}

\begin{figure}[!htbp]
\centering
\includegraphics[width=0.45\linewidth]{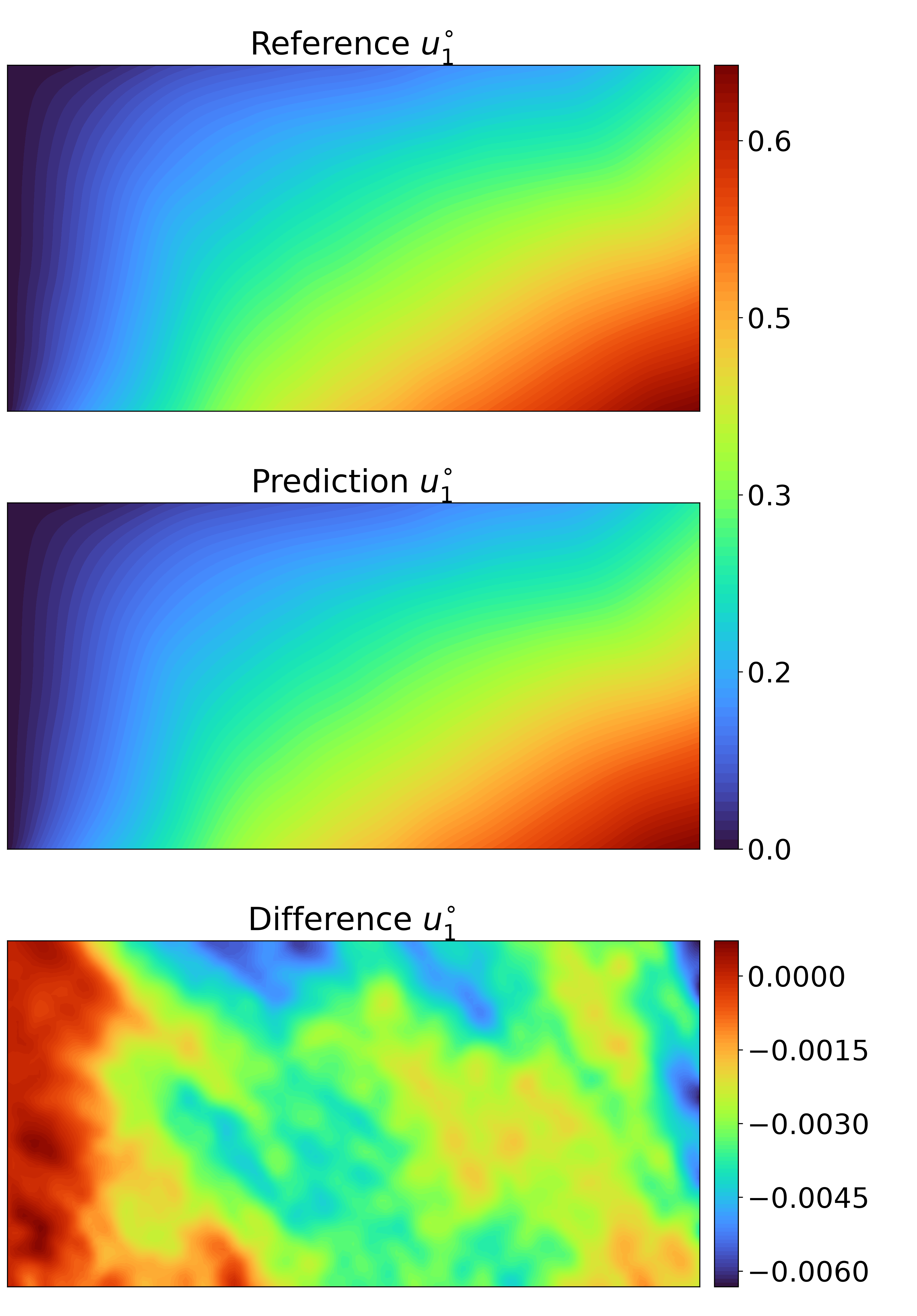}\hspace{0.04\linewidth}
\includegraphics[width=0.45\linewidth]{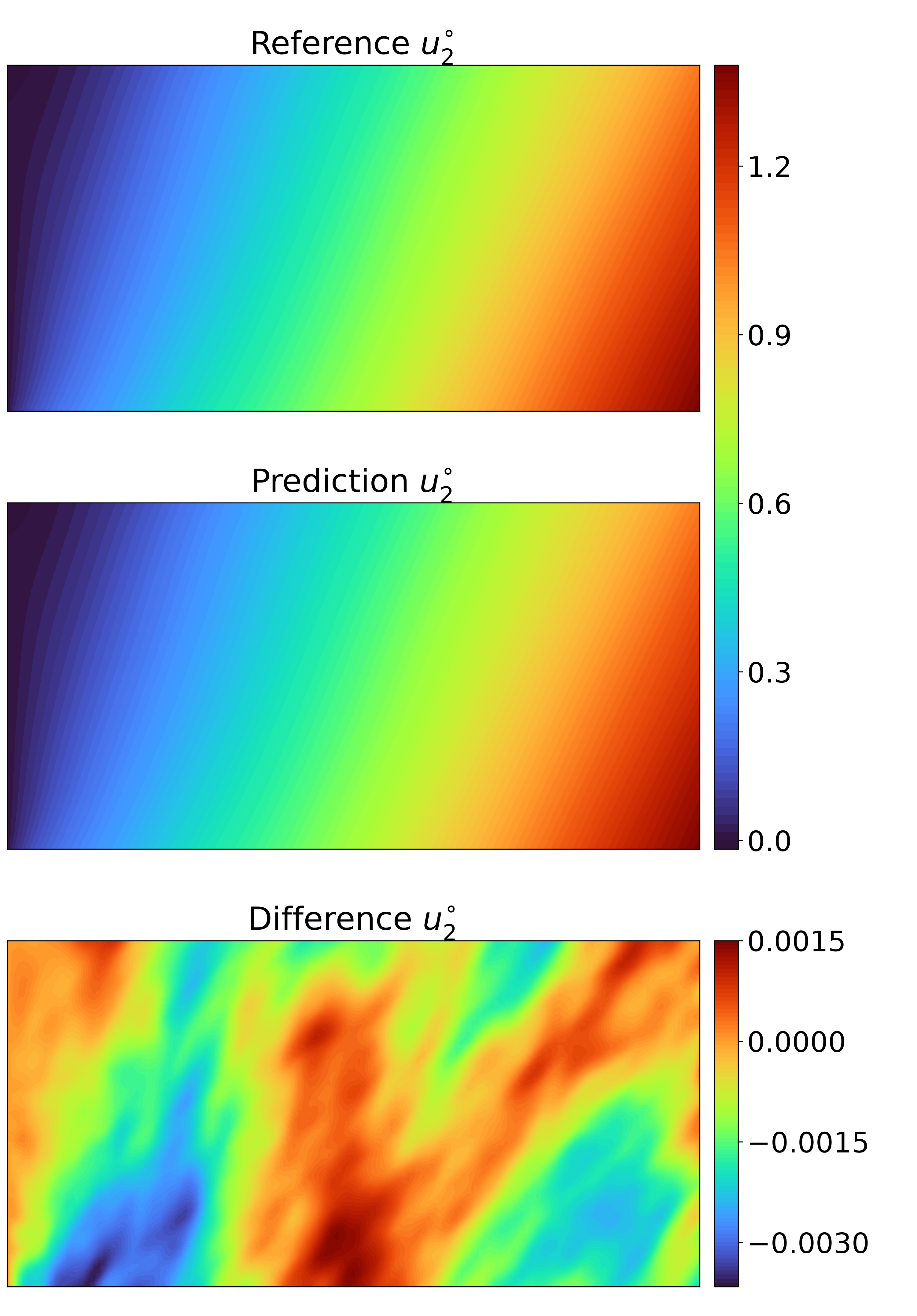}\\[6pt]
\includegraphics[width=0.23\linewidth]{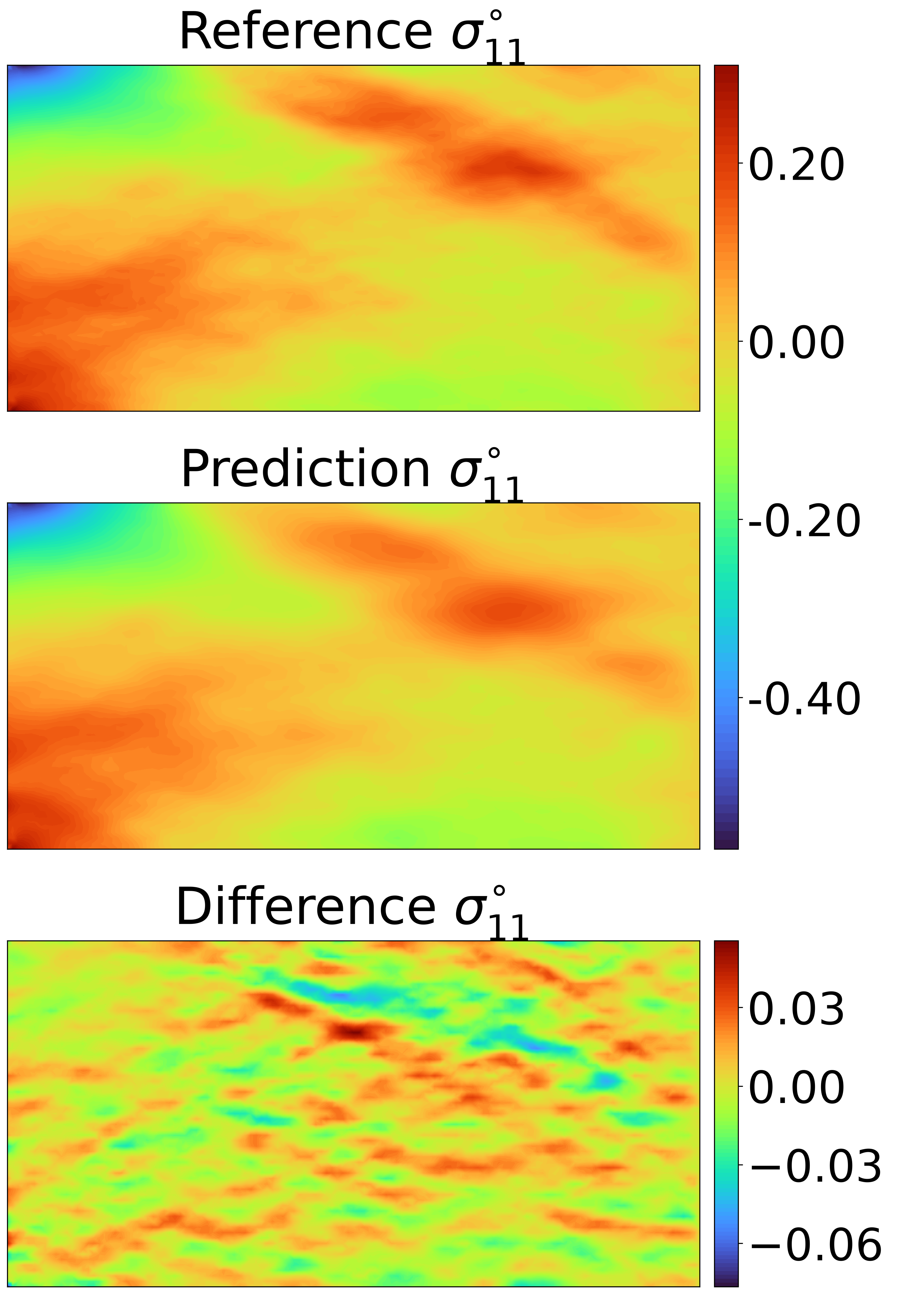}\hfill
\includegraphics[width=0.23\linewidth]{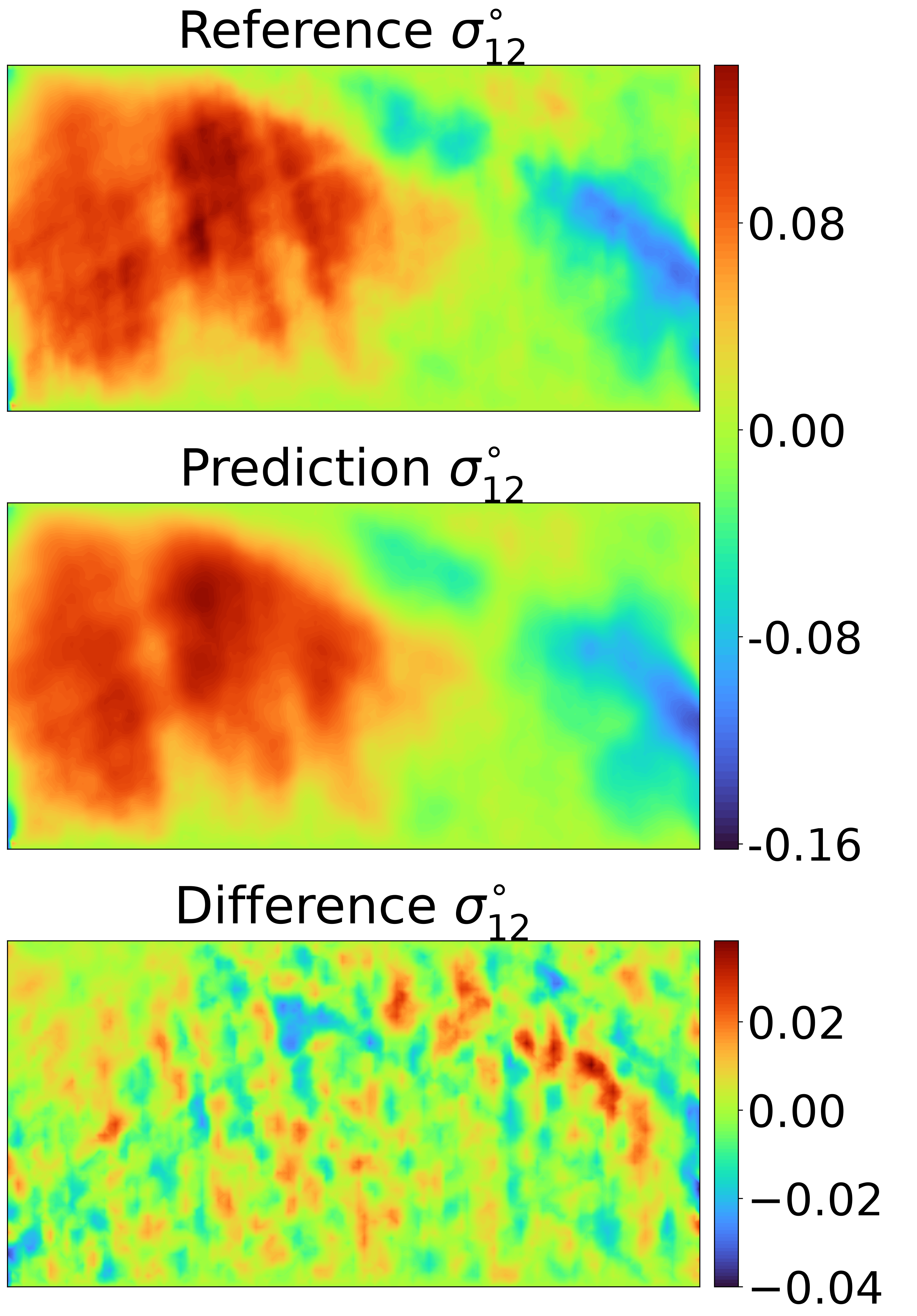}\hfill
\includegraphics[width=0.23\linewidth]{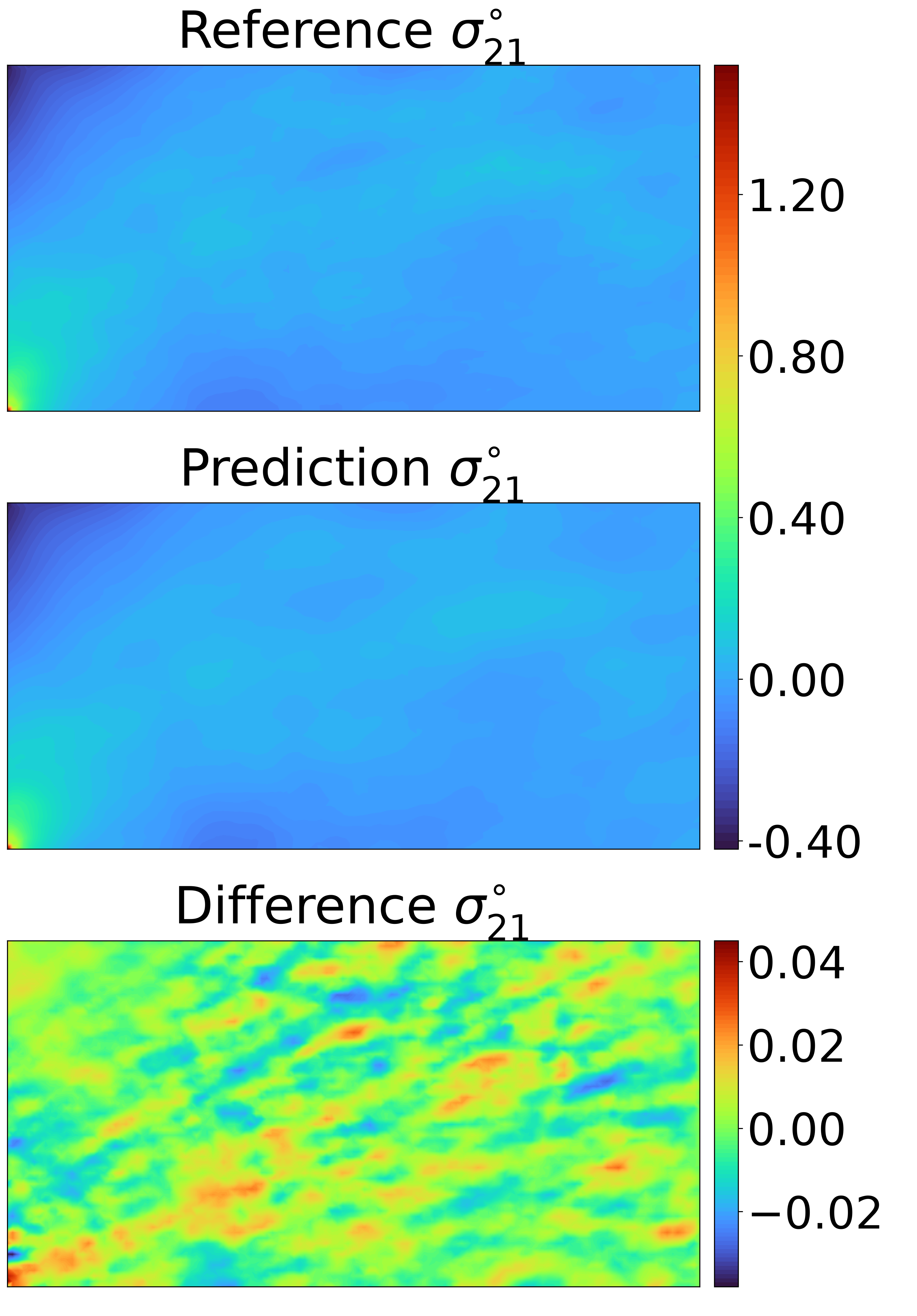}\hfill
\includegraphics[width=0.23\linewidth]{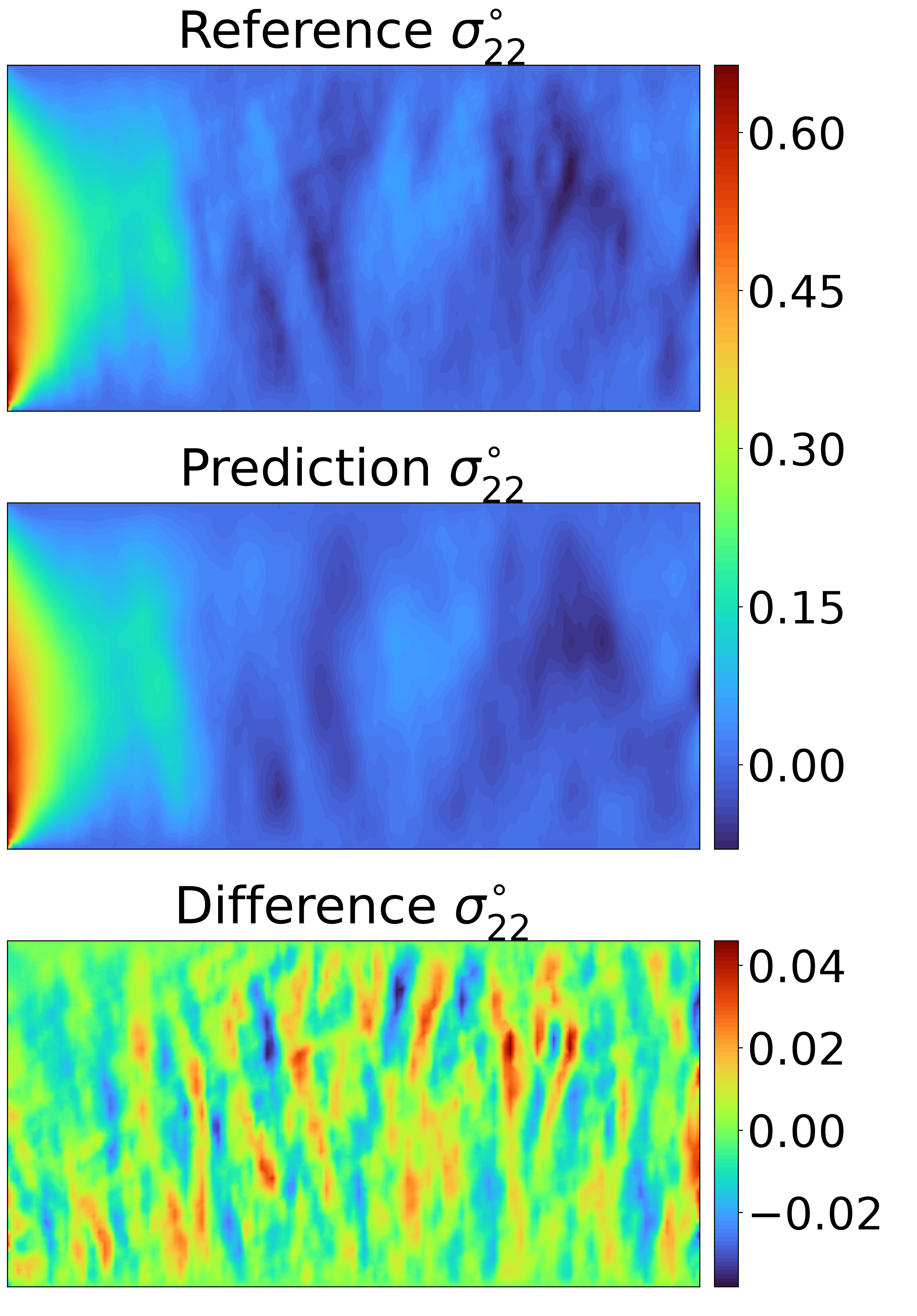}
\caption{Reference, RBNO prediction, and difference (top to bottom) of the displacement components $u_1^{\circ}$ and $u_2^{\circ}$ (top row) and of the stress components $\sigma_{11}^{\circ}$, $\sigma_{12}^{\circ}$, $\sigma_{21}^{\circ}$, $\sigma_{22}^{\circ}$ (bottom row, left to right) for the compressible hyperelasticity benchmark.}\label{fig:ref_pred_diff_hyperelastic}
\end{figure}

\Cref{fig:ref_pred_diff_semilinear,fig:ref_pred_diff_quasilinear,fig:ref_pred_diff_hyperelastic} compare the finite element solution (reference) with the RBNO prediction for the test parameters of \cref{fig:benchmarks}. All functions are evaluated on the uniform grid, so the figures show the $\mathrm{CG}_1\times\mathrm{CG}_1$ representation of the states.

\FloatBarrier
\subsection{Other neural operators}\label{app:baselines}
\begin{figure}[!htbp]
\centering
\includegraphics[width=.32\linewidth]{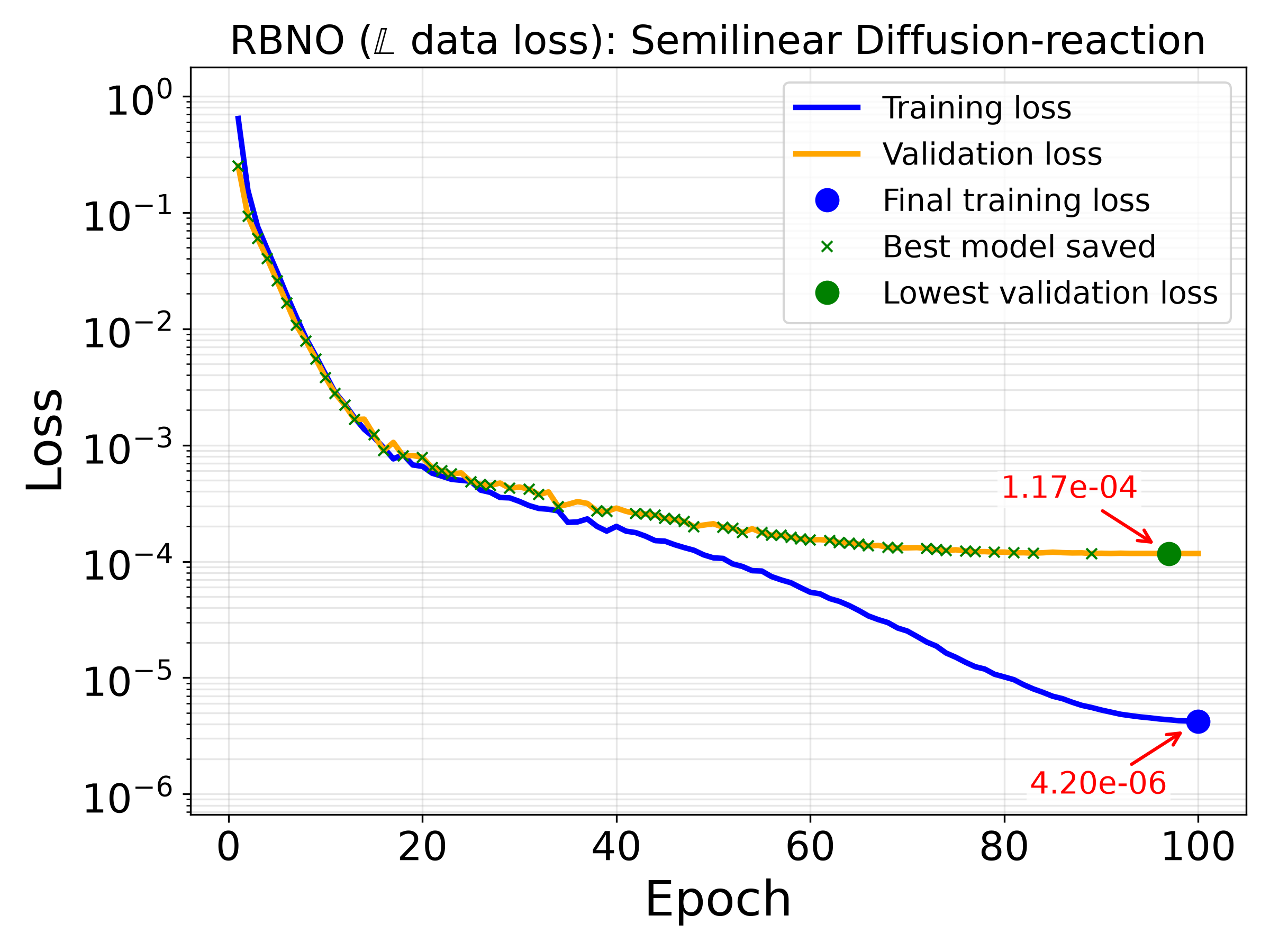}\hfill
\includegraphics[width=.32\linewidth]{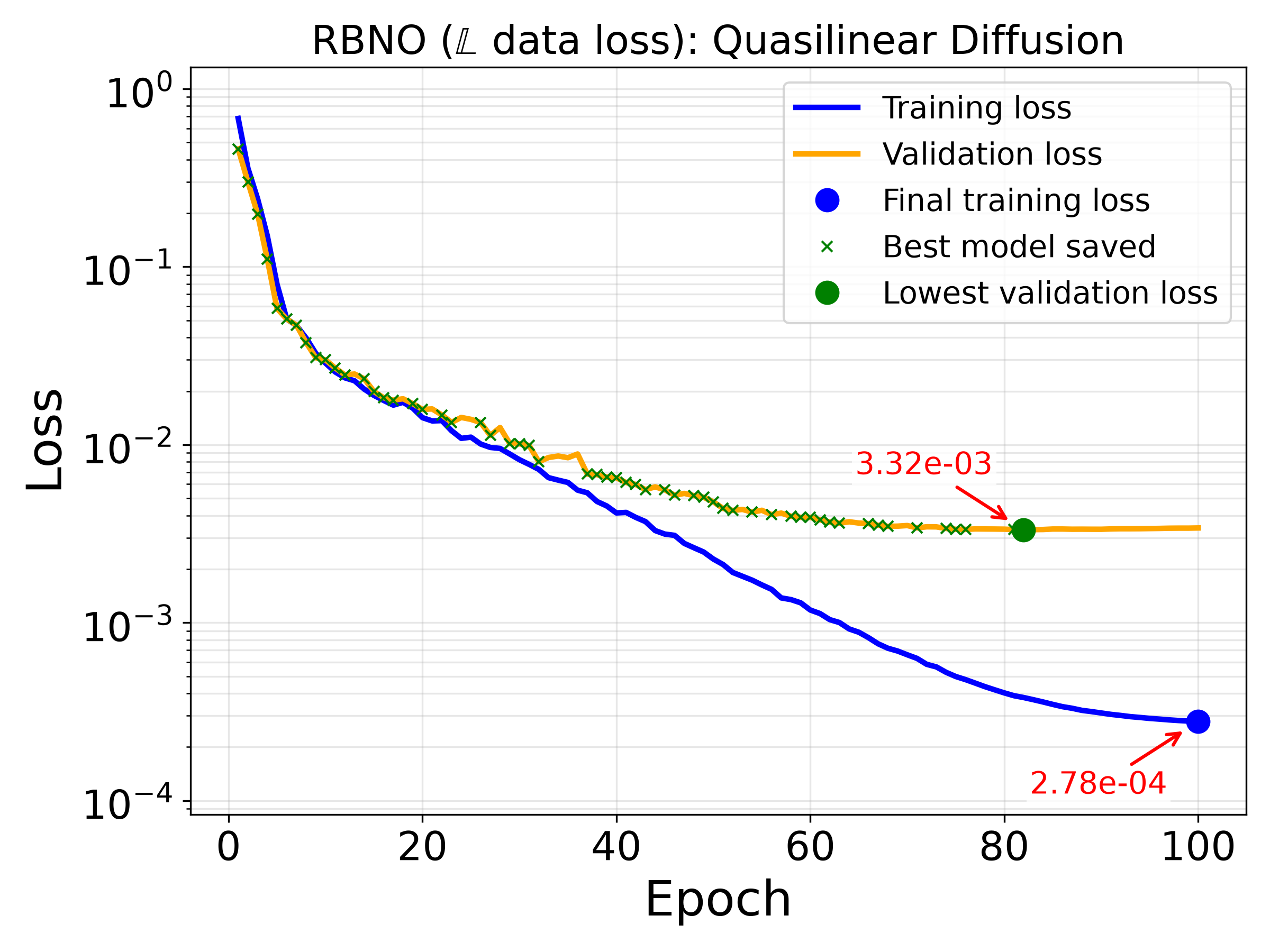}\hfill
\includegraphics[width=.32\linewidth]{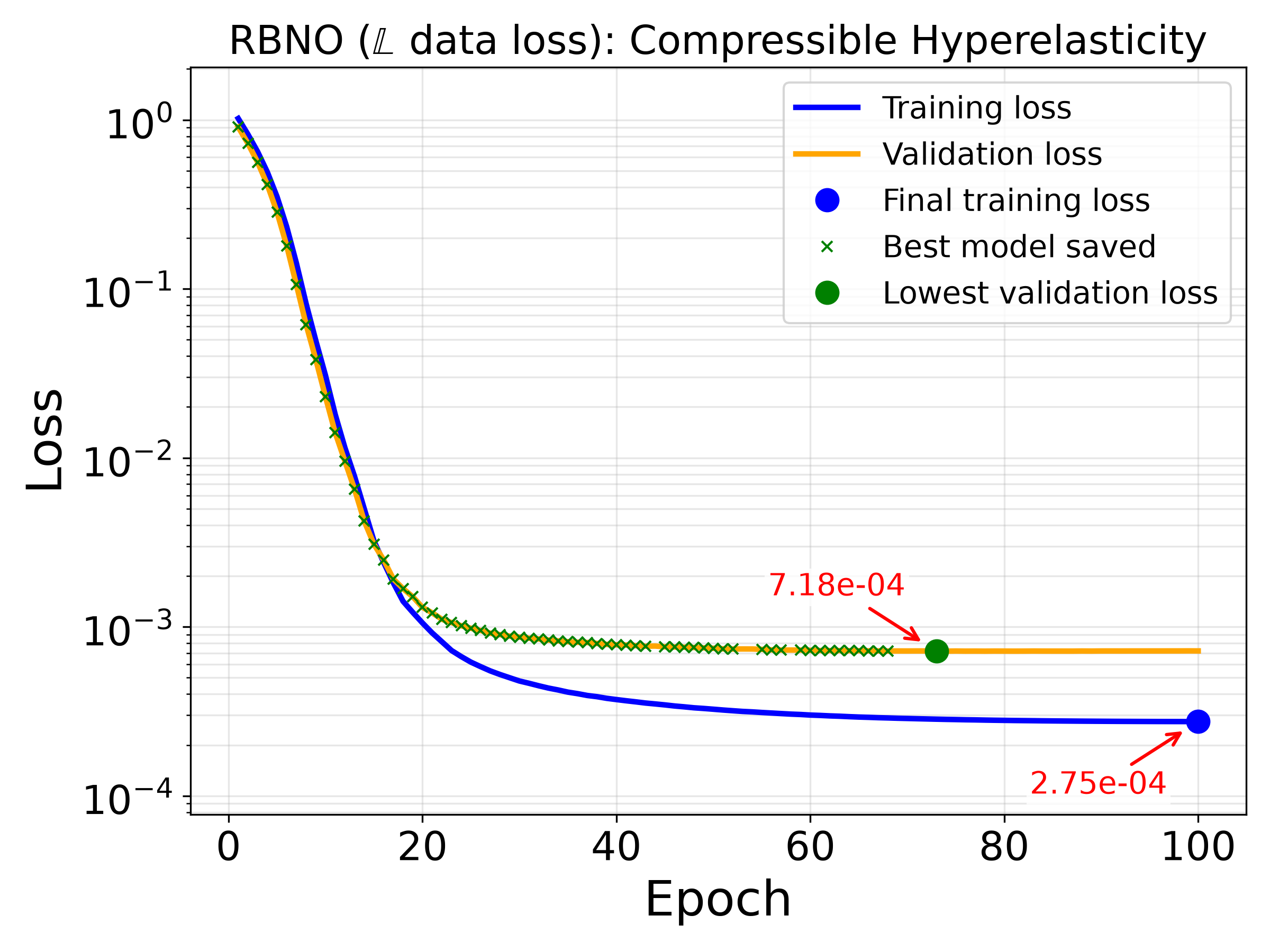}\\[2pt]
\includegraphics[width=.32\linewidth]{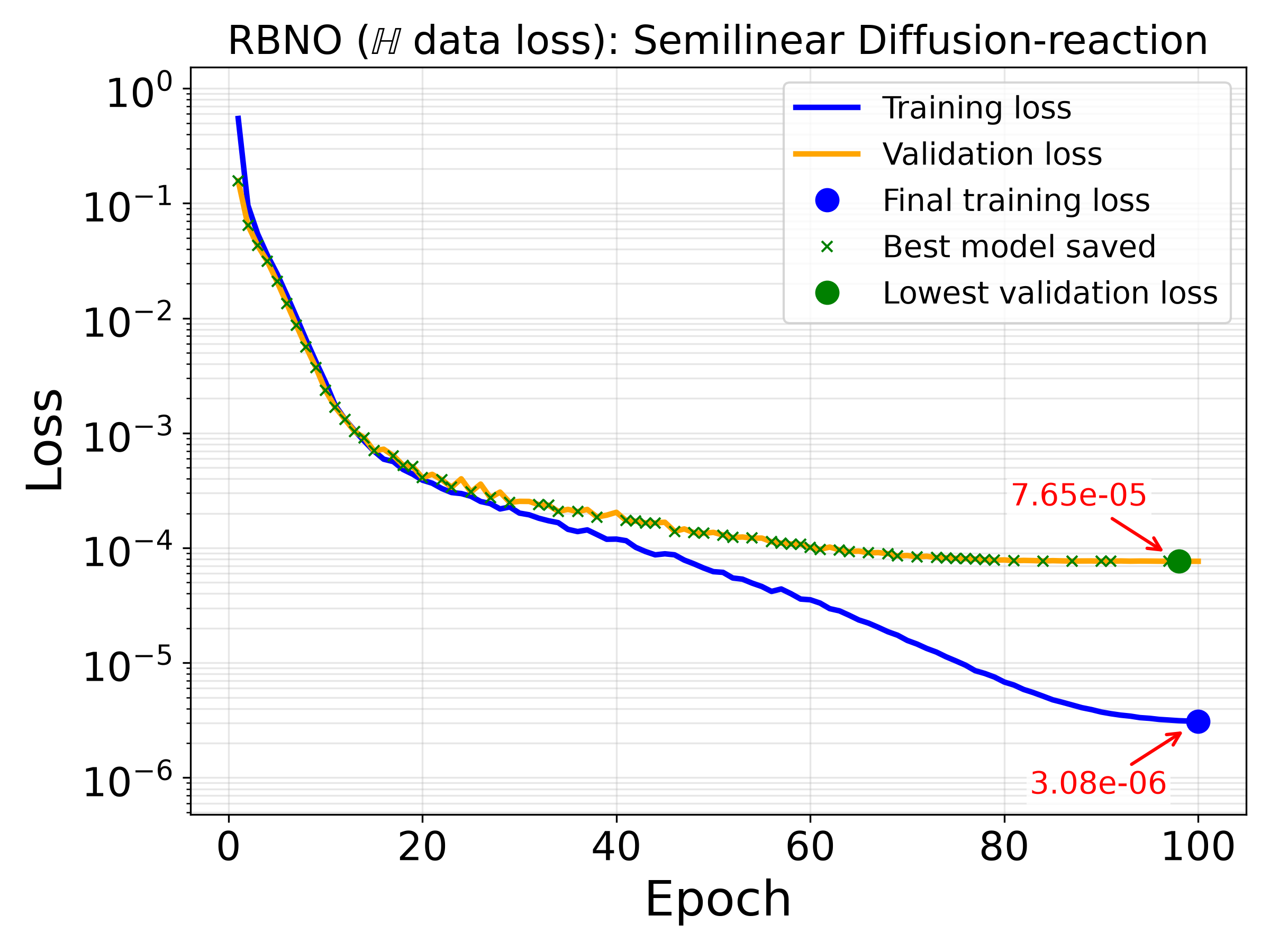}\hfill
\includegraphics[width=.32\linewidth]{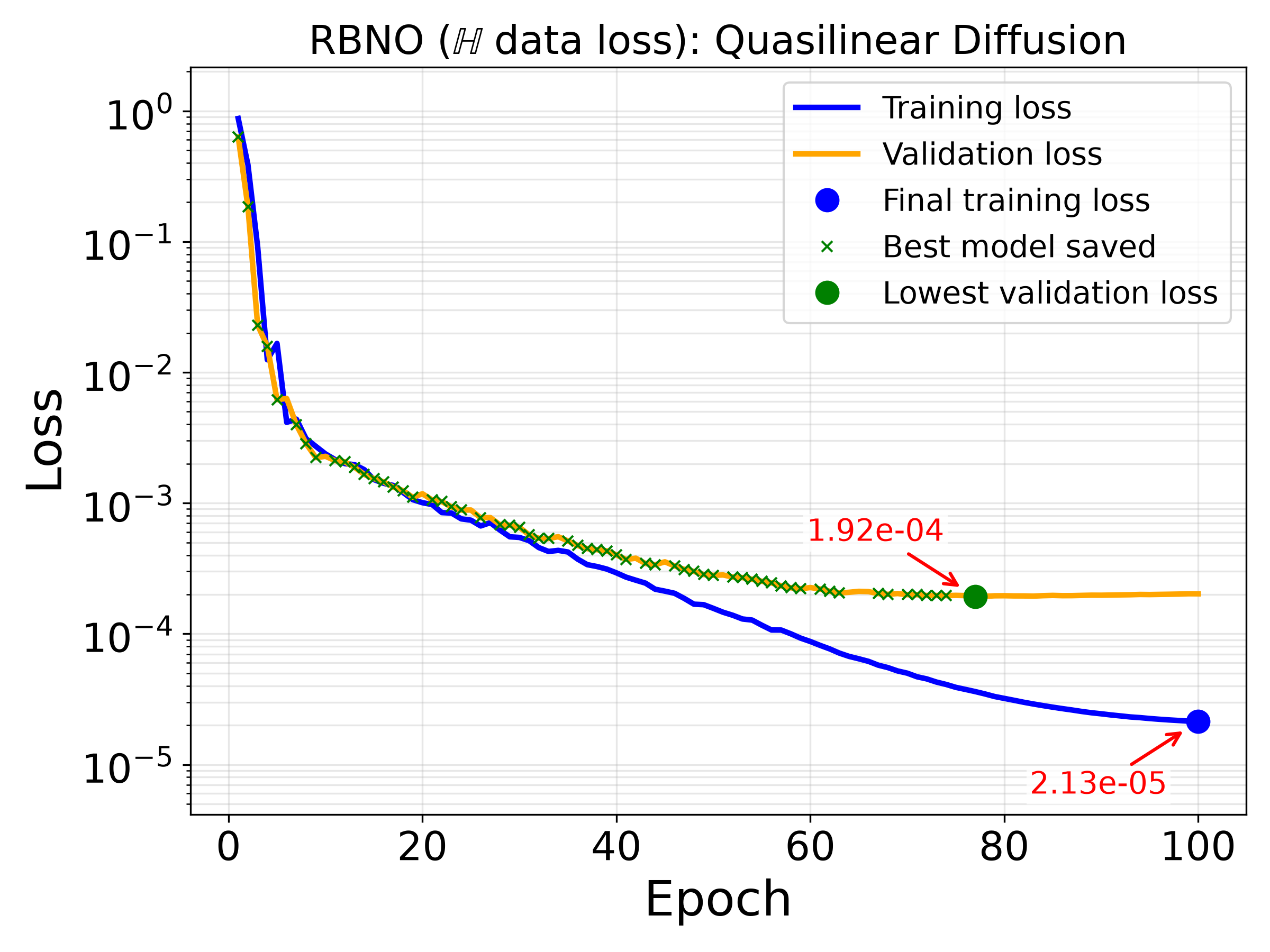}\hfill
\includegraphics[width=.32\linewidth]{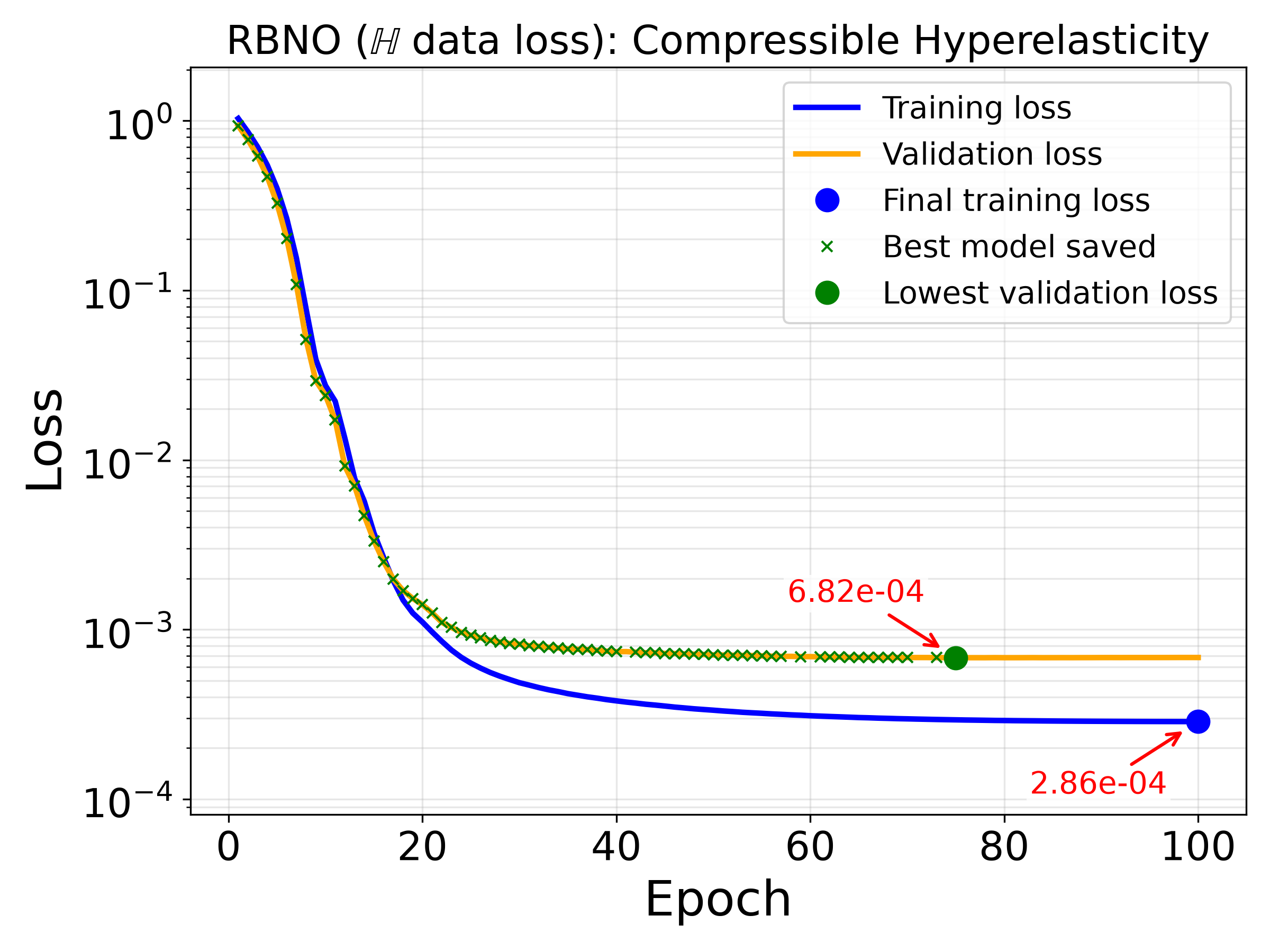}
\caption{Training histories of the two variants of RBNO trained on data in the $\LL^2$ norm (top) and in the $\HH$ norm (bottom); semilinear diffusion--reaction ($\kappa=10^3$), quasilinear diffusion, and compressible hyperelasticity benchmarks from left to right.}\label{fig:baseline_training_rbno}
\end{figure}

\begin{figure}[!htbp]
\centering
\includegraphics[width=.32\linewidth]{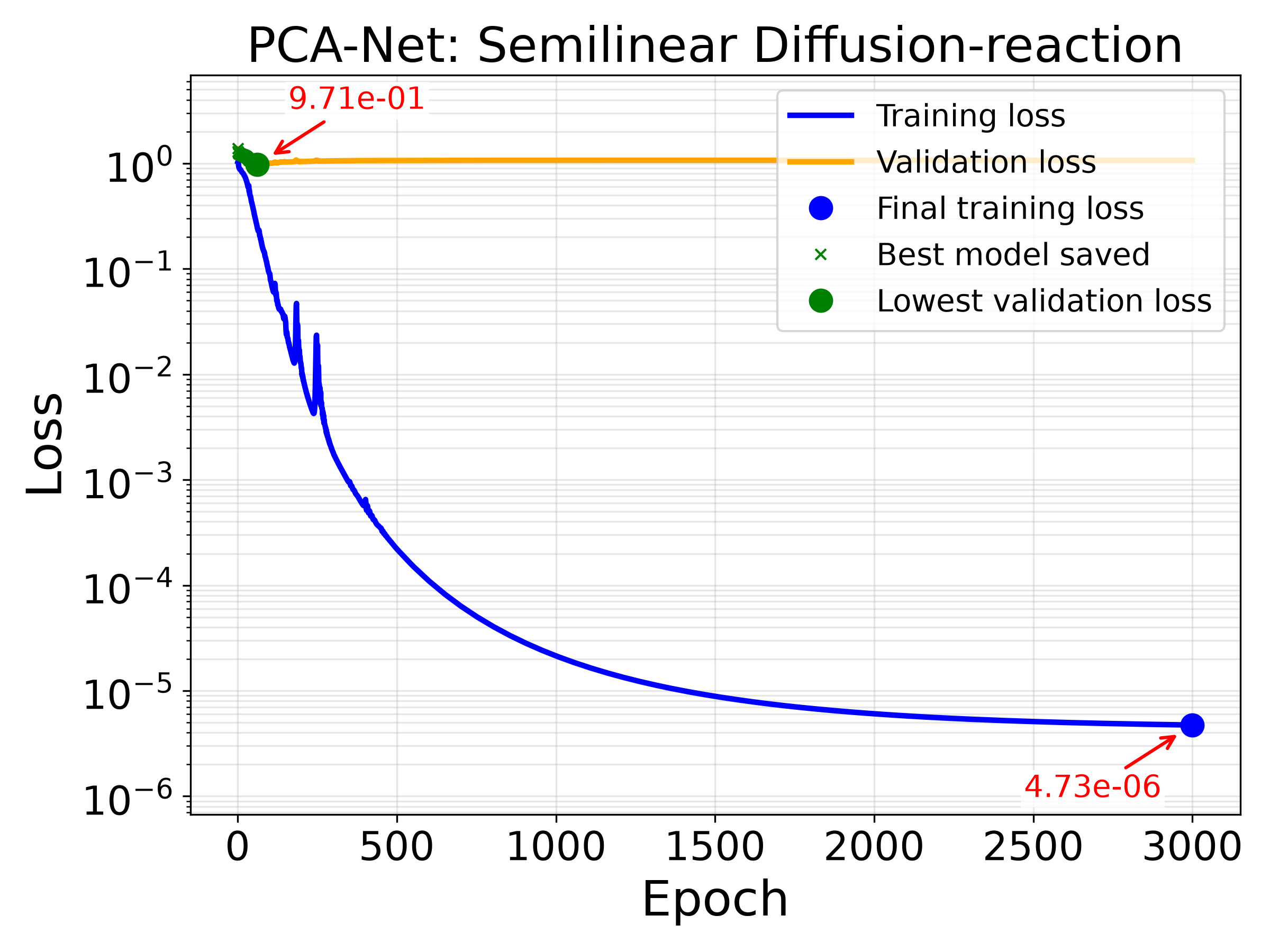}\hfill
\includegraphics[width=.32\linewidth]{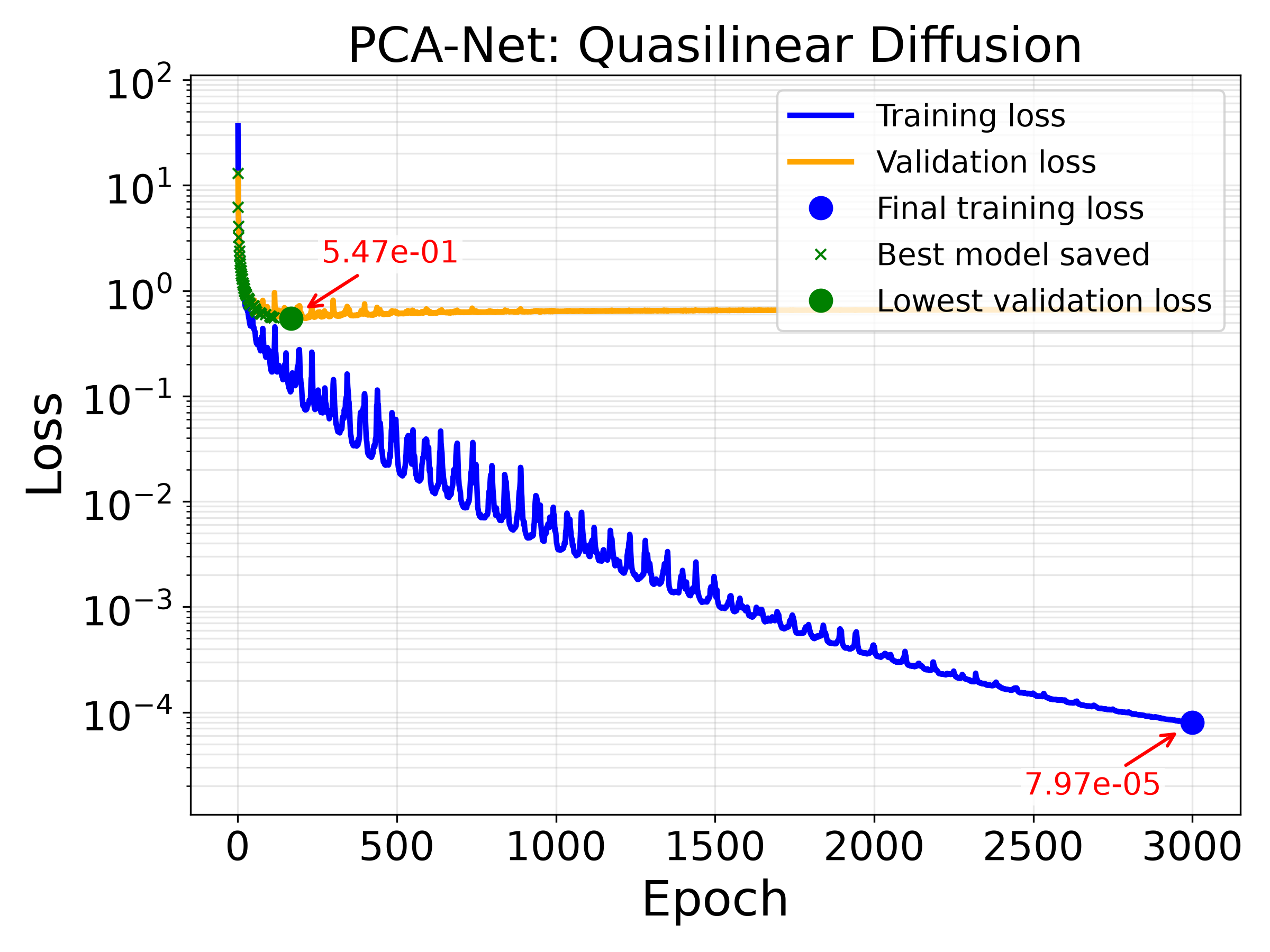}\hfill
\includegraphics[width=.32\linewidth]{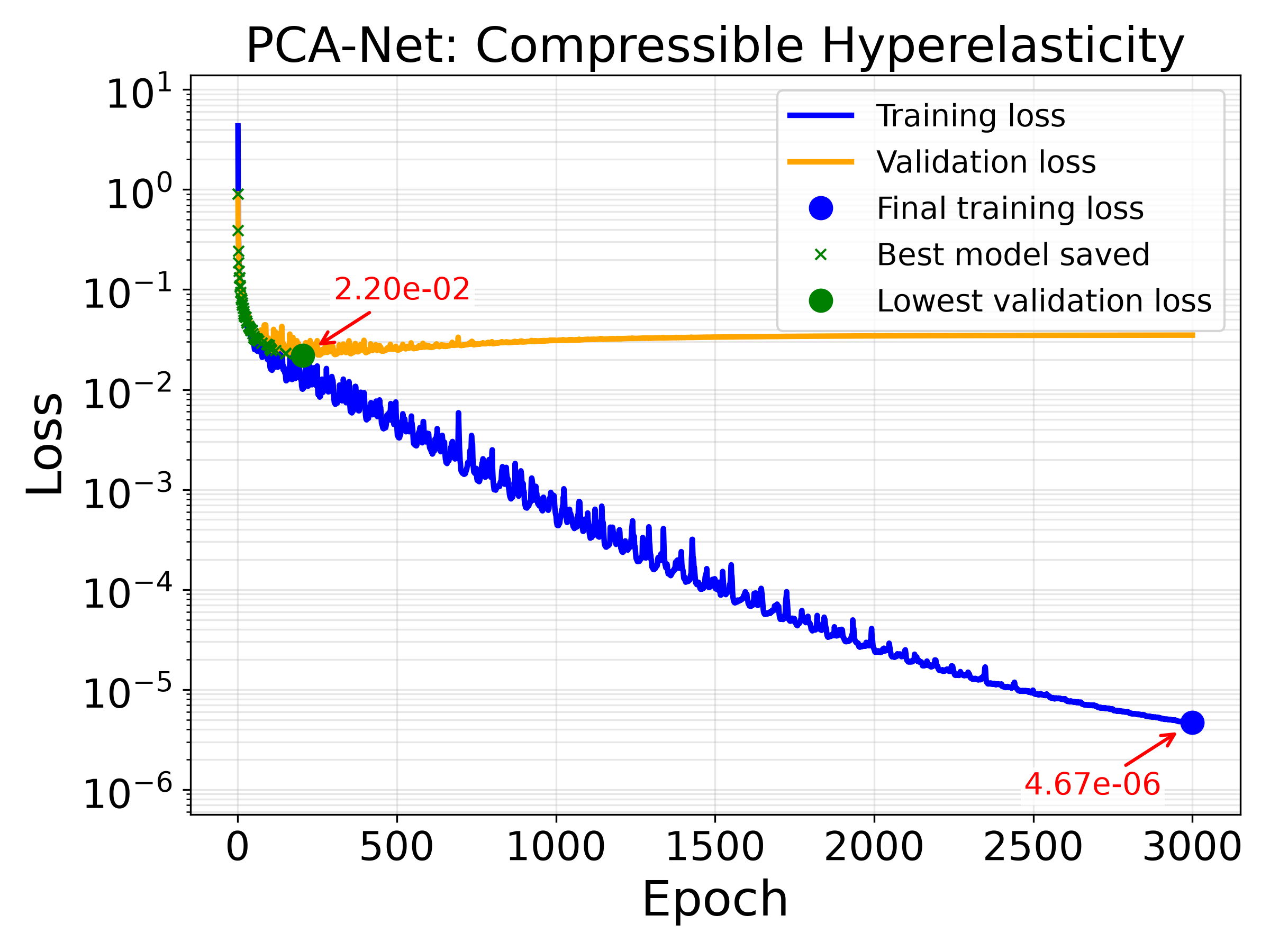}\\[2pt]
\includegraphics[width=.32\linewidth]{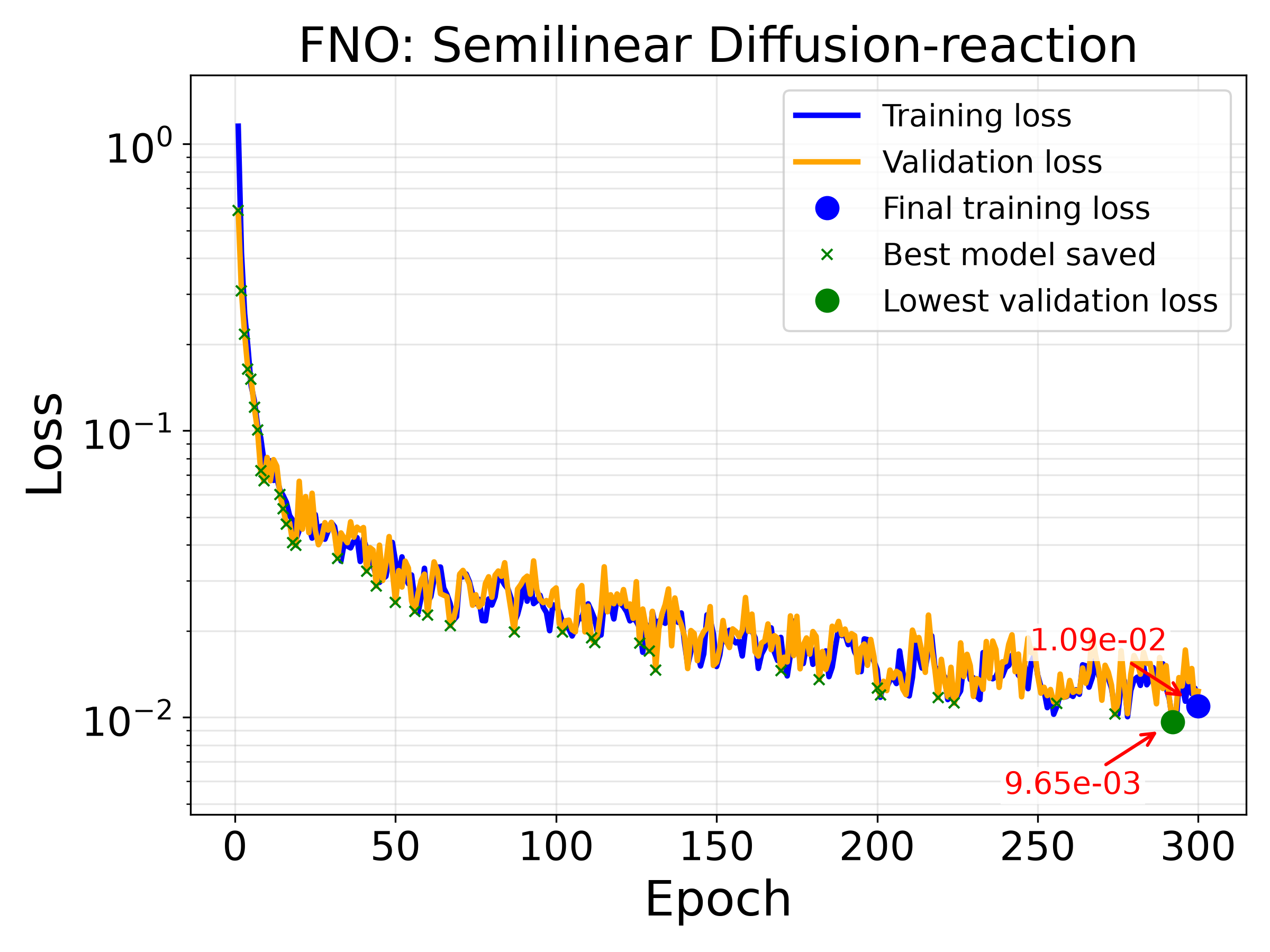}\hfill
\includegraphics[width=.32\linewidth]{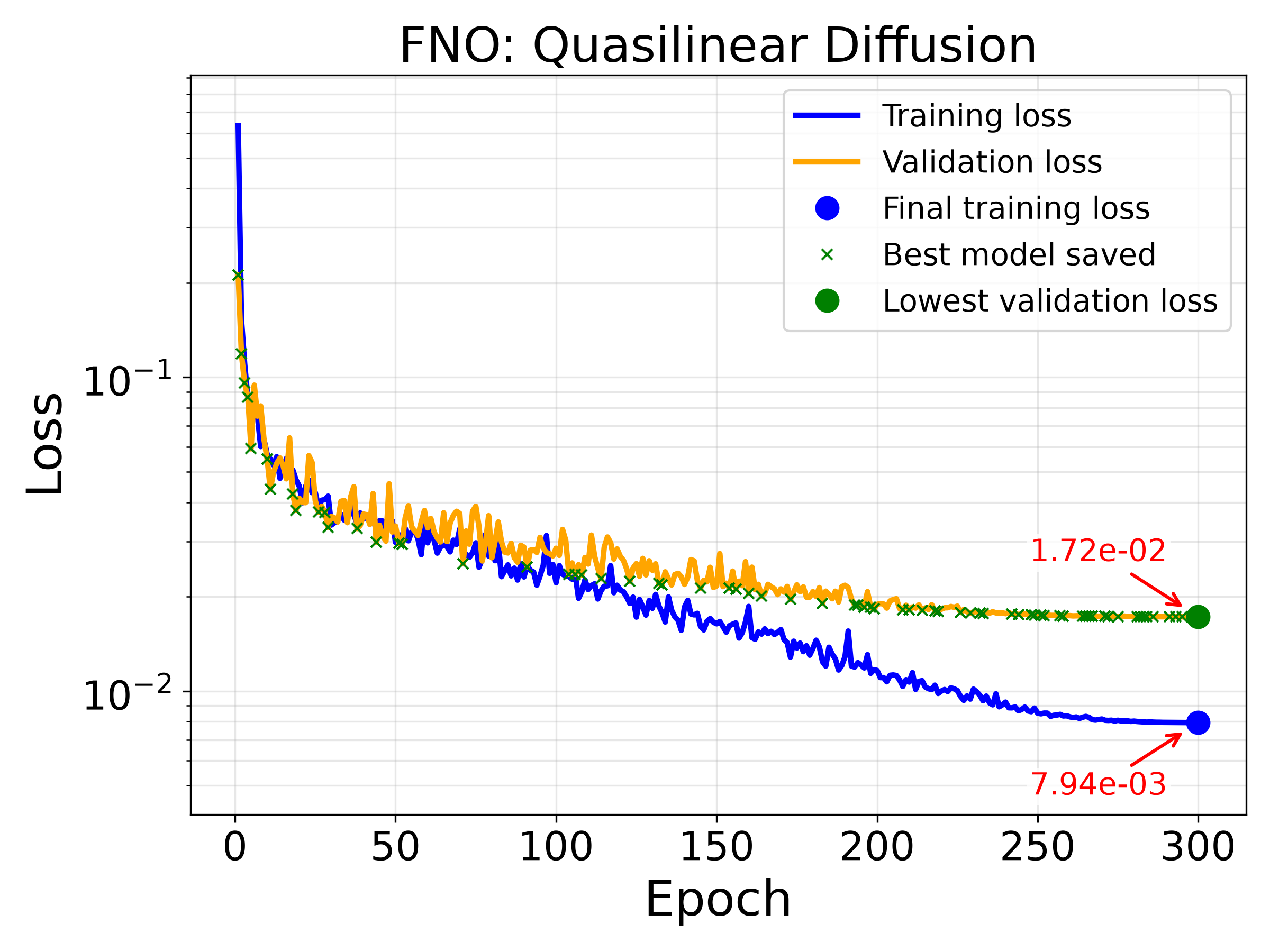}\hfill
\includegraphics[width=.32\linewidth]{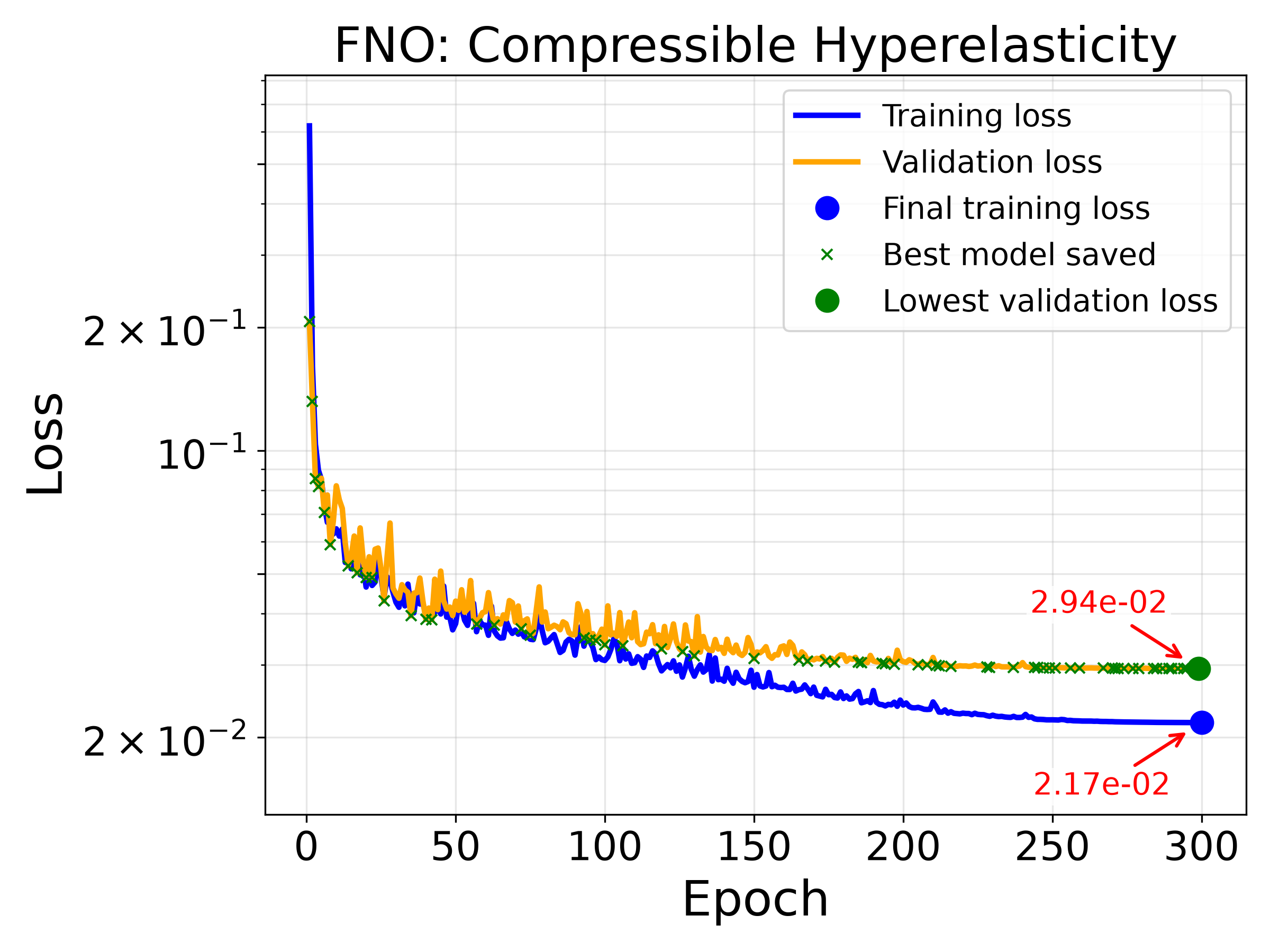}\\[4pt]
\includegraphics[width=\linewidth]{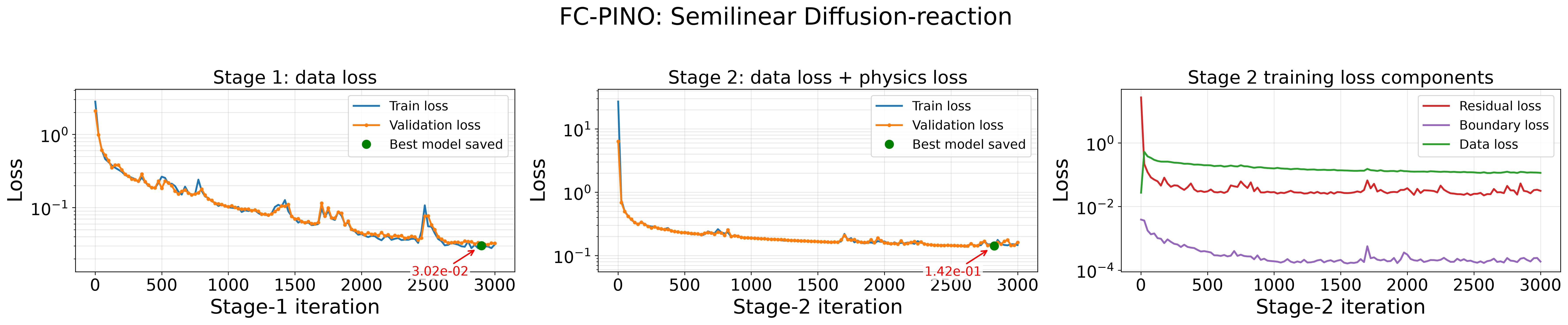}\\[2pt]
\includegraphics[width=\linewidth]{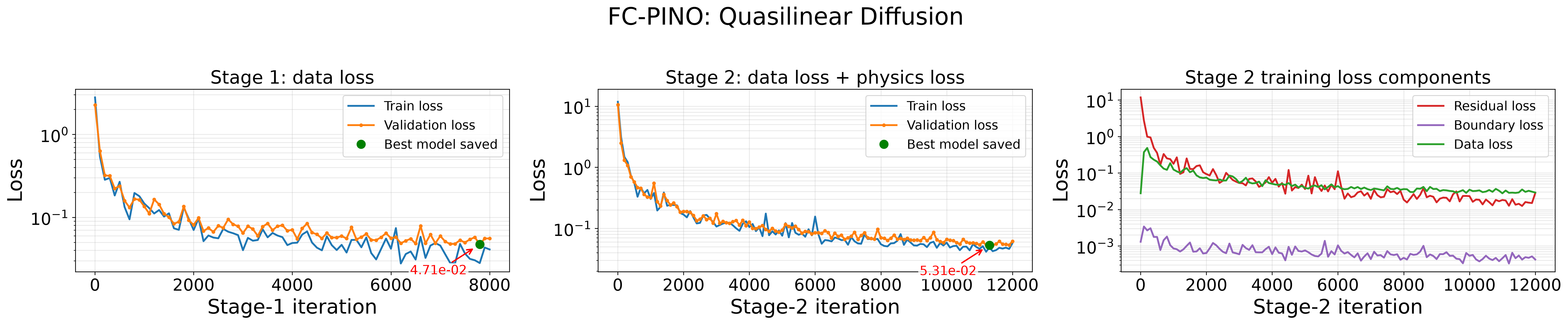}\\[2pt]
\includegraphics[width=\linewidth]{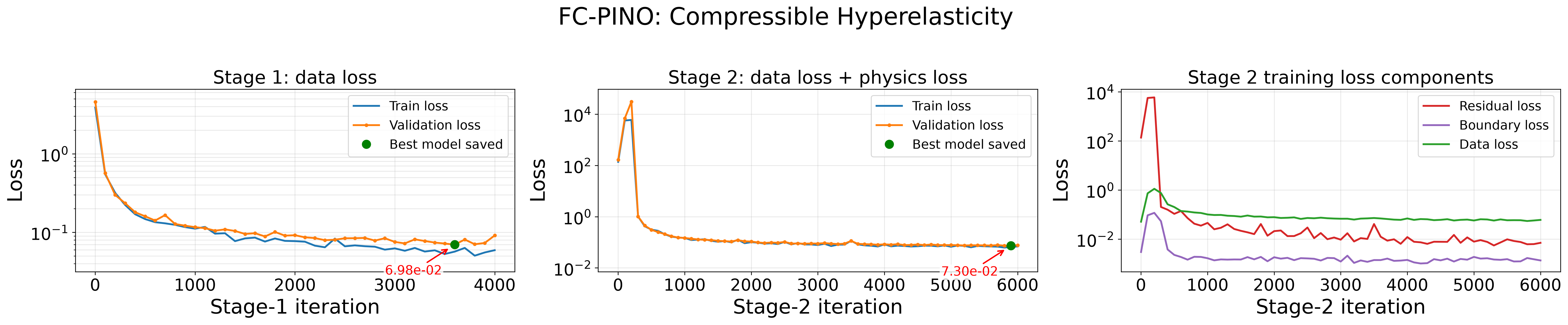}
\caption{Training histories of the grid-based baselines for the semilinear diffusion--reaction ($\kappa=10^3$), quasilinear diffusion, and compressible hyperelasticity benchmarks. First row: PCA-Net; second row: FNO; bottom rows: FC-PINO.}\label{fig:baseline_training}
\end{figure}

We summarize the architectures and training settings of PCA-Net, FNO, and FC-PINO in \cref{tab:architecture_training_details_pcanet,tab:architecture_training_details_fno,tab:architecture_training_details_fc_pino}. The training data for all compared models are described in \cref{subsec:exp_baselines}. In the second stage of FC-PINO training, we assign equal weights to the data loss, the lifted FOSLS residual loss evaluated on the grid, and the penalty for homogeneous boundary conditions. We present the training histories of the two RBNO variants trained with data losses and the grid-based neural operators in \cref{fig:baseline_training_rbno,fig:baseline_training}.

\begin{table}[!htbp]
\centering
\setlength{\tabcolsep}{4pt}
\renewcommand{\arraystretch}{0.95}
\begin{footnotesize}
\begin{tabular}{@{}lccc@{}}
\toprule
 & Semilinear diffusion--reaction & Quasilinear diffusion & Compressible hyperelasticity \\
\midrule
input/output dimension & $20/128$ & $512/512$ & $128/256$ \\
hidden dimensions & \multicolumn{3}{c}{$(1024,\,2048,\,1024)$} \\
activation function & \multicolumn{3}{c}{SELU} \\
PCA coefficients & standardized & \multicolumn{2}{c}{not standardized} \\
numerical precision & \multicolumn{3}{c}{float32 (PCA bases computed in float64)} \\
number of trainable parameters & $4.35$\,M & $5.25$\,M & $4.59$\,M \\
\midrule
\multicolumn{4}{@{}l}{\textit{Training}} \\
loss & \multicolumn{3}{c}{mean squared error of the output PCA coefficients} \\
\# training / validation samples & $500$ / $500$ & $2000$ / $200$ & $2000$ / $200$ \\
epochs & \multicolumn{3}{c}{$3000$} \\
optimizer & \multicolumn{3}{c}{AdamW~\cite{loshchilov2019decoupled} (lr $=10^{-3}$, weight decay $=10^{-2}$)} \\
lr scheduler & \multicolumn{3}{c}{StepLR (step $=50$ epochs, gamma $=0.9$)} \\
batch size & $500$ (full batch) & $100$ & $100$ \\
\bottomrule
\end{tabular}
\end{footnotesize}
\caption{Architecture and training settings of PCA-Net.}
\label{tab:architecture_training_details_pcanet}
\end{table}

\begin{table}[!htbp]
\centering
\setlength{\tabcolsep}{4pt}
\renewcommand{\arraystretch}{0.95}
\begin{footnotesize}
\begin{tabular}{@{}lccc@{}}
\toprule
 & Semilinear diffusion--reaction & Quasilinear diffusion & Compressible hyperelasticity \\
\midrule
number of modes & \multicolumn{3}{c}{$(32,\,32)$} \\
in/out channels ($\pp$, $x_1$, $x_2$) & $3/3$ & $3/3$ & $3/6$ \\
hidden channels & \multicolumn{3}{c}{$128$} \\
FNO layers & \multicolumn{3}{c}{$4$} \\
activation & \multicolumn{3}{c}{GELU} \\
numerical precision & \multicolumn{3}{c}{float32} \\
\# trainable parameters & \multicolumn{3}{c}{$35.9$\,M} \\
\midrule
\multicolumn{4}{@{}l}{\textit{Training}} \\
loss & \multicolumn{3}{c}{relative $\ell^2$ error} \\
\# training / validation samples & $500$ / $50$ & $2000$ / $200$ & $2000$ / $200$ \\
epochs & \multicolumn{3}{c}{$300$} \\
optimizer & \multicolumn{3}{c}{AdamW (lr $=10^{-3}$, weight decay $=10^{-2}$)} \\
lr scheduler & StepLR & \multicolumn{2}{c}{CosineAnnealingLR} \\
 & (step $=50$ epochs, gamma $=0.9$) & \multicolumn{2}{c}{($T_{\max}=300$, $\eta_{\min}=10^{-6}$)} \\
batch size & $10$ & $20$ & $20$ \\
\bottomrule
\end{tabular}
\end{footnotesize}
\caption{Architecture and training settings of FNO.}
\label{tab:architecture_training_details_fno}
\end{table}

\begin{table}[!htbp]
\centering
\setlength{\tabcolsep}{4pt}
\renewcommand{\arraystretch}{0.95}
\begin{footnotesize}
\begin{tabular}{@{}lccc@{}}
\toprule
& Semilinear diffusion--reaction & Quasilinear diffusion & Compressible hyperelasticity \\
\midrule
\multicolumn{4}{@{}l}{\textit{Architecture}} \\
number of modes & \multicolumn{3}{c}{$(32,\,32)$} \\
in/out channels ($\pp$, $x_1$, $x_2$) & $3/3$ & $3/3$ & $3/6$ \\
hidden channels & \multicolumn{3}{c}{$128$} \\
FNO layers & \multicolumn{3}{c}{$4$} \\
FNO activation & \multicolumn{3}{c}{GELU} \\
projection activation & \multicolumn{3}{c}{$\tanh$} \\
Fourier continuation & \multicolumn{3}{c}{Legendre} \\
continuation polynomial degree & $6$ & $3$ & $3$ \\
continuation points & \multicolumn{3}{c}{$50$} \\
numerical precision & \multicolumn{3}{c}{float64} \\
\# trainable parameters & \multicolumn{3}{c}{$35.9$\,M} \\
\midrule
\multicolumn{4}{@{}l}{\textit{Stage 1 (data loss: relative $\ell^2$ error)}} \\
\# parameter--solution pairs & \multicolumn{3}{c}{$1000$} \\
\# validation pairs & $500$ & $100$ & $100$ \\
iterations & $3000$ & $8000$ & $4000$ \\
\midrule
\multicolumn{4}{@{}l}{\textit{Stage 2 (data loss $+$ residual loss $+$ boundary loss, weights $=(1,1,1)$)}} \\
\# parameter instances & $3000$ & $8000$ & $5000$ \\
\# validation pairs & $500$ & $100$ & $100$ \\
iterations & $3000$ & $12000$ & $6000$ \\
\midrule
\multicolumn{4}{@{}l}{\textit{Stages 1 and 2}} \\
optimizer & \multicolumn{3}{c}{AdamW (lr $=10^{-3}$, weight decay $=10^{-2}$)} \\
lr scheduler & \multicolumn{3}{c}{StepLR (step $=600$, $2000$, $1000$ iterations; gamma $=0.8$)} \\
batch size & \multicolumn{3}{c}{$10$} \\
gradient clipping & \multicolumn{3}{c}{maximal $\ell^2$ norm $=1$} \\
\bottomrule
\end{tabular}
\end{footnotesize}
\caption{Architecture and training settings of FC-PINO.}
\label{tab:architecture_training_details_fc_pino}
\end{table}

\FloatBarrier
\subsection{Pointwise residuals of FC-PINO and RBNO}\label{app:residual_maps}
\begin{figure}[!htbp]
\centering
\begin{subfigure}[b]{0.49\linewidth}\centering
\includegraphics[width=\linewidth]{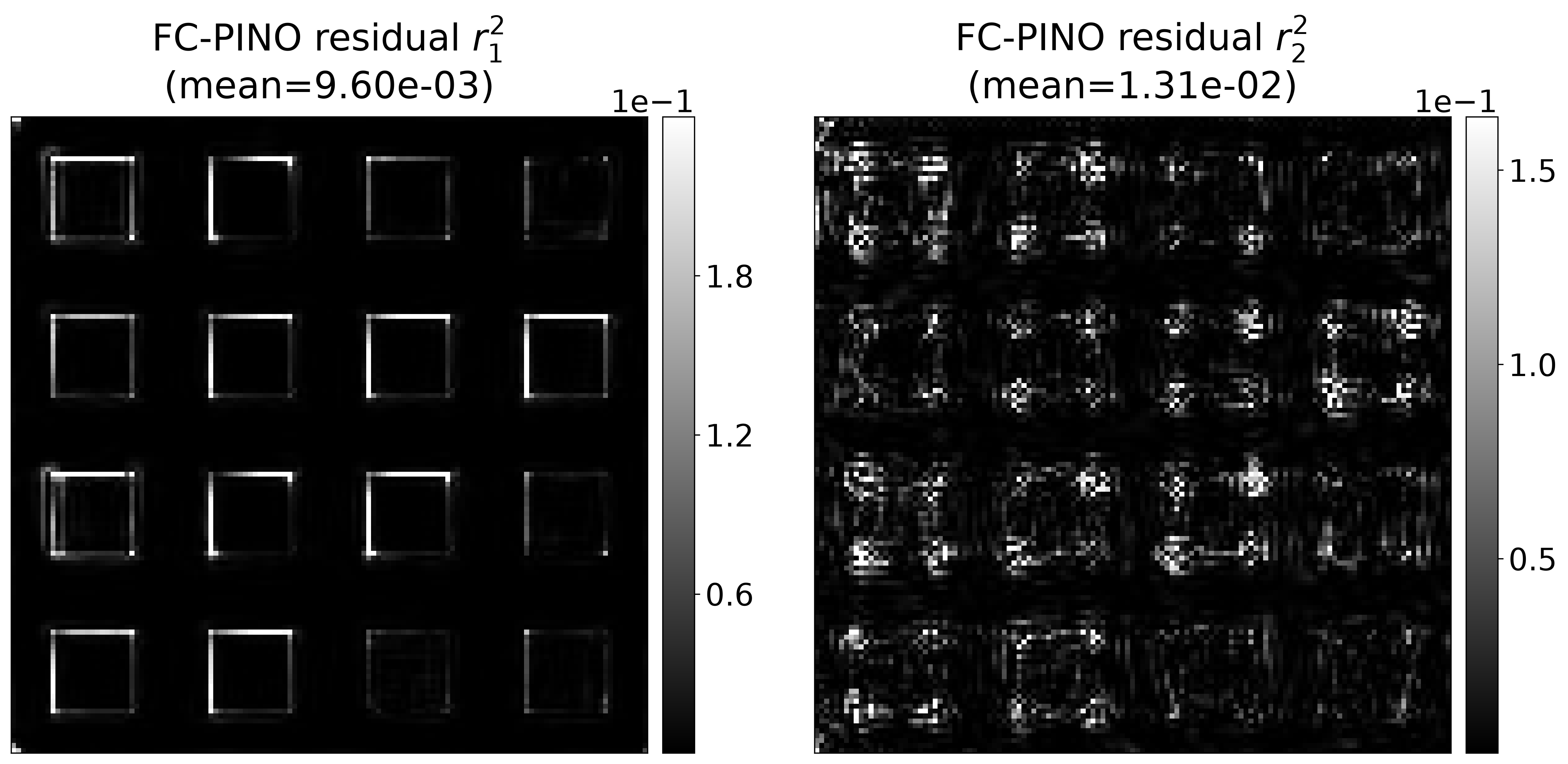}
\caption{Semilinear diffusion--reaction (FC-PINO)}
\end{subfigure}
\hfill
\begin{subfigure}[b]{0.49\linewidth}\centering
\includegraphics[width=\linewidth]{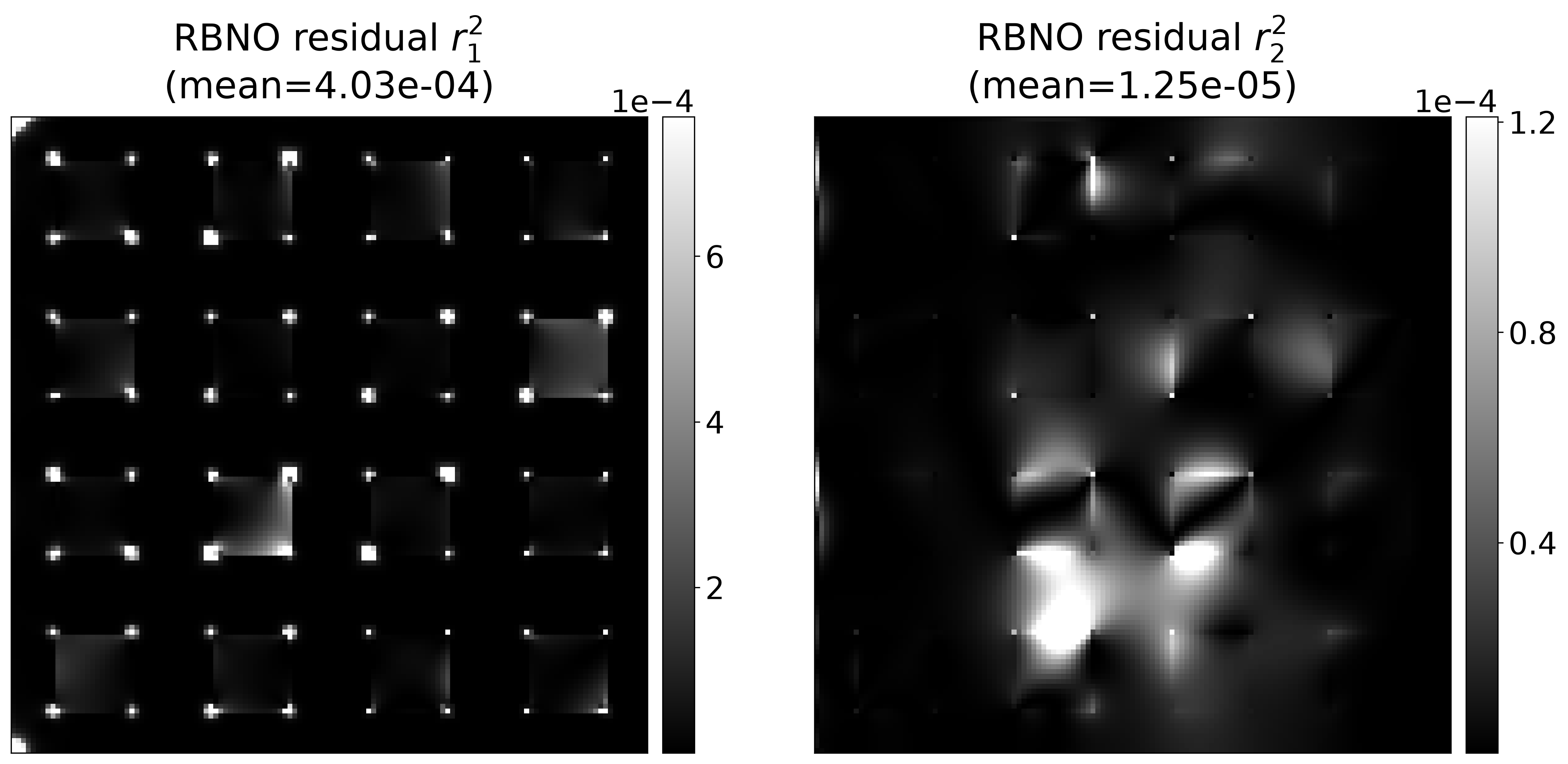}
\caption{Semilinear diffusion--reaction (RBNO)}
\end{subfigure}
\\[6pt]
\begin{subfigure}[b]{0.49\linewidth}\centering
\includegraphics[width=\linewidth]{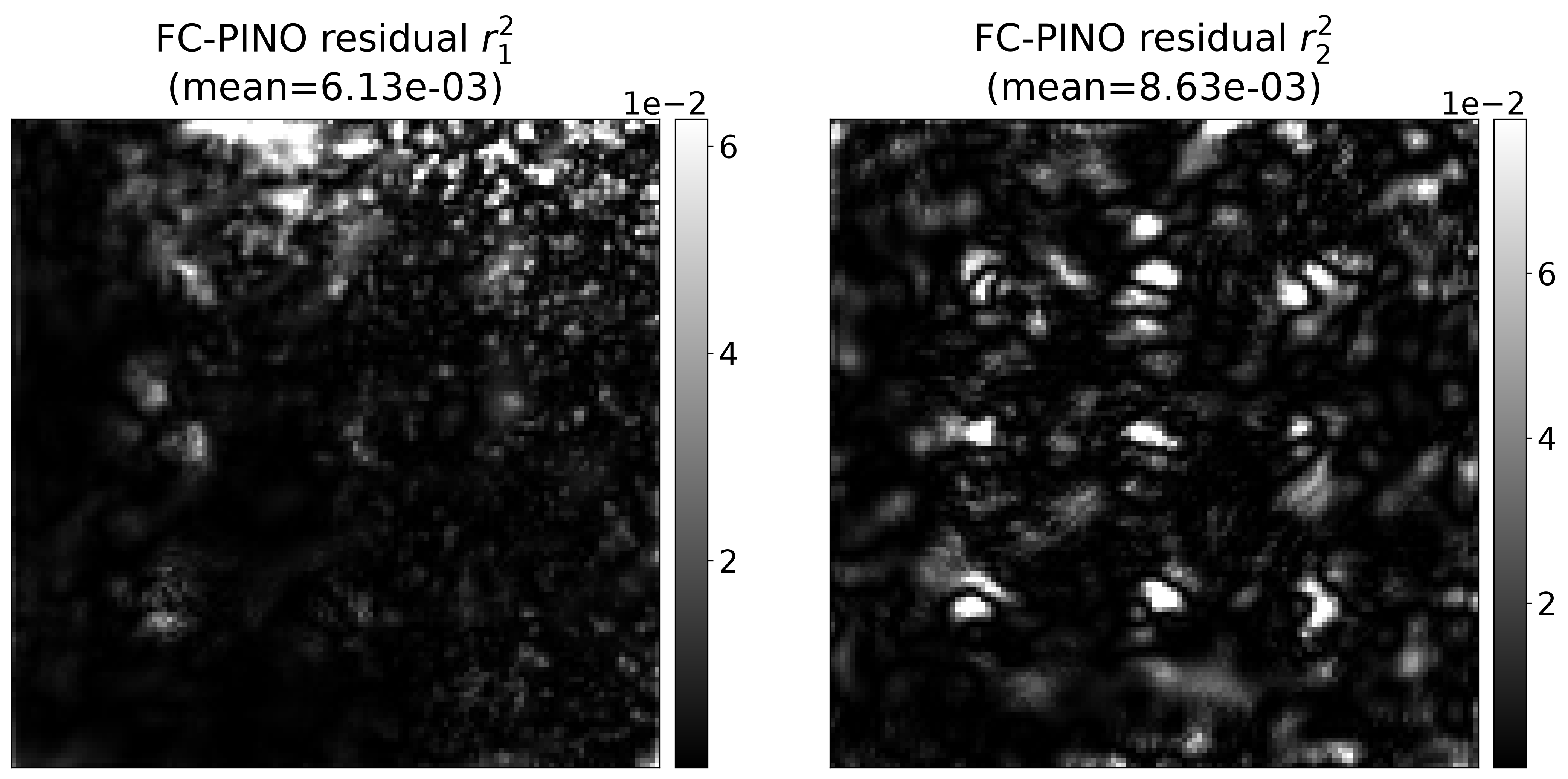}
\caption{Quasilinear diffusion (FC-PINO)}
\end{subfigure}
\hfill
\begin{subfigure}[b]{0.49\linewidth}\centering
\includegraphics[width=\linewidth]{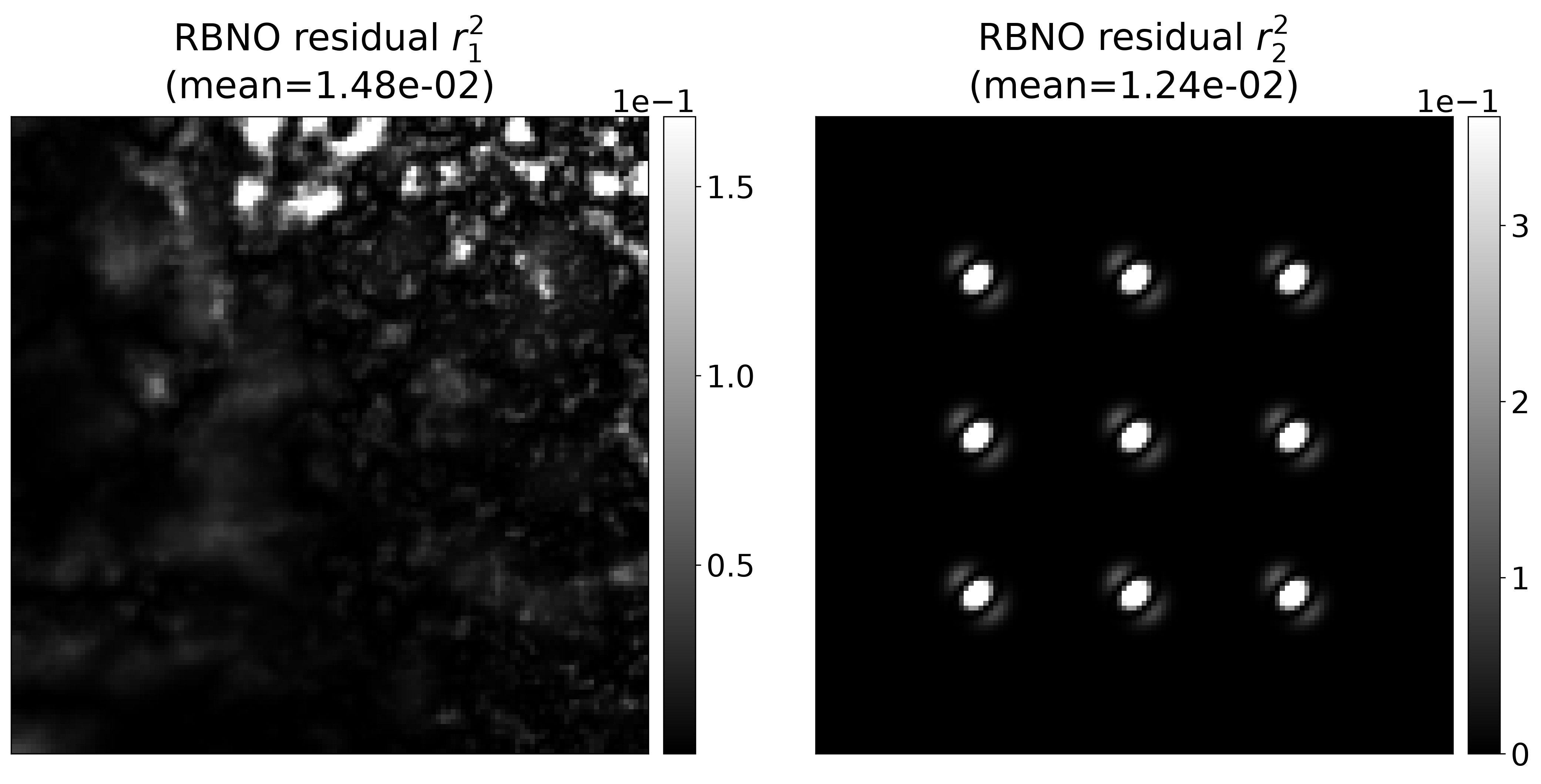}
\caption{Quasilinear diffusion (RBNO)}
\end{subfigure}
\\[6pt]
\begin{subfigure}[b]{0.4\linewidth}\centering
\includegraphics[width=\linewidth]{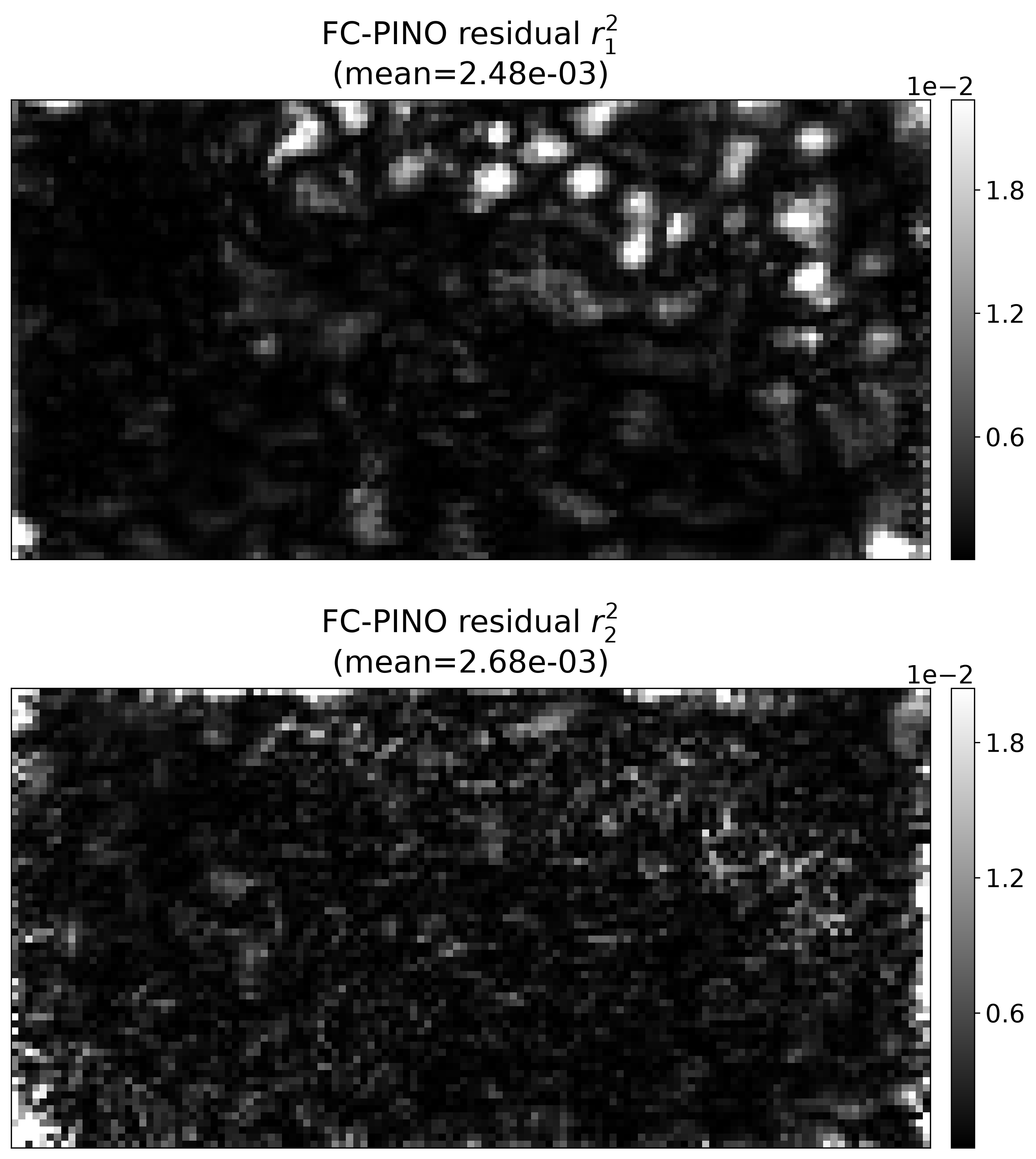}
\caption{Compressible hyperelasticity (FC-PINO)}
\end{subfigure}
\hspace{0.06\linewidth}
\begin{subfigure}[b]{0.4\linewidth}\centering
\includegraphics[width=\linewidth]{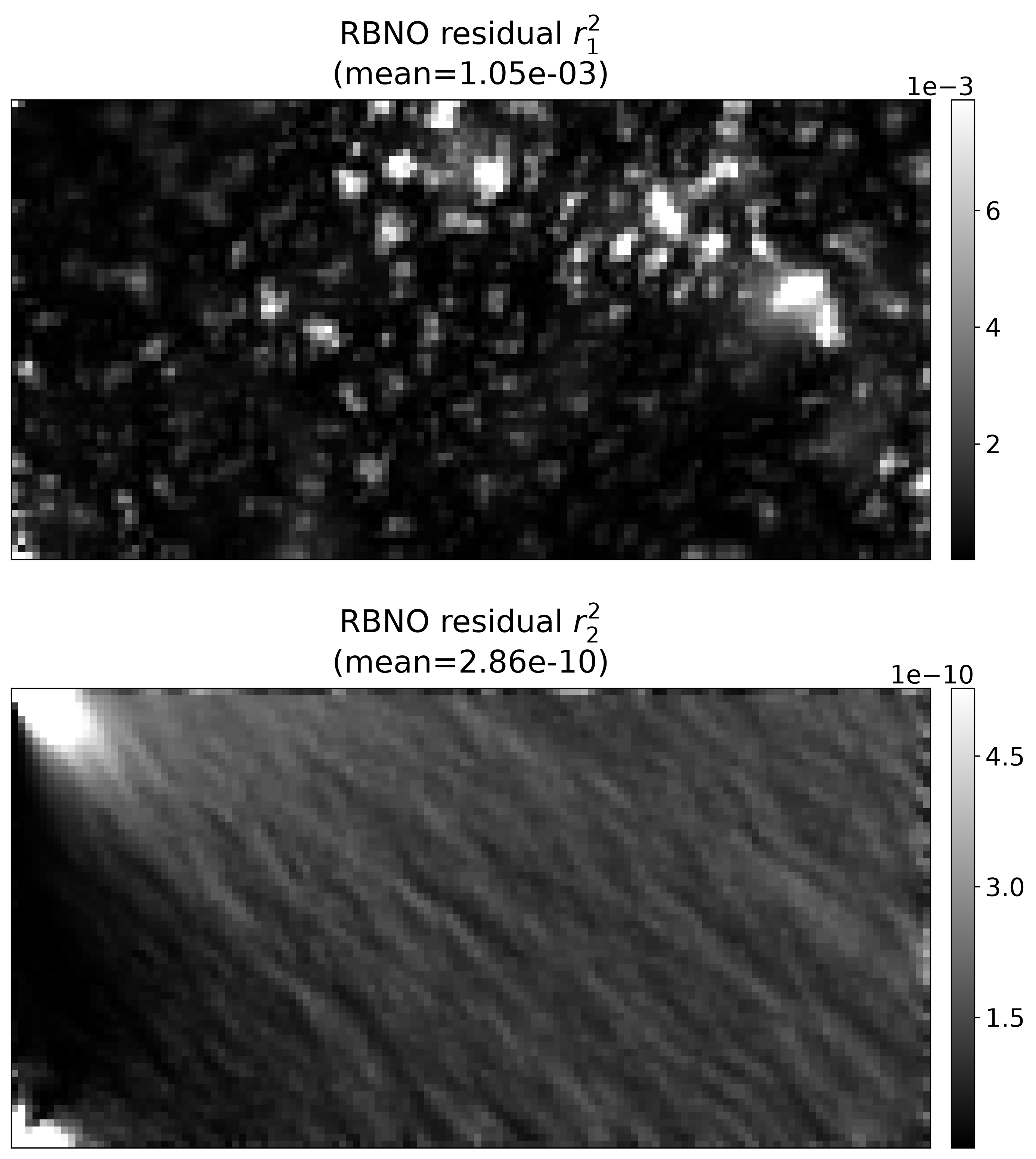}
\caption{Compressible hyperelasticity (RBNO)}
\end{subfigure}
\caption{Pointwise squared residual components $r_1^2$ and $r_2^2$ for FC-PINO (left column) and RBNO (right column), using the test parameters shown in \cref{fig:benchmarks}. Panels (a, b) show semilinear diffusion--reaction with $\kappa=10^3$, and panels (c, d) show quasilinear diffusion, with $r_1^2$ on the left and $r_2^2$ on the right within each panel. Panels (e, f) show compressible hyperelasticity on the $2\times1$ domain, with $r_1^2$ above $r_2^2$. For FC-PINO, we evaluate the residual at the grid points using spectral differentiation with Fourier continuation, as in training. For RBNO, we evaluate the finite element residual at the same grid points. We set the color limits separately for each image at the 1st and 99th percentiles and report the mean of the plotted values in each title. For RBNO, these grid averages differ from the integrated squared residuals $\norm{r_i}^2$.}\label{fig:residual_maps}
\end{figure}

\Cref{fig:residual_maps} shows the two components $r_1^2$ and $r_2^2$ of the pointwise FOSLS residual of FC-PINO and of RBNO for the test parameters of \cref{fig:benchmarks}.

\FloatBarrier

\end{document}